\documentclass{amsart}

\usepackage{graphicx}
\usepackage{amsmath}
\usepackage{amscd}
\usepackage{amsthm}
\usepackage{amsfonts}
\usepackage{amssymb}
\usepackage{hyperref}
\usepackage[all]{xy}

\newcommand{\chap}[1]{\section{#1}}
\newcommand{\sect}[1]{\subsection{#1}}
\newcommand{\firstpage}{}
\newcommand{\chapword}{section}
\newcommand{\Chapword}{Section}

\newtheorem{thm}{Theorem}[section]
\newtheorem*{nonum}{Theorem}

\newtheorem{claim}[thm]{Claim}

\newtheorem{corollary}[thm]{Corollary}

\newtheorem*{thmMainTrue}{Theorem \ref{thm.main}}
\newtheorem*{thmMainGen}{Theorem \ref{thm.q.main}}
\newtheorem*{thmMainLag}{Theorem \ref{t.Lag.Perm}}
\newtheorem*{corClassify}{Corollary \ref{cor.classify}}

\theoremstyle{definition}
\newtheorem{exam}[thm]{Example}
\newtheorem{defn}[thm]{Definition}
\theoremstyle{plain}

\newtheorem{lem}[thm]{Lemma}

\newtheorem{prop}[thm]{Proposition}
\newtheorem{rem}[thm]{Remark}

\newcommand{\term}[1]{\textbf{#1}}
\newcommand{\mtrx}[1]{\begin{matrix} #1 \end{matrix}}
\newcommand{\pmtrx}[1]{\begin{pmatrix} #1 \end{pmatrix}}
\newcommand{\cmtrx}[1]{\left\{\begin{matrix}#1\end{matrix}\right\}}
\newcommand{\arry}[2]{\begin{array}{#1} #2 \end{array}}
\newcommand{\tbl}[2]{\begin{tabular}{#1} #2 \end{tabular}}
\newcommand{\RHScase}[1]{\left\{\begin{array}{ll} #1 \end{array}\right.}

\newcommand{\CC}{\mathbb{C}}
\newcommand{\NN}{\mathbb{N}}
\newcommand{\RR}{\mathbb{R}}

\newcommand{\ZZ}{\mathbb{Z}}

\newcommand{\AAA}{\mathcal{A}}
\newcommand{\BBB}{\mathcal{B}}
\newcommand{\CCC}{\mathcal{C}}
\newcommand{\HHH}{\mathcal{H}}
\newcommand{\LLL}{\mathcal{L}}
\newcommand{\OOO}{\mathcal{O}}
\newcommand{\QQQ}{\mathcal{Q}}
\newcommand{\RRR}{\mathcal{R}}
\newcommand{\TTT}{\mathcal{T}}

\newcommand{\eps}{\varepsilon}

\newcommand{\dset}[1]{\{1,\dots, #1\}}
\newcommand{\ddset}[2]{\{#1,\dots, #2\}}
\newcommand{\inv}[1]{#1^{-1}}

\renewcommand{\AA}{\mathbb{A}}
\newcommand{\sS}{\mathbb{S}}

\newcommand{\perm}{\mathfrak{S}}
\newcommand{\genperm}{\mathfrak{Q}}
\newcommand{\irr}{\perm^0}
\newcommand{\genirr}{\genperm^0}
\newcommand{\RClass}{\RRR}

\newcommand{\Heven}{\HHH^{even}}
\newcommand{\Hodd}{\HHH^{odd}}
\newcommand{\Hhyp}{\HHH^{hyp}}
\newcommand{\Hnonhyp}{\HHH^{nonhyp}}

\newcommand{\tT}{\tilde T}
\newcommand{\ti}{i^\pi}
\newcommand{\td}{d^\pi}
\newcommand{\tdd}{(2d)^\pi}

\newcommand{\hI}{\hat I}

\newcommand{\mA}{\mathbf{A}}
\newcommand{\mB}{\mathbf{B}}
\newcommand{\mD}{\mathbf{D}}

\newcommand{\mP}{\mathbf{P}}
\newcommand{\mQ}{\mathbf{Q}}
\newcommand{\mS}{\mathbf{S}}
\newcommand{\mT}{\mathbf{T}}
\newcommand{\mU}{\mathbf{U}}
\newcommand{\mV}{\mathbf{V}}
\newcommand{\mW}{\mathbf{W}}
\newcommand{\mX}{\mathbf{X}}
\newcommand{\mY}{\mathbf{Y}}

\newcommand{\mO}{\mathbf{0}}

\newcommand{\ee}{\mathbf{e}}
\newcommand{\vv}{\mathbf{v}}
\newcommand{\ww}{\mathbf{w}}

\newcommand{\mone}{\mathbf{1}}

\newcommand{\msp}{\mbox{ }}

\newcommand{\IET}{IET}
\newcommand{\veech}{one-row}
\newcommand{\yoccoz}{two-row}

\newcommand{\rtt}{1}
\newcommand{\rtb}{0}
\newcommand{\rt}{r_{\rtt}}
\newcommand{\rb}{r_{\rtb}}
\newcommand{\reps}{r_\eps}
\newcommand{\ropp}{r_{1-\eps}}
\newcommand{\rsub}[1]{r_{#1}}
\newcommand{\sig}{\sigma}

\newcommand{\inx}{\cdot}

\newcommand{\EVEN}{\mU}
\newcommand{\ODD}{\mV}
\newcommand{\PAIR}{\mW}
\newcommand{\SPACE}{\mS}
\newcommand{\BLAH}{\mD}

\newcommand{\POLES}{\mP}
\newcommand{\INSERT}{\mQ}
\newcommand{\TWOFOURS}{\mT}
\newcommand{\JUSTFOURS}{\mX}
\newcommand{\TWOTWOS}{\mY}

\newcommand{\lA}{A}
\newcommand{\lB}{Z}
\newcommand{\LL}[2]{\mtrx{#1 \\ #2}}

\newcommand{\LLcdots}{\LL{\cdots}{\cdots}}
\newcommand{\AB}{\LL{\lA}{\lB}}
\newcommand{\BA}{\LL{\lB}{\lA}}
\newcommand{\stdperm}[1]{\cmtrx{\AB & #1 & \BA}}

\newcommand{\spaN}{\mathrm{span}}

\newcommand{\dsum}[1]{\underset{#1}{\sum}}
\newcommand{\dusum}[2]{\overset{#2}{\underset{#1}{\sum}}}

\begin{document}


\title[Self-Inverses in Rauzy Classes]{Self-inverses in Rauzy Classes}

\author[J. Fickenscher]{Jonathan Fickenscher}

\address{Department of Mathematics, Rice University, Houston, TX~77005, USA}

\email{jonfick@rice.edu}


\date{\today}

\maketitle

\pagestyle{headings}

\begin{abstract}
\noindent Thanks to works by M. Kontsevich and A. Zorich followed by C. Boissy, we have a classification of all Rauzy Classes of any given genus. It follows from these works that Rauzy Classes are closed under the operation of inverting the permutation. In this paper, we shall prove the existence of self-inverse permutations in every Rauzy Class by giving an explicit construction of such an element satisfying the sufficient conditions. As a corollary, we will give another proof that every Rauzy Class is closed under taking inverses. In the case of generalized permutations, generalized Rauzy Classes have been classified by works of M. Kontsevich, H. Masur and J. Smillie, E. Lanneau, and again C. Boissy. We state the definition of self-inverse for generalized permutations and prove a necessary and sufficient condition for a generalized Rauzy Class to contain self-inverse elements.
\end{abstract}

\setcounter{tocdepth}{2}
\tableofcontents



\chap{Introduction}\label{chap.intro}

\firstpage

Interval exchange transformations (\IET s) are encoded by a permutation $\pi$ and length vector $\lambda$. In \cite{c.R79}, Rauzy introduces Rauzy induction, a first return map of an \IET\ on a specific subinterval. This induction takes one of two forms on the space of \IET s and therefore descends to two different maps on the set of permutations. Therefore permutations are divided into \term{Rauzy Classes}, minimal sets closed under the two types of induction maps. We dedicate Sections \ref{sec.IET} and \ref{sec.RC} to providing some basic background and well known results concerning both \IET s and Rauzy Classes.\\

From another direction, we consider the \term{moduli space of Abelian differentials}. By the zippered rectangle construction in \cite{c.Ve82}, Veech shows that a generic \IET\ is uniquely ergodic (a result independently proved by Masur in \cite{c.Mas82}). This construction establishes a relationship between an \IET\ and flat surfaces with oriented measured foliations. We present an equivalent construction, called a \term{suspension}, in Section \ref{sec.surface}. Using suspensions, we assign properties to a permutation $\pi$: its \term{signature} (see Definition \ref{def.signature}), which is related to the singularities of these suspensions, and its \term{type} (see Section \ref{sec.classify}), which represents any other necessary data from its suspensions. The crucial result in this \chapword\ is the following:
\begin{corClassify}
    Every Rauzy class is uniquely determined by signature and type. So given Rauzy class $\RClass$, if $\pi\in\irr$ has the same signature and type as $\RClass$, then necessarily $\pi\in\RClass$.
\end{corClassify}
\noindent This immediately follows from \cite{c.KZ} and \cite{c.B09}. In Sections \ref{sec.hyperelliptic} and \ref{sec.spin}, we discuss \term{hyperelliptic} surfaces and the \term{parity} of a surface's \term{spin structure}. These discussions give us the necessary information to determine a permutation's type. In Sections \ref{sec.gen.perm}-\ref{sec.gen.classes}, we discuss a generalization of \IET s, called \term{linear involutions}. Introduced in \cite{c.DN90}, these give rise to \term{generalized permutations}. Such linear involutions arise from flat surfaces with non-orientable measured foliations. We discuss Rauzy-Veech induction on linear involutions and \term{generalized Rauzy Classes}. We conclude in Section \ref{sec.gen.classify} that Corollary \ref{cor.classify} applies to these classes as well, as a direct result from the works \cite{c.K97},\cite{c.MS} and \cite{c.L08} with the final result from \cite{c.B09}.\\

In Equations (2.2), Veech (see \cite{c.Ve84_I}) shows a definition of Rauzy induction on permutations. It is clear from this definition that the map $\pi\mapsto\inv\pi$ conjugates one type of induction with the other. This natural relationship conjures two natural questions:
\begin{enumerate}
    \item Are Rauzy Classes closed under taking inverses?
    \item Do all Rauzy classes contain self-inverse permutations?
\end{enumerate}
\noindent The work leading up to Corollary \ref{cor.classify} in Sections \ref{sec.classify} and \ref{sec.gen.classify} provides an affirmative to the first question: any suspensions of $\pi$ and $\inv\pi$ have the same signature and type and therefore $\pi$ and $\inv\pi$ belong to the same class. However, proving a positive result for the second question would naturally imply one for the first also. This work answers the second question.
\begin{thmMainTrue}
	Every (true) Rauzy Class contains a permutation $\pi$ such that $\pi=\inv{\pi}$.
\end{thmMainTrue}
\noindent In Section $\ref{sec.blocks}$, we form patterns of letters, or \term{blocks}, that we may use to construct a self-inverse $\pi$ such that $\pi\in\RClass$ by Corollary \ref{cor.classify}. This method follows in the spirit of \cite{c.Z08}. In that paper, Zorich constructs permutations with desired properties. He then shows that these permutations belong to the desired Rauzy Class in a fashion similar to Corollary \ref{cor.classify}.\\

In \Chapword\ \ref{chapGen}, we prove the generalized Rauzy Class analogue to Theorem \ref{thm.main}.
\begin{thmMainGen}
    A generalized Rauzy Class $\RClass\subseteq\genirr_\AAA$ with signature $(\ell_1,\dots,\ell_m)$ contains a permutation $\pi$ such that $\pi=\inv{\pi}$ if and only if the following conditions hold:
        \begin{subequations}
            \begin{equation*}
                \ell_1=0,
            \end{equation*}
            \begin{equation*}
                \#\{i:\ell_i=j,~1\leq i \leq m\}\mbox{ is even}, \mbox{ for all odd }j\geq-1.
            \end{equation*}
        \end{subequations}
\end{thmMainGen}
\noindent These conditions are dependent only on the signature $\sig(\RClass)$. We still explicitly construct a self-inverse $\pi$ for any given $\RClass$ satisfying these conditions, but this construction is no longer obtained by a simple concatenation of blocks. Rather, we insert blocks into generalized permutations and verify this effect on the signature (see Section \ref{sec.gen.blocks}).\\

We consider the topic of Lagrangian subspaces of suspensions in \Chapword\ \ref{chapLag}. We call a permutation $\pi$ \term{Lagrangian} if the vertical trajectories of any suspension of an \IET\ $T\sim(\pi,\mone)$, where $\mone=(1,\dots,1)$, span a $g$-dimensional subspace in homology, where $g$ is the genus of $\pi$. We prove the following:
\begin{thmMainLag}
        Suppose $\pi\in\irr_d$ is self-inverse. Then $\pi$ is Lagrangian.
\end{thmMainLag}
\noindent This theorem provides an alternative proof of Forni's Lemma 4.4 in \cite{c.For02}. We present this proof as Corollary \ref{cor.Lag.For}. In this lemma, Forni shows that the set of $q\in\HHH_g$ (the moduli space of Abelian differentials of genus $g$) such that
\begin{enumerate}
    \item The vertical trajectories of $q$ are (almost all) periodic,
    \item These trajectories span a $g$-dimensional subspace in homology,
\end{enumerate}
is a dense set in $\HHH_g$. Corollary \ref{cor.Lag.For} uses Theorem \ref{t.Lag.Perm} and the fact that the Teichm\"uller geodesic flow is generically dense in each connected component of $\HHH_g$. We further show that the permutations we construct in \Chapword\ \ref{chapTrue} need only consider the transposition pairs (letters interchanged by the permutation) to form such a basis.

\sect{Interval Exchange Transformations}\label{sec.IET}

    Let $\perm_d$ be the set of permutations on $\dset{d}$. $\pi\in\perm_d$ is \term{irreducible} if $\pi(\dset{k})=\dset{k}$ only when $k=d$. The set of all irreducible permutations on $\dset{d}$ is $\irr_d$. The notation $\pi\in\irr$ shall be used to indicate that $\pi$ is irreducible when it is not necessary to state $d$. $\pi\in\perm_d$ is \term{standard} if $\pi(d)=1$ and $\pi(1)=d$. Note that a standard permutation is necessarily irreducible.\\
	
    When we refer to a (sub)interval, we mean open on the right and closed on the left (i.e. of the form $[a,b)$, for some $a<b$). Let $\RR_+^d$ be the cone of positive length vectors in $\RR^d$. For $\lambda\in\RR_+^d$, let $|\lambda|:=\dusum{i=1}{d}\lambda_i$, $I:=[0,|\lambda|)$, and define subintervals $I_i\subseteq I$, $i\in\dset{d}$ as $I_i:=[\dsum{j<i} \lambda_j,\dsum{j\leq i} \lambda_j)$.
	
	\begin{defn}\label{def.IET}
		An \term{interval exchange transformation} (\IET) is a bijective map $T:I\rightarrow I$ such that there exists a
		partition of $I$ into subintervals $I_1\dots I_d$ such that for each $i$, $T|_{I_i}$ is a translation.
	\end{defn}

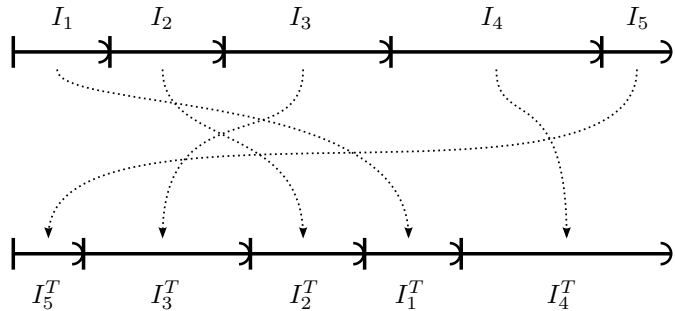
\begin{figure}[h]
        \begin{center}
           \setlength{\unitlength}{250pt}
            \begin{picture}(1,.4)
                \put(0,0){\includegraphics[width=\unitlength]{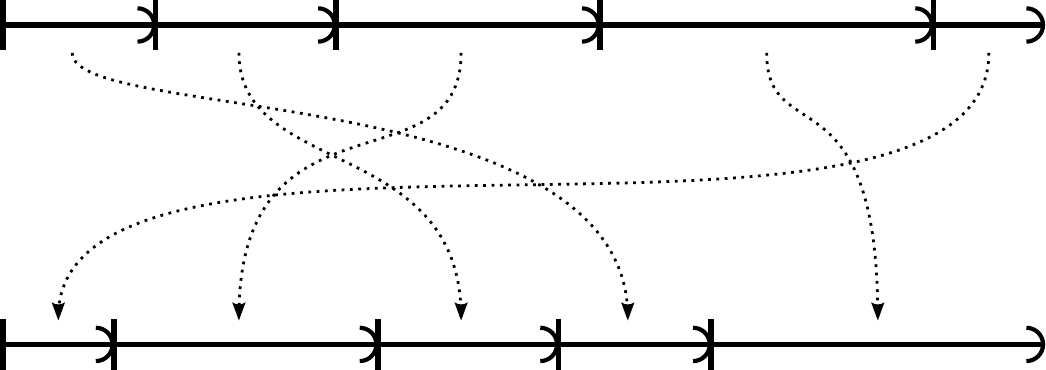}}

                \put(0.03,-.05){$I^T_5$}
                \put(0.21,-.05){$I^T_3$}
                \put(0.42,-.05){$I^T_2$}
                \put(0.58,-.05){$I^T_1$}
                \put(0.81,-.05){$I^T_4$}

                \put(.06,.37){$I_1$}
                \put(.21,.37){$I_2$}
                \put(.42,.37){$I_3$}
                \put(.71,.37){$I_4$}
                \put(.93,.37){$I_5$}
            \end{picture}
       \end{center}
       \caption{The mapping of \IET\ $T$ with $\pi=(5,3,2,1,4)$.}\label{fig.IETveech}
    \end{figure}

    Because an \IET\ is a piecewise translation, it takes subintervals and reorders them in $I$. If we choose these $d$ subintervals, we shall encode this reordering as $\pi\in\perm_d$ in such a way that $\pi(i)$ indicates the position of $I_i$ after the map. Let $I^T_i:=T(I_i)$. If we start with our $I_i$'s in order, they are moved into the order $I^T_{\inv\pi(1)},\dots,I^T_{\inv\pi(d)}$ under $T$ (see Figure \ref{fig.IETveech}).

    By convention, we label $\pi\in\perm_d$ as $\pi=(\pi^{-1}(1), \dots, \pi^{-1}(d))$ to indicate the ordering of the original subintervals after the \IET. Given an \IET\ $T$, we can choose $\lambda\in\RR_+^d$ such that $\lambda_i=|I_i|$ for each $i$. Given the pair $(\pi,\lambda)$, we may describe the \IET\ $T$ explicitly by
    	$$ T(x) = x + \omega_i,~~x\in I_i$$
    where $\omega_i=\sum_{j:\pi(j)<\pi(i)}\lambda_j - \sum_{j: j<i}\lambda_j $. This vector $\omega=(\omega_1,\dots\omega_d)$ is called the \term{translation vector} for $T$ and may be described by a matrix $\Omega_\pi$ by $\omega=\Omega_\pi\lambda$ where,
    \begin{equation}\label{eq.omega_pi}
        (\Omega_\pi)_{i,j}=\RHScase{1,& \mbox{if }i<j~\&~\pi(i)>\pi(j),\\
								-1,& \mbox{if }i>j~\&~\pi(i)<\pi(j),\\
								0,& \mathrm{otherwise.}}
    \end{equation}

    \begin{rem}
        The matrix $\Omega_\pi$ is the same as the matrix $L^\pi$ seen in \cite{c.Ve78} and $M$ in \cite{c.R79}. For an example of the notation in Equation \eqref{eq.omega_pi}, see \cite{c.Vi06}. This applies to the other definition of $\Omega_\pi$ given in Equation \eqref{eq.omega_pi2} as well.
    \end{rem}
	
	\begin{defn}\label{def.veech}
		Using the discussed notation, an \IET\ $T$ may be described by a pair $(\pi,\lambda)$, which
		shall be denoted $T\sim(\pi,\lambda$).
		When this convention is used, we shall denote this as \term{\veech} notation.
	\end{defn}
		
\begin{figure}[b]
        \begin{center}
           \setlength{\unitlength}{250pt}
            \begin{picture}(1,.36)
                \put(0,0){\includegraphics[width=\unitlength]{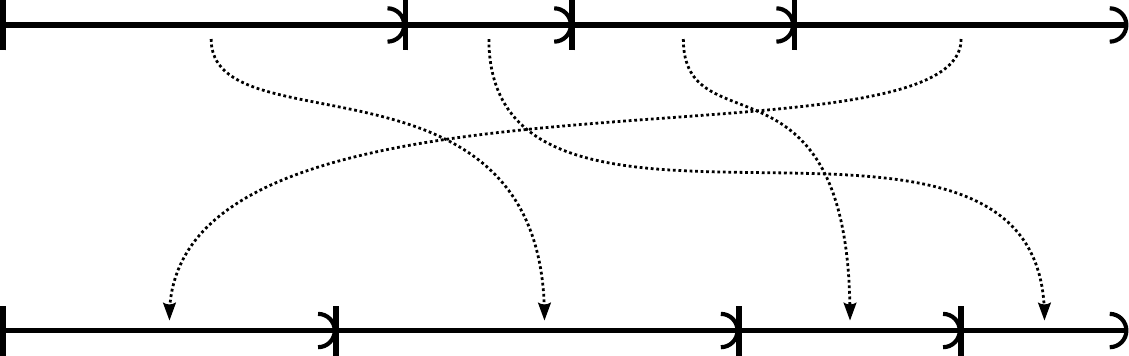}}

                \put(0.18,.33){$I^0_a$}
                \put(0.41,.33){$I^0_b$}
                \put(0.59,.33){$I^0_c$}
                \put(0.83,.33){$I^0_d$}

                \put(.13,-.04){$I^1_d$}
                \put(.46,-.04){$I^1_a$}
                \put(.73,-.04){$I^1_c$}
                \put(.90,-.04){$I^1_b$}
            \end{picture}
       \end{center}
       \caption{An \IET\ on $\AAA=\{a,b,c,d\}$.}\label{fig.IETyoccoz}
    \end{figure}

    We now discuss another notation for representing \IET's (see \cite{c.Ker85}, \cite{c.MMY05} and \cite{c.Bu06}). Let $\AAA$ be a finite alphabet and $d:=\#\AAA$. Let $\pi_i:\AAA\rightarrow\{1,\dots,d\}$, $i\in\{0,1\}$, be bijections. Now we let $\pi_0$ describe the ordering of subintervals $I^0_\alpha$, $\alpha\in\AAA$, and $\pi_1$ shall indicate what position the image of $I^0_\alpha$ takes under $T$, which we shall denote as $I^1_\alpha$. Consider the cone $\RR_+^\AAA\subset\RR^\AAA$, and for $\lambda\in\RR_+^\AAA$, define subintervals $I^\eps_\alpha$ of $I:=[0,|\lambda|)$, $\alpha\in\AAA$ and $\eps\in\{0,1\}$, by $I^\eps_\alpha = [\dsum{\beta: \pi_\eps(\beta)<\pi_\eps(\alpha)}\lambda_\beta, \dsum{\beta: \pi_\eps(\beta)\leq \pi_\eps(\alpha)}\lambda_\beta)$. The \IET\ $T$ translates the intervals by $T(I^0_\alpha)=I^1_\alpha$ as indicated in Figure \ref{fig.IETyoccoz}.
    We shall denote $(\pi_0,\pi_1)$ by
        $$(\pi_0,\pi_1)=\cmtrx{\pi_0^{-1}(1)&\dots & \pi_0^{-1}(d)\\ \pi_1^{-1}(1)&\dots & \pi_1^{-1}(d)}$$
    indicating the orders of the subintervals before and after the application of $T$. As in the case of \veech\ notation, we have a translation vector $\omega\in\RR^\AAA$ such that
	$$ T(x)=x+\omega_\alpha,~~x\in I^0_\alpha,$$
	$\omega$ may again be described by a matrix $\Omega_{(\pi_0,\pi_1)}$ by $\omega=\Omega_{(\pi_0,\pi_1)}\lambda,$ where
    \begin{equation}\label{eq.omega_pi2}
        (\Omega_{(\pi_0,\pi_1)})_{\alpha,\beta}=\left\{\begin{array}{rl} 1,& \mbox{if }\pi_0(\alpha)<\pi_0(\beta)~\&~\pi_1(\alpha)>\pi_1(\beta),\\
								-1,& \mbox{if }\pi_0(\alpha)>\pi_0(\beta)~\&~\pi_1(\alpha)<\pi_1(\beta),\\
								0,& \mathrm{otherwise.}\end{array}\right.
    \end{equation}
	
	\begin{defn}\label{def.yoccoz}
		When we associate to $T$ a triple $(\pi_0,\pi_1,\lambda)$, we denote this by $T\sim(\pi_0,\pi_1,\lambda)$.
		When this convention is used, we refer to this as \term{\yoccoz} notation.
	\end{defn}
	
    It should be noted that there is a natural association between pairs $(\pi_0,\pi_1)$ and permutations $\pi$. For pairs $(\pi_0,\pi_1)$ we have a map onto $\perm_d$, the \term{monodromy invariant}, by $(\pi_0,\pi_1)\mapsto \pi_1\circ\pi_0^{-1}$. Given any bijection $b:\AAA'\rightarrow\AAA$ such that $\#\AAA=\#\AAA'$, $(\pi_0\circ b,\pi_1\circ b)$ and $(\pi_0,\pi_1)$ have the same monodromy invariant. We use this as a natural equivalence relation between pairs: two pairs are equivalent if and only if they have the same monodromy invariant. We then use the notation $\pi=(\pi_0,\pi_1)$ when $(\pi_0,\pi_1)\sim(\pi_0\circ\pi_0^{-1},\pi_1\circ\pi_0^{-1})\sim(id,\pi)$ where $id$ is the identify function on the alphabet $\dset{d}$. This association applies also to $(\pi,\lambda)$ and $(\pi_0,\pi_1,\lambda')$. $(\pi,\lambda)=(\pi_0,\pi_1,\lambda')$ if $\pi=(\pi_0,\pi_1)$ as above and $\lambda_i=\lambda'_{\pi_0^{-1}(i)}$. In this paper, we do work and calculations in \yoccoz\ notation, while understanding the underlying statements for \veech\ notation. For example, $\pi=(\pi_0,\pi_1)$ is \term{standard} if the pair is of the form $ \pi=\cmtrx{\alpha & \dots & \beta\\\beta& \dots & \alpha}$. If $\pi=(\pi_0,\pi_1)$, then $\pi^{-1}=(\pi_1,\pi_0)$. So a pair is self inverse if $(\pi_0,\pi_1)\sim(\pi_1,\pi_0)$, i.e. if $\alpha$ appears above $\beta$ in the pair then $\beta$ also appears above $\alpha$. So $\pi=\cmtrx{a&b&c&d\\d&b&c&a}$ is self-inverse, while $\pi=\cmtrx{a&b&c\\c&a&b}$ is not.\\

\sect{Rauzy Classes}\label{sec.RC}

    In this section, we define a family of maps on irreducible permutations, known as Rauzy induction. Introduced in \cite{c.R79}, this is realized as a first return map of an \IET\ on appropriate subintervals. These moves partition each set $\irr_d$ into equivalence classes under induction. We then state a relationship between induction and the map $\pi\mapsto\inv\pi$ in Claim \ref{cor.2}, observed by Veech.\\

    Given an \IET\ $T\sim(\pi_0,\pi_1,\lambda)$, $\pi=(\pi_0,\pi_1)\in\perm_d$, let $\alpha_\eps=\pi_\eps^{-1}(d)$. Define  $I':=[0,|\lambda|-\min\{\lambda_{\alpha_0},\lambda_{\alpha_1}\})$. Then the first return map $T'$ of $T$ on $I'$ is again an \IET\ with $T'\sim(\pi_0',\pi_1',\lambda')$ defined by the following rules:\\
	
    \begin{defn}\label{def.RV}
        We define Rauzy induction on $T\sim(\pi_0,\pi_1,\lambda)$ as $T'\sim(\pi_0',\pi_1',\lambda')$ by the following:
        \begin{itemize}
            \item Assume $\lambda_{\alpha_0}>\lambda_{\alpha_1}$. We call this \term{Rauzy induction of type $\rtb$}. Then $\pi'=(\pi_0',\pi_1')$ is defined by the following rules:

                $$ \pi_0'=\pi_0,\mbox{ and }\pi_1'(\alpha)=\RHScase{
							\pi_1(\alpha), & \mbox{if }\pi_1(\alpha)\leq\pi_1(\pi_0^{-1}(d)),\\
							\pi_1(\alpha)+1, & \mbox{if }\pi_1(\pi_0^{-1}(d))<\pi_1(\alpha)<d,\\
							\pi_1(\pi_0^{-1}(d))+1, &\mbox{if }\pi_1(\alpha)=d,}$$
                or by the following diagram
                    $$ \xymatrix{{\pi=\cmtrx{\dots& & & \dots& \alpha_0\\
    			     \dots& \alpha_0& \beta&\dots&\alpha_1}}
    		          \ar[r]^\rtb &
    		      {\cmtrx{\dots& & & \dots& \alpha_0\\
    			             \dots& \alpha_0&\alpha_1& \beta&\dots}=\pi',}}$$
    	       and $\lambda$ is related to $\lambda'$ by
        	       $$ \lambda_\alpha'=\RHScase{
    					\lambda_{\alpha_0}-\lambda_{\alpha_1},&\mbox{if }\alpha=\alpha_0,\\
    					\lambda_{\alpha}, & \mathrm{otherwise.}}$$
				
            \item Now assume $\lambda_{\alpha_0}<\lambda_{\alpha_1}$. This is \term{Rauzy induction of type $\rtt$}. Then $\pi'=(\pi_0',\pi_1')$ is defined by the following rules:
            	$$ \pi_1'=\pi_1,\mbox{ and }\pi_0'(\alpha)=\RHScase{
							\pi_0(\alpha), & \mbox{if }\pi_0(\alpha)\leq\pi_0(\alpha_1),\\
							\pi_0(\alpha)+1, & \mbox{if }\pi_0(\alpha_1)<\pi_0(\alpha)<d,\\
							\pi_0(\alpha_1)+1, &\mbox{if }\alpha=\alpha_0,}$$
            	or by the following diagram
            	   $$ \xymatrix{{\pi=\cmtrx{\dots& \alpha_1& \beta&\dots&\alpha_0\\
            			\dots& & & \dots& \alpha_1}} \ar[r]^\rtt &
            	       {\cmtrx{\dots& \alpha_1&\alpha_0& \beta&\dots\\
            			\dots& & & \dots& \alpha_1}=\pi',}}$$
            	and $\lambda$ is related to $\lambda'$ by
            	   $$ \lambda_\alpha'=\RHScase{
			     		\lambda_{\alpha_1}-\lambda_{\alpha_0},&\mbox{if }\alpha=\alpha_1,\\
    					\lambda_{\alpha}, & \mathrm{otherwise.}}$$
        \end{itemize}
				
        We shall denote $\pi'$ as $\rb\pi$ or $\rt\pi$ if the induction was of type $\rtb$ or $\rtt$, respectively.
    \end{defn}

    \begin{rem}
        The case $\lambda_{\alpha_0}=\lambda_{\alpha_1}$ does not have a valid definition, as the resulting induced transformation is over $(d-1)$ symbols. However, such $\lambda$'s form a codimension one (therefore Lebesgue measure zero) set in $\RR_+^\AAA$.
    \end{rem}

    \begin{defn}\label{def.Keane}
        Assume $\AAA=\dset{d}$ and $\pi=(id,\pi)$. Let $T\sim(\pi,\lambda)$ and $\partial I_i$ denote the left endpoint of subinterval $I_i$ for $i\in\dset{d}$. Then $T$ satisfies the \term{Keane Condition} if
            \begin{equation}\label{eq.connexion}
                \underbrace{T\circ\cdots\circ T}_{m}(\partial I_i) = T^m(\partial I_i)\neq \partial I_j
            \end{equation}
        for all $i,j\in\dset{d}$, $j>1$ and $m\geq 1$.
    \end{defn}

    \begin{rem}
        Each violation of the Keane condition satisfies Equation \eqref{eq.connexion} for a certain triple $(i,j,m)$. However, each of these conditions is a codimension one set in $\RR_+^d$. So given $\pi\in\irr_d$, we see that the Keane property is satisfied for Lebesgue almost everywhere in $\RR_+^d$.
    \end{rem}

    \begin{prop}\label{prop.Keane_is_inducible}
        Let $T^{(n)}$ denote the $n^{th}$ iteration of induction on \IET\ $T$. Then the following are equivalent:
        \begin{itemize}
            \item $T$ satisfies the Keane condition.
            \item $T^{(n)}$ is defined for all $n\geq 0$.
        \end{itemize}
    \end{prop}

	\begin{claim}
        Let $\eps\in\{0,1\}$. Then
            $$\pi\in\irr_d \iff \reps\pi\in\irr_d.$$
	\end{claim}
	
	\begin{proof}
        Suppose $\pi=(\pi_0,\pi_1)\in\perm_d\setminus\irr_d$, and fix type $\eps$. Then there is a proper subset $\AAA'\subset\AAA$, $\#\AAA=k<d$, such that $\pi_i(\AAA')=\{1,\dots,k\}$, $i\in\{0,1\}$. Most importantly, $\alpha_i=\pi_i^{-1}(d)\notin\AAA'$. So our induction must only move elements of $\AAA\setminus\AAA'$ as every element of $\AAA'$ appears before every element of $\AAA\setminus\AAA'$ in both rows. Namely $\reps\pi_i(\AAA')=\{1,\dots,k\}$ as well, or $\reps\pi\in\perm_d\setminus\irr_d$. The argument above applies if we first assume $\reps\pi\in\perm_d\setminus\irr_d$ and evaluate $\pi$.
	\end{proof}
	
	So Rauzy induction is a closed operation in the set $\irr_d$.
	
	\begin{defn}\label{def.rauzyclass}
        Given $\pi\in\irr_d$, the \term{Rauzy Class} of $\pi$, $\RClass(\pi)\subseteq\irr_d$, is the orbit of $\rb$ and $\rt$ moves on $\pi$. The \term{Rauzy Graph} of $\pi$ is the graph with vertices in $\RClass(\pi)$ and directed edges corresponding to the inductive moves.
	\end{defn}
	
	\begin{exam}
        Start with permutation $\pi=(3,2,1)$, we have the following two other elements
            $$ \rb\pi=\cmtrx{1 & 2 & 3 \\ 3 & 1 & 2}=(3,1,2)\mbox{, } \rt\pi=\cmtrx{1 & 3 & 2 \\ 3 & 2 & 1}=\cmtrx{1 & 2 & 3 \\ 2 & 3 & 1}=(2,3,1)$$
	   The Rauzy Graph for $\RClass(\pi)$ is listed in Figure \ref{fig.RC321}.
    \end{exam}
    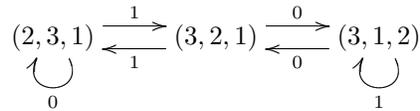
\begin{figure}[b]
	$$ \xymatrix{ (2,3,1) \ar@(dr,dl)^\rtb \ar@<1ex>[r]^\rtt& (3,2,1) \ar@<1ex>[l]^\rtt \ar@<1ex>[r]^\rtb & (3,1,2)
		\ar@(dr,dl)^\rtt \ar@<1ex>[l]^\rtb}$$
    \caption{The Rauzy Graph for $(3,2,1)$.}\label{fig.RC321}
    \end{figure}
	The given definition of a Rauzy Class is dependent on the choice of $\pi$, but the next claim shows that being in the
	same Rauzy Class is an equivalence condition and not dependent on our choice of representative.
	\begin{claim}\label{cor.1}
		For any $\pi^{(1)},\pi^{(2)}\in\RClass(\pi)$, there exists a directed path from $\pi^{(1)}$ to $\pi^{(2)}$ in
		the Rauzy Graph.
	\end{claim}
	
	\begin{proof}
        It suffices to show that for a permutation $\tilde{\pi}\in\RClass(\pi)$, there exists a path from $\reps\tilde{\pi}$ to $\tilde{\pi}$ for $\eps\in\{0,1\}$. By Definition \ref{def.RV}, there exists $n>0$ such that $\reps^n\tilde{\pi}=\tilde{\pi}$. So $n-1$ moves of type $\eps$ form a path from $\reps\tilde{\pi}$ to $\tilde{\pi}$.
	\end{proof}
	
    So if $\tilde{\pi}\in\RClass(\pi)$, then $\RClass(\pi)=\RClass(\tilde{\pi})$. The next result is used in Sections \ref{sec.hyperelliptic} and \ref{sec.hyperelliptic}.

    \begin{claim}\label{cor.std_in_RC}
        Every Rauzy Class $\RClass\subset\irr_d$ contains a standard permutation (i.e. $\pi$ such that $\pi(d)=1$ and $\pi(1)=d$).
    \end{claim}

    \begin{proof}
        Consider any $\pi\in\RClass$. Denote by $\alpha_\eps=\inv\pi_\eps(d)$ and $\beta_\eps=\inv\pi_\eps(1)$, for $\eps\in\{0,1\}$, or
        $$ \pi = \cmtrx{\LL{\beta_0}{\beta_1}~\LL{\dots}{\dots}~\LL{\alpha_0}{\alpha_1}}.$$ Let $n = \min\{\pi_0(\alpha_1), \pi_1(\alpha_0)\}$. Suppose $n=1$ and choose $\eps$ such that $\pi_{1-\eps}(\alpha_\eps)=1$. In this case, $\alpha_\eps = \beta_{1-\eps}$, and if we perform $m$ inductive moves of type $\eps$, the resulting permutation is standard, where $m=d-\pi_{1-\eps}(\beta_\eps)$. Consider the following diagram for $\eps=\rtb$:
        $$ \xymatrix{{\cmtrx{\LL{\beta_0}{\alpha_0} ~\LL{\dots}{\dots}~\LL{\dots}{\beta_0~\delta} ~\LL{\dots}{\dots}~\LL{\alpha_0}{\alpha_1}}} \ar^{\rtb^m}[r] &
            {\cmtrx{\LL{\beta_0}{\alpha_0} ~\LL{\dots}{\delta\dots}~\LL{\dots}{\alpha_1} ~\LL{\dots}{\dots}~\LL{\alpha_0}{\beta_0}}}}.$$

        If $n>1$ then we may fix $\eps\in\{0,1\}$ and find $\gamma\in\AAA$ such that $\pi_\eps(\gamma) < n < \pi_{1-\eps}(\gamma)$. If no such $\gamma$ exists, then $\inv{\pi_0}(\{n,\dots,d\})=\inv{\pi_1}(\{n,\dots,d\})$ and $\pi$ is not irreducible. So let $m=d-\pi_{1-\eps}(\gamma)$, and perform $m$ iterations of type $\eps$. Call this new permutation $\pi'$ and note that $\alpha'_{1-\eps}=\inv{\pi'}_{1-\eps}(d) = \gamma$ and $\alpha'_\eps = \inv{\pi'}_\eps(d)=\alpha_\eps$. Therefore $n' = \min\{\pi'_0(\alpha'_1), \pi'_1(\alpha'_0)\}<n$. Consider the following diagram for $\eps=\rtt$:
        $$ \xymatrix{\cmtrx{\LL{\dots}{\dots\gamma} ~\LL{\alpha_1\dots}{\dots\alpha_0} ~\LL{\gamma\dots}{\dots} ~\LL{\alpha_0}{\alpha_1}} \ar^{\rtt^m}[r] &
            \cmtrx{\LL{\dots}{\dots\gamma}~\LL{\alpha_1\dots\alpha_0}{\dots\alpha_0}~\LL{\dots}{\dots}~\LL{\gamma}{\alpha_1}}}.$$

        Repeat the above argument for $\pi'$ until $n'=1$ and we may derive a standard permutation.
    \end{proof}

    Consider one more observation that is used in Corollary \ref{cor.main}. This result is evident from Equations (2.2) in \cite{c.Ve84_I}.
	
	\begin{claim}\label{cor.2}
		For $\eps\in\{0,1\}$ and $\pi\in\irr_d$, $\reps\inv{\pi}=\inv{(\ropp\pi)}$.
	\end{claim}
	
	\begin{proof}
        We will show that $\inv{(\rb\inv\pi)}=\rt\pi$ as it will prove the claim for all cases. Let $\pi =(\pi_0,\pi_1)$.
        Then $\inv\pi = (\pi_1,\pi_0)$. Now let $\rb\inv\pi=(\pi_0^{\bullet},\pi_1^{\bullet})$. By Definition \ref{def.RV},
        $$ \pi_0^{\bullet} = \pi_1 \mbox{ and } \pi_1^{\bullet}(\alpha) =
            \RHScase{\pi_0(\alpha), & \mbox{if }\pi_0(\alpha)\leq\pi_0(\inv\pi_1(d)),\\
                \pi_0(\alpha)+1, & \mbox{if }\pi_0(\inv\pi_1(d)) < \pi_0(\alpha) < d, \\
                \pi_0(\inv\pi_1(d))+1, & \mbox{if }\pi_0(\alpha)=d.}$$
        Then $\inv{(\rb\inv\pi)} = (\pi_1^\bullet,\pi_0^\bullet)$. By checking Definition \ref{def.RV}, we conclude that $\rt\pi = (\pi_1^\bullet,\pi_0^\bullet) = \inv{(\rb\inv\pi)}$.
	\end{proof}

	So the action of taking the inverse permutation conjugates with the Ruazy moves on $\pi$ by sending them to
	the opposite move on $\pi^{-1}$.

    \begin{rem}\label{rem.RV_is_2to1}
         Let $T\sim(\pi,\lambda)$ for $\pi=(\pi_0,\pi_1)\in\irr_d$ and $\lambda\in\RR_+^\AAA$. We also let $\alpha_\eps=\inv\pi_\eps(d)$ define the last letter of the rows of $\pi$. From the proof of Claim \ref{cor.1}, we can define $\pi^\eps$, for $\eps\in\{0,1\}$, such that $\reps\pi^\eps=\pi$. Also for $\eps\in\{0,1\}$, let $\lambda^\eps\in \RR_+^\AAA$ be defined by
             $$ \lambda^\eps_\alpha = \RHScase{ \lambda_{\alpha_0}+\lambda_{\alpha_1}& \mbox{if }\alpha=\alpha_\eps,\\ \lambda_\alpha& \mathrm{otherwise.}}$$
         It follows that if $T_\eps\sim(\pi^\eps,\lambda^\eps)$, then $T_\eps' = T$ and induction on $T_\eps$ is type $\eps$. So almost everywhere on the set $\irr_d\times\RR_+^\AAA$, Rauzy induction is a $2$ to $1$ map.
    \end{rem}

\sect{Suspended Surfaces for Interval Exchanges}\label{sec.surface}

    In \cite{c.Ve82}, Veech introduced the zippered rectangle construction, which allows us to associate to an \IET\ a flat surface with an Abelian differential. We present an equivalent construction, presented for example by Viana in \cite{c.Vi06}, of suspended surfaces over an \IET. We discuss Rauzy-Veech induction on these surfaces and introduce the moduli space of Abelian differentials.\\
	
\begin{figure}[h]
        \begin{center}
           \setlength{\unitlength}{200pt}
            \begin{picture}(1,.57)
                \put(0,0){\includegraphics[width=\unitlength]{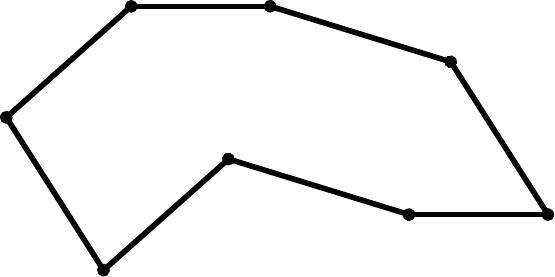}}

                \put(0.09,.41){$1$}
                \put(0.35,.51){$2$}
                \put(0.66,.46){$3$}
                \put(0.92,.26){$4$}

                \put(.06,.09){$4$}
                \put(.32,.07){$1$}
                \put(.55,.10){$3$}
                \put(.85,.05){$2$}
            \end{picture}
       \end{center}
        \caption{A suspension surface for $(4,1,3,2)$.}\label{fig.suspend4132}
    \end{figure}
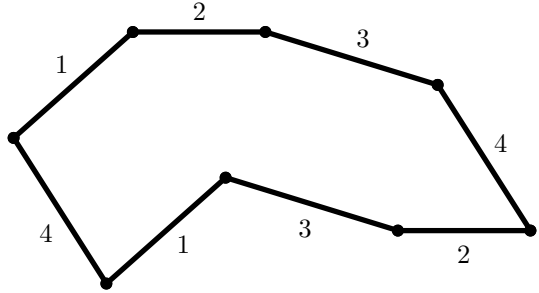

\begin{figure}[h]
        \begin{center}
           \setlength{\unitlength}{200pt}
            \begin{picture}(1,.80)
                \put(0,0){\includegraphics[width=\unitlength]{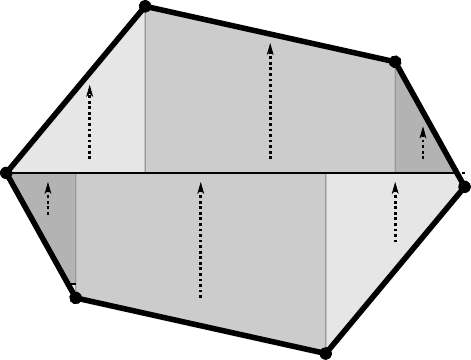}}

                \put(.13,.62){$1$}
                \put(.61,.73){$2$}
                \put(.94,.52){$3$}

                \put(.04,.21){$3$}
                \put(.38,.00){$2$}
                \put(.86,.14){$1$}
            \end{picture}
       \end{center}
        \caption{The first return of $I_S$ in the suspension is the original \IET.}\label{fig.firstreturnisIET}
    \end{figure}
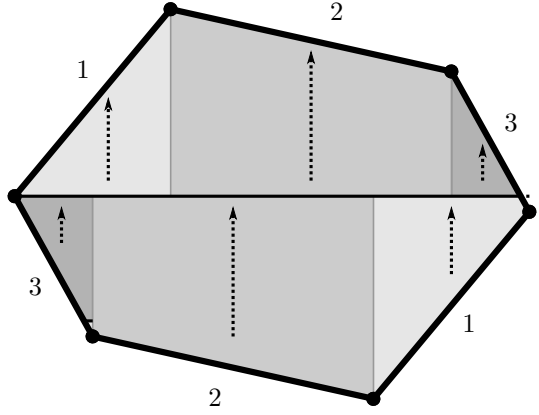

	Fix $\pi=(\pi_0,\pi_1)\in\irr_d$, $\pi_\eps:\AAA\rightarrow\dset{d}$, and $\lambda\in\RR_+^\AAA$. Let
    \begin{equation}\label{eq.Tpi}
        \TTT_\pi:=\left\{\tau\in\RR^\AAA: \sum_{\alpha: \pi_0(\alpha)\leq k}\tau_\alpha>0,
            \sum_{\alpha: \pi_1(\alpha)\leq k}\tau_\alpha<0\mbox{, for all }1\leq k < d\right\}.
    \end{equation}
    Define vectors $\vec{\zeta}_\alpha:=(\lambda_\alpha,\tau_\alpha)$ and segments $\zeta_\alpha^\eps$, $\eps\in\{0,1\}$, as the segment starting at $\sum_{\beta: \pi_\eps(\beta) <\pi_\eps(\alpha)}\vec{\zeta}_\beta$ and ending at $\sum_{\beta: \pi_\eps(\beta) \leq\pi_\eps(\alpha)}\vec{\zeta}_\beta$, noting that $\zeta_\alpha^0$ and $\zeta_\alpha^1$ are parallel (they are just translations of vector $\vec{\zeta}_\alpha$). Let $S:=S(\pi,\lambda,\tau)$ be the surface bounded by all $\zeta_\alpha^\eps$ with each $\zeta_\alpha^0$ and $\zeta_\alpha^1$ identified by translation. For example, if $\pi=(4,1,3,2)$, one suspension is given by Figure \ref{fig.suspend4132}.
    To avoid cumbersome notation, we denote the segments $\zeta_\alpha^\eps$ simply by $\alpha$ in a suspension. By definition, the leftmost endpoint is $(0,0)$. Define $I_S:=[0,|\lambda|)\times\{0\}$. With the exception of the points of discontinuity, the \IET\ $T\sim(\pi,\lambda)$ is realized by the first return of the positive vertical direction of $S$ on $I_S$, as is illustrated Figure \ref{fig.firstreturnisIET}.
			
    Each of the identifications on these surfaces is a translation. Therefore the standard form $dz$ in the polygon descends to a holomorphic $1$-form on the surface with zeroes, if any, at the vertex equivalence classes. Each vertex class is called a singularity of degree $k$, where $k$ is the degree of the corresponding zero of the differential and the total angle around the singularity is $2\pi(k+1)$.\\

    In order to give an explicit way to determine the degree of the singularities in a surface $S$, let us label the endpoints of our segments by $(\alpha,\eps,\imath)$, where $\alpha\in\AAA$, $\eps\in\{0,1\}$ and $\imath\in\{L,R\}$, to denote the left or right endpoint of segment $\zeta_\alpha^\eps$. We have the natural identification rules:
    \begin{enumerate}
        \item For $1\leq i < d$ and $\eps\in\{0,1\}$, $(\inv\pi_\eps(i),\eps,R)\sim(\inv\pi_\eps(i+1),\eps,L)$.
        \item $(\inv\pi_0(1),0,L)\sim(\inv\pi_1(1),1,L)$ and $(\inv\pi_0(d),0,R)\sim(\inv\pi_1(d),1,R)$.
        \item For $\alpha\in\AAA$ and $\imath\in\{L,R\}$, $(\alpha,0,\imath)\sim(\alpha,1,\imath)$.
    \end{enumerate}
    The equivalence sets determine the identified singularities in our surface $S$. The first rule lets us consider only the vertices of the form $(\alpha,\eps,L)$. With the exception of $(\inv\pi_0(1),0,L)$, every other vertex of this form has a downward direction in $S$. If a singularity is of degree $k$, it must have $k+1$ different vertices in its equivalence class, to ensure the total angle of $2\pi(k+1)$. Therefore, if we have $n$ vertices identified in our surface on the top (or bottom) row, it is a singularity of degree $n-1$.

\begin{figure}[h]
        \begin{center}
           \setlength{\unitlength}{200pt}
            \begin{picture}(1,.58)
                \put(0,0){\includegraphics[width=\unitlength]{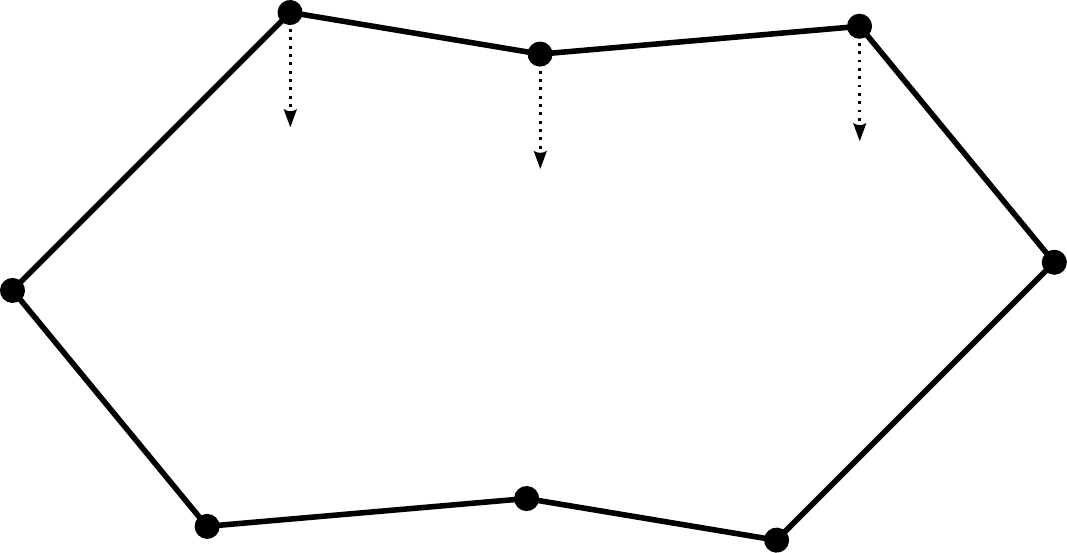}}

                \put(.12,.41){$1$}
                \put(.38,.52){$2$}
                \put(.64,.51){$3$}
                \put(.91,.41){$4$}

                \put(.07,.08){$4$}
                \put(.35,-.02){$3$}
                \put(.60,-.03){$2$}
                \put(.87,.09){$1$}
            \end{picture}
       \end{center}
    \caption{A suspension for $(4,3,2,1)$ has one singularity of degree $2$.}\label{fig.sing4321}
    \end{figure}
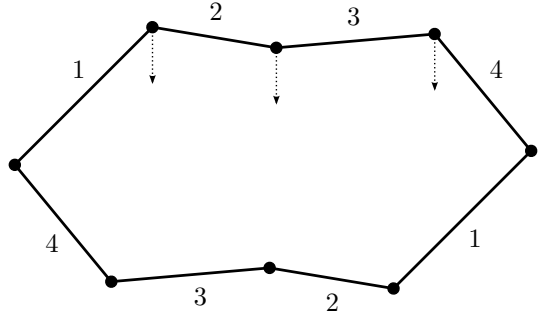

    \begin{exam}
         A suspension of $(4,3,2,1)$ has one singularity in Figure \ref{fig.sing4321}, which has $3$ copies on the top row. So it is a singularity of degree $2$.
    \end{exam}
    Suppose our surface has $m$ singularities of degrees $\ell_1,\dots,\ell_m$. If $s=\dusum{i=1}{m} \ell_i$, the number of edges in our surface is $d=m+s+1$. So the Euler characteristic is $\chi(S)=m-(m+s+1)+1=-s$. The genus of the surface is then $g(S) = \frac{2-\chi(S)}{2} = 1+\frac{s}{2}$. The number and degrees of singularities do not depend on our choice of $\lambda$ or $\tau$, only on $\pi$. Therefore the genus of $\pi$, $g(\pi)$, is well defined.\\
			
\begin{figure}
        \begin{center}
           \setlength{\unitlength}{250pt}
            \begin{picture}(1,.57)
                \put(0,0){\includegraphics[width=\unitlength]{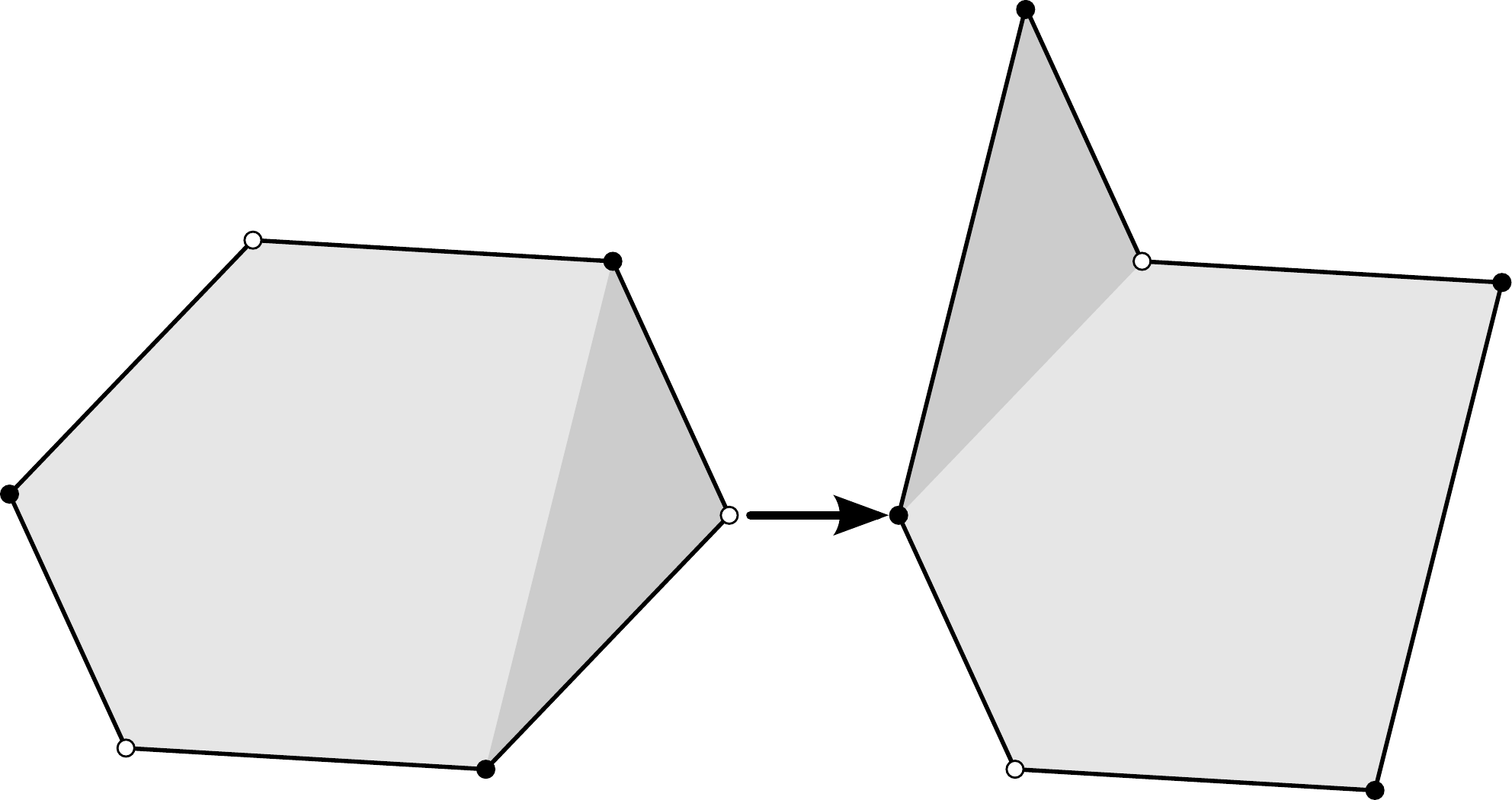}}

                \put(.22,.18){$S$}

                \put(.78,.18){$S'$}

                \put(.06,.30){$a$}
                \put(.27,.38){$b$}
                \put(.46,.28){$c$}

                \put(.01,.10){$c$}
                \put(.19,-.02){$b$}
                \put(.41,.08){$a$}

                \put(.60,.09){$c$}
                \put(.78,-.03){$b$}
                \put(.97,.16){$a$}

                \put(.60,.36){$a$}
                \put(.72,.46){$c$}
                \put(.86,.37){$b$}

                \put(.52,.22){$\rtt$}
            \end{picture}
       \end{center}
    \caption{A move of Rauzy-Veech induction on suspension $S$ for $\pi=(3,2,1)$.}\label{fig.suspensionRV}
    \end{figure}
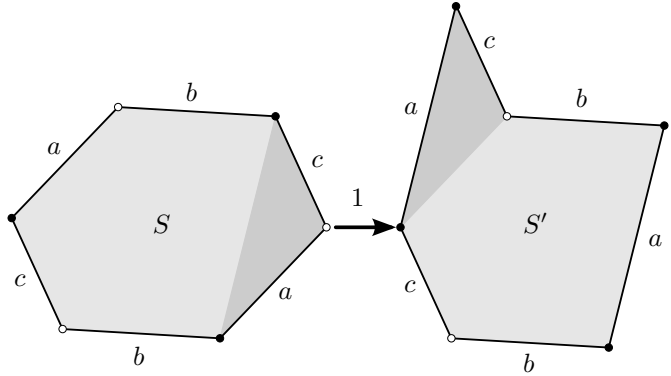

    Rauzy induction may be extended to these surfaces as well and is called \term{Rauzy-Veech (R-V) induction}. Let $\eps\in\{0,1\}$ be such that $\lambda_{\alpha_\eps}>\lambda_{\alpha_{1-\eps}}$. We can define a new surface by $S'=S(\pi',\lambda',\tau')$ where $(\pi',\lambda')$ is defined as in Section \ref{sec.RC} and $\tau'$ is defined as
	\begin{equation}\label{eq.Tpi_RV}
         \tau_\alpha':=\RHScase{
				\tau_\alpha, & \mbox{if }\alpha\neq\alpha_\eps,\\
				\tau_\alpha-\tau_{\alpha_{1-\eps}}, & \mbox{if }\alpha=\alpha_\eps.}
    \end{equation}
    This procedure is a ``cut and paste" by translation from $S=S(\pi,\lambda,\tau)$ to $S':=S(\pi',\lambda',\tau')$, as shown for $\pi=(3,2,1)$ and induction type $\rtt$ in Figure \ref{fig.suspensionRV}.
    Note that, as opposed to the case of $(4,1,3,2)$, $\pi=(3,2,1)$ has two singularities of degree zero. The induced permutation, $\rt\pi$ (see Definition \ref{def.RV}), has the same number and degrees of singularities as $\pi$. This is a general fact.
	
	\begin{prop}
		The number and degrees of singularities, and consequently the genus, are constant over a Rauzy Class.
	\end{prop}
	\begin{proof}
        This follows from counting before and after each type of inductive move to verify that the number and degrees of singularities do not change.
	\end{proof}
	
    While some singularities may be permuted by R-V induction, it is clear that the leftmost singularity remains fixed in the entire class. We shall call this singularity the \term{marked singularity}.
	
\begin{figure}[h]
        \begin{center}
           \setlength{\unitlength}{250pt}
            \begin{picture}(1,.3)
                \put(0,0){\includegraphics[width=\unitlength]{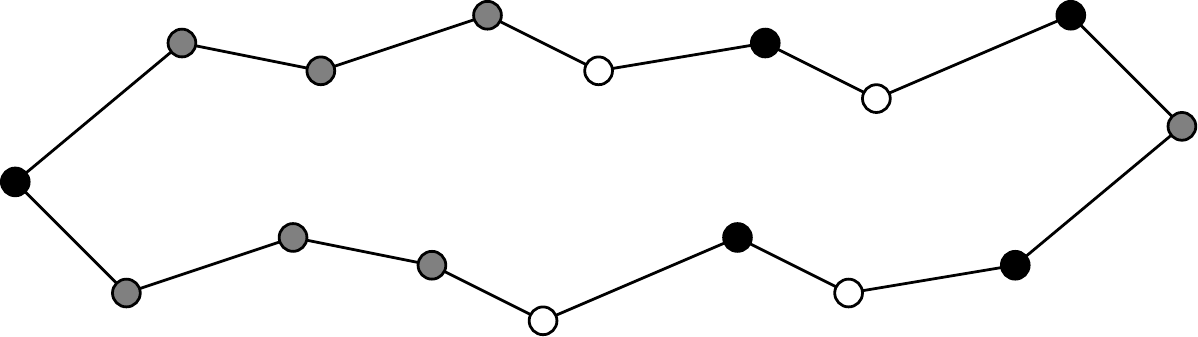}}

                \put(.065,.22){$1$}
                \put(.2,.25){$2$}
                \put(.325,.255){$3$}
                \put(.445,.26){$4$}
                \put(.55,.25){$5$}
                \put(.68,.245){$6$}
                \put(.8,.25){$7$}
                \put(.95,.235){$8$}

                \put(.015,.035){$8$}
                \put(.16,.01){$3$}
                \put(.285,.015){$2$}
                \put(.38,-.005){$4$}
                \put(.525,-.005){$7$}
                \put(.64,0){$6$}
                \put(.765,-.005){$5$}
                \put(.94,.065){$1$}
            \end{picture}
       \end{center}
       \caption{$\pi=(8,3,2,4,7,6,5,1)$ has signature $\sig(\pi)=(1,2,1)=(1,1,2)$.}\label{fig._83247651}
    \end{figure}
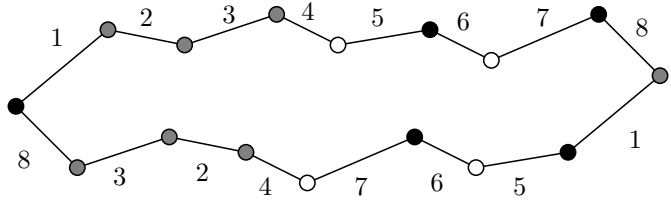

	\begin{defn}\label{def.signature}
        For a Rauzy Class $\RClass$, let $\ell_i$ denote the degrees of the $m$ singularities of $\pi$ with repetition. The $m$-tuple $(\ell_1,\dots,\ell_m)$, where $\ell_1$ is the degree of the marked singularity, is the \term{singularity signature} (or signature) of $\RClass$, denoted as $\sig=\sig(\RClass)$. If $\pi\in\RClass$, then $\sig(\pi)=\sig(\RClass)$.
	\end{defn}
    While the choice of $\ell_1$ is clear in Definition \ref{def.signature}, the other $\ell_i$'s may be in any order we wish. For example, the signature for $\pi=(8,3,2,4,7,6,5,1)$ can be written as $(1,1,2)$ or $(1,2,1)$ (see Figure \ref{fig._83247651}).\\

    Let $M_d:=\{S(\pi,\lambda,\tau):\pi=(\pi_0,\pi_1)\in\irr_d, \lambda\in\RR_+^\AAA, \tau\in\TTT_\pi, \mathrm{area}(S)=1\}$ minus the zero measure set where R-V induction is not well defined for all forward and backward iterates. Here the natural measure, $\mu$, is the product measure on $\irr_d$, $\RR_+^\AAA$, and $\TTT_\pi$, where the first is a counting measure and the last two inherit Lebesgue measure from $\RR^\AAA$. Let $R$ denote the action of R-V induction on each $S\in M_d$.

\begin{figure}[h]
        \begin{center}
           \setlength{\unitlength}{300pt}
            \begin{picture}(1,.5)
                \put(0,0){\includegraphics[width=\unitlength]{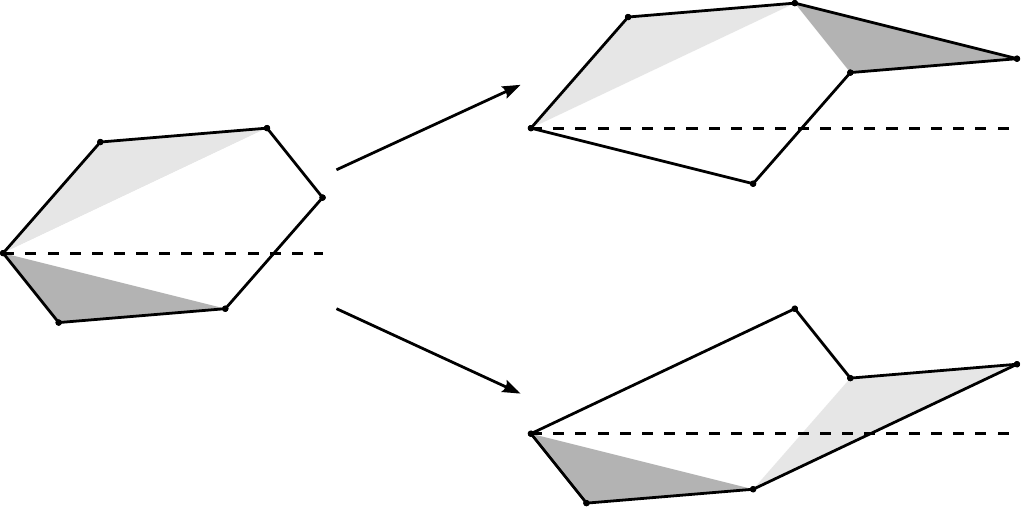}}
                \put(.15,.275){$S$}
                \put(.675,.4055){$S_\rtb$}
                \put(.695,.1){$S_\rtt$}
                \put(.37,.39){``$\rb^{-1}$"}
                \put(.36,.09){``$\rt^{-1}$"}
            \end{picture}
        \end{center}
        \caption{A suspension $S$ with $|\tau|>0$. $S_\rtb$ is not a valid suspension.}\label{fig.RV_IS_1TO1}
    \end{figure}
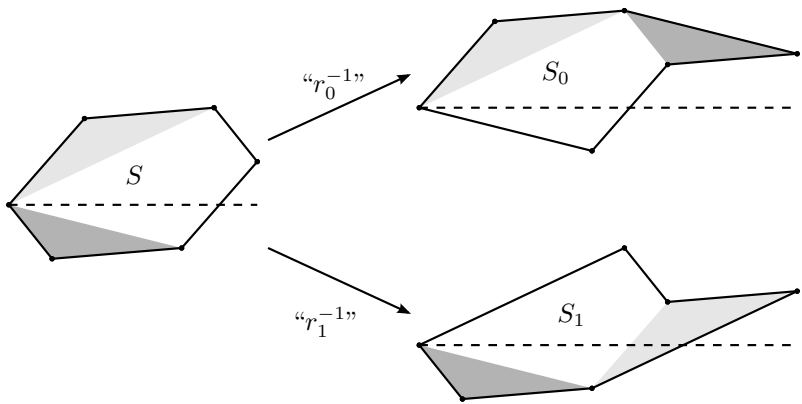

    \begin{rem}\label{rem.RV_is_now_1to1}
        Let $S=S(\pi,\lambda,\tau)$ and suppose $|\tau|=\dsum{\alpha\in\AAA}\tau_\alpha > 0$. Recall that $\alpha_\eps = \inv\pi_\eps(d)$. Consider $\pi^\rtt$ and $\lambda^\rtt$ from Remark \ref{rem.RV_is_2to1}. We see that $RS_\rtt=S$ for $S_\rtt=(\pi^\rtt,\lambda^\rtt,\tau^\rtt)$ where
            $$ \tau^\rtt_\alpha = \RHScase{\tau_{\alpha_0}+\tau_{\alpha_1}, & \mbox{if }\alpha = \alpha_\rtt,\\ \tau_\alpha, & \mathrm{otherwise.}}$$
        In this case the induction is type $\rtt$. Now let's attempt to construct $S_\rtb=S(\pi^\rtb,\lambda^\rtb,\tau^\rtb)$ such that $RS_\rtb=S$ by inductive move of type $\rtb$. We have $\pi^\rtb$ and $\lambda^\rtb$ as before. However Equation \eqref{eq.Tpi_RV} would require $\tau^\rtb$ to be defined by
            $$ \tau^\rtb_\alpha = \RHScase{\tau_{\alpha_0}+\tau_{\alpha_1}, & \mbox{if }\alpha = \alpha_\rtb,\\ \tau_\alpha, & \mathrm{otherwise.}}$$
        But then
            $$ \dsum{\alpha:\pi_\rtt(\alpha)\leq d-1}\tau^\rtb_\alpha = \dsum{\alpha:\alpha\neq \alpha_0,\pi_\rtt(\alpha)\leq d-1}\tau_\alpha+(\tau_{\alpha_0}+\tau_{\alpha_1})= |\tau| > 0.$$
        By Equation \eqref{eq.Tpi}, it follows that $\tau^\rtb \notin \TTT_{\pi^\rtb}$ (see Figure \ref{fig.RV_IS_1TO1}). $S_\rtt$ is therefore the unique suspension such that $RS_\rtt=S$. If instead $|\tau|<0$, we can similarly show that $S_\rtb$ exists while $S_\rtt$ does not. So we see that, as opposed to Rauzy induction on \IET's (see Remark \ref{rem.RV_is_2to1}), R-V induction is almost everywhere 1 to 1 on the space of suspensions.
    \end{rem}

    Let the map $F_t:M_d\rightarrow M_d$ be the flow defined by
	   $$ F_t(\pi,\lambda,\tau)=(\pi,e^t\lambda,e^{-t}\tau)$$
    Denote by $M_d^0$ the quotient space of $M_d$ under the equivalence $S\sim RS$ for $S\in M_d$. Then a fundamental domain for $M_d^0$ is
        $$\{ S(\pi,\lambda,\tau)\in M_d: 1\leq|\lambda|\leq|\lambda'|^{-1}\}.$$
    The flow $F_t$ is well defined on $M_d^0$ as $F_tR=RF_t$. There exists an $F_t$-invariant probability measure on $M_d^0$ which is absolutely continuous with respect to $\mu$ (Veech \cite{c.Ve82}), and $F_t$ is ergodic on this space with respect to this measure. This action $F_t$ on $M_d^0$ is called the \term{Teichm\"uller flow}. So each $M_d^0$ has a natural mapping into the moduli space of Abelian differentials. Denote by $\HHH(\ell_1,\dots, \ell_m)$ the stratum of differentials with $m$ zeros of degrees $\ell_1,\dots,\ell_m$. As opposed to Rauzy classes, the ordering of the $\ell_i$'s is completely arbitrary in terms of the strata of Abelian differentials.

\sect{Classification of True Rauzy Classes}\label{sec.classify}

    Each stratum, $\HHH(\ell_1,\dots, \ell_m)$, can generally be divided further into connected components, which correspond to \term{Extended Rauzy Classes} (see \cite{c.Ve90}). The following theorems completely categorize every connected component for all strata. A stratum is \term{hyperelliptic} if a Riemann surface with differential in the stratum is hyperelliptic (see Section \ref{sec.hyperelliptic}). A stratum with all singularities of even degree has a flow invariant $\ZZ_2$-valued property called the \term{parity} of its \term{spin structure}. Details on this and calculations will be presented in Section \ref{sec.spin}.\\

    If the genus of $\RClass\subseteq\irr_d$ is $1$, we conclude from Sections \ref{sec.hyperelliptic} and \ref{sec.removable} that
    $$\pi=(d,2,\dots,d-1,1)$$
    belongs to $\RClass$. The following theorem categorizes all strata of genus $2$ and $3$.
	
	\begin{nonum}(M. Kontsevich and A. Zorich \cite{c.KZ})
		The moduli space of Abelian differentials on a complex curve of genus $g=2$ contains two strata: $\HHH(1,1)$ and
		$\HHH(2)$. Each of them is connected and hyperelliptic.\\
		
		Each stratum $\HHH(2,2)$, $\HHH(4)$ of the moduli space of Abelian differentials on a complex curve of genus $g=3$
		has two connected components: the hyperelliptic one, and one having odd spin structure. The other strata are
		connected for genus $g=3$.
	\end{nonum}
	
	The following theorem categorizes the connected components for each stratum of genera $4$ or greater.
	
	\begin{nonum}(M. Kontsevich and A. Zorich \cite{c.KZ})
		All connected components of any stratum of Abelian differentials on a complex curve of genus $g\geq 4$ are
		described by the following list:
		\begin{itemize}
            \item The stratum $\HHH(2g-2)$ has three connected components: the hyperelliptic one, $\Hhyp(2g-2)$, and components $\Heven(2g-2)$ and $\Hodd(2g-2)$ corresponding to even and odd spin structures.
            \item The stratum $\HHH(2\ell,2\ell)$, $\ell\geq 2$ has three connected components: $\Hhyp(2\ell,2\ell)$, $\Heven(2\ell,2\ell)$ and $\Hodd(2\ell,2\ell)$.
            \item All the other strata of the form $\HHH(2\ell_1,\dots,2\ell_m)$, where all $\ell_i\geq 1$, have two connected components: $\Heven(2\ell_1,\dots,2\ell_m)$ and $\Hodd(2\ell_1,\dots,2\ell_m)$.
			\item The strata $\HHH(2\ell-1,2\ell-1)$, $\ell\geq 2$, have two connected components; one of them,
				$\Hhyp(2\ell-1,2\ell-1)$, is hyperelliptic; the other one, $\Hnonhyp(2\ell-1,2\ell-1)$, is not.
			\item All other strata of Abelian differentials on complex curves of genera $g\geq 4$ are nonempty and
				connected.
		\end{itemize}
	\end{nonum}
	
    We are given a full classification of each connected component by the above results. To each connected component, we denote the \term{type} by the information other than the singularities. The type takes one of the following values \{-, even, odd, hyperelliptic, nonhyperelliptic\} as applicable. This however is not enough to calculate what Rauzy Class a permutation $\pi\in\irr$ belongs to, only the Extended Rauzy Class.

    \begin{exam}
        Consider
        $$ \pi= \cmtrx{\LL{a}{z}~\LL{b}{c}~\LL{c}{b}~\LL{d}{d}~\LL{e}{w}~\LL{f}{f}~\LL{w}{e}~\LL{x}{y}~\LL{y}{x}~\LL{z}{a}}$$
        and
        $$ \pi'=\cmtrx{\LL{a}{z}~\LL{b}{c}~\LL{c}{b}~\LL{d}{d}~\LL{e}{f}~\LL{f}{e}~\LL{w}{y}~\LL{x}{x}~\LL{y}{w}~\LL{z}{a}}.$$
        Both $\pi$ and $\pi'$ have three singularities one each of degrees $1$,$2$ and $3$. Therefore both $\pi$ and $\pi'$ belong to the stratum $\HHH(1,2,3)$. However the marked singularity of $\pi$ is of degree $3$, while the marked singularity of $\pi'$ is of degree $1$. Because the degree of the marked singularity is fixed throughout a Rauzy Class, $\RClass(\pi)\neq\RClass(\pi')$.
    \end{exam}

    It becomes clear that in order to distinguish Rauzy Classes, the degree of the marked singularity must be considered. Indeed, the following theorem shows the addition of this final invariant completes the classification of all Rauzy Classes:
	
	\begin{nonum}(C. Boissy \cite{c.B09})
		$\pi_1,\pi_2\in\irr_d$ belong to the same Rauzy class if and only if they belong to the same connected component
		and their marked singularities are the same degree.
	\end{nonum}
	
	We restate the above information in a different form and make an observation that, while clear from
	everything above, is crucial to our main result.
	
	\begin{corollary}\label{cor.classify}
		Every Rauzy class is uniquely determined by signature and type. So given Rauzy class $\RClass$,
		if $\pi\in\irr$
		has the same signature and type as $\RClass$, then necessarily $\pi\in\RClass$.
	\end{corollary}
	
\sect{Hyperelliptic Surfaces}\label{sec.hyperelliptic}

    \begin{defn}\label{def.hyp_perm}
        A surface with quadratic differential $(M,q)$ of genus $g$ is \term{hyperelliptic} if there exists a map $h:M\to M$ such that
        \begin{itemize}
            \item $h=\inv{h}$,
            \item $h_*q = -q$,
            \item $h$ fixes $2g+2$ points,
        \end{itemize}
         and such an $h$ is called a \term{hyperelliptic involution}. A permutation $\pi$ is \term{hyperelliptic} if every suspension of $\pi$ is hyperelliptic.
    \end{defn}

    \begin{rem}
        Because $h$ is a well defined map on the differential, $h$ must take singularities to singularities. Also, $h$ must take geodesics to geodesics. Therefore $h$ maps saddle connections (geodesics with endpoints that are singularities) to saddle connections. In this case, we do not consider removable singularities.
    \end{rem}

    \begin{rem}\label{rem.how_to_check_hyp}
        Consider $\pi\in\irr_d$. Any given suspension $S$ is represented by a polygon in $\CC$ whose differential is represented by the standard $dz$ in its interior. In this case, the only possible candidates for a hyperelliptic involution on $S$ are of the form $z\mapsto -z + c$ for some constant $c\in\CC$. These maps automatically satisfy the first two conditions in Definition \ref{def.hyp_perm}.
    \end{rem}

    \begin{defn}
        For $d\geq 2$, let $\pi_{(d)}$ be the permutation such that $\pi_{(d)}(i)=d-i+1$ for all $i\in\dset{d}$.
    \end{defn}

    \begin{lem}\label{lem.pi_d_is_hyp}
        $\pi_{(d)}$ is hyperelliptic.
    \end{lem}

    \begin{proof}
        Consider any suspension $S=S(\pi_{(d)},\lambda,\tau)$. We will construct $h$ and show that it satisfies Definition \ref{def.hyp_perm}. Let $h(z) = -z + |\lambda|+\imath|\tau|$. The first two conditions are satisfied. In order to show that $h$ is the appropriate map, we will define the vertices $p^\eps_k$ for $k\in\{0,\dots,d\}$ and $\eps\in\{0,1\}$ by
        $$ \begin{array}{rcl}
            p^0_k &=&\dusum{j=1}{k}\lambda_j+\imath\dusum{j=1}{k}\tau_j, \mbox{ and}\\
            p^1_k &=&\dusum{j=1}{k}\lambda_{d-j+1}+\imath\dusum{j=1}{k}\tau_{d-j+1}.
            \end{array}$$
        We note that the top segment labeled $k$ in $S$ has endpoints $p^0_{k-1}$ and $p^0_k$ while the bottom segment labeled $k$ in $S$ has endpoints $p^1_{d-k}$ and $p^1_{d-k+1}$. Because $h$ is an isometry, it maps segments to segments. We examine mapping the endpoints $\pi_k^\eps$ under $h$. For $k\in\{0,\dots,d\}$,
        $$ \begin{array}{rcl}
                h(p^0_k) & = & -p^0_k + |\lambda|+\imath|\tau| \\
                    & = & \dusum{j=1}{d}\lambda_j - \dusum{\ell=1}{k}\lambda_\ell + \imath\left(\dusum{j=1}{d}\tau_j - \dusum{\ell=1}{k}\tau_\ell\right)\\
                    & = & \dusum{j=k+1}{d}\lambda_j = \imath\dusum{j=k+1}{d}\tau_j \\
                    & = & \dusum{j'=1}{d-k}\lambda_{d-j'+1} + \imath\dusum{j'=1}{d-k}\tau_{d-j'+1}\\
                    & = & p^1_{d-k}.
            \end{array}$$
         Because $h=\inv h$, we conclude that any segment labeled $k$ is mapped to the other segment labeled $k$. As these segments are identified, $h$ fixes these segments. Now it remains to count the fixed points.\\

        If $d=2m$ is even, there is one singularity of degree $2m-2$, the genus is $m$ and there should be $2m+2 = d+2$ fixed points. There are $d$ segments each with a fixed midpoint. The point $\frac{1}{2}(|\lambda|+\imath|\tau|)$ is fixed, and the singularity represented by the class of all $p^\eps_k$'s is fixed. Therefore $h$ fixes $2g+2$ points.\\

        If $d=2m+1$ is odd, there are two singularities each of degree $m-1$, the genus is $m$ and there should be $2m+2=d+1$ fixed points. We note that this time, the two singularities, one represented by all $p_k^\eps$'s with even $k$'s and the other by all odd $k$'s, are interchanged by $h$. However $h$ fixes the $d$ midpoints of the labeled segments and the point $\frac{1}{2}(|\lambda|+\imath|\tau|)$. Therefore $h$ fixes $2g+2$ points.\\

        So we see that in either case, $h$ is the hyperelliptic involution for $S$.
    \end{proof}

    \begin{prop}\label{prop.hyp}
        Let $\pi\in\irr_d$ be standard with no removable singularities. If $\pi$ is hyperelliptic, then $\pi=\pi_{(d)}$.
    \end{prop}

    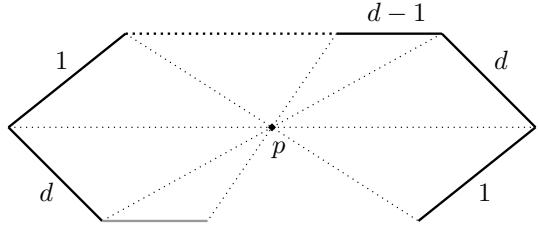
\begin{figure}[t]
        \begin{center}
           \setlength{\unitlength}{200pt}
            \begin{picture}(1,.42)
                \put(0,0){\includegraphics[width=\unitlength]{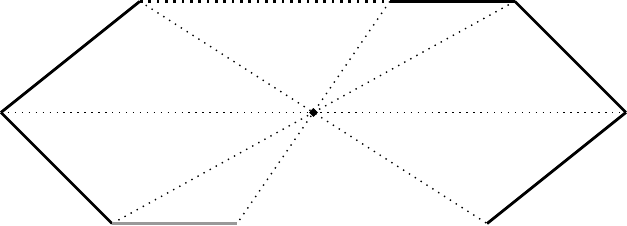}}

                \put(.09,.29){$1$}
                \put(.68,.38){$d-1$}
                \put(.92,.29){$d$}

                \put(.50,.13){$p$}

                \put(.06,.04){$d$}
                \put(.89,.04){$1$}
            \end{picture}
       \end{center}
       \caption{If the hyperelliptic involution must fix the point $p$, then the grey segment must be identified with $d-1$.}\label{fig.std_hyp}
    \end{figure}

    \begin{proof}
        Let $\pi$ be a hyperelliptic standard permutation. We will explicitly construct a suspension for $\pi$ that excludes any possibility but $\pi=\pi_{(d)}.$ Assume $\AAA=\dset{d}$ and $\pi=(\pi_0,\pi_1)$ where $\pi_0(i)=i$. Fix any $m\in(0,\frac{1}{4})$ and let $\lambda\in\RR_+^d$ be defined by $\lambda_i = 1+m^i$ for $i\in\dset{d}$. Also, let $\tau$ be defined by $\tau = (1,0,\dots,0,-1)$. From Equation \eqref{eq.Tpi}, we see that $S=S(\pi,\lambda,\tau)$ is a valid suspension for $\pi$. Now consider the hyperelliptic involution for $S$, denoted as $h$. By construction, the segments labeled by $1$ and $d$ must be interchanged under $h$, as no other saddle connections exist of the appropriate length. We conclude that $h(z) = -z + |\lambda|+\imath|\tau| = |\lambda|-z$ and see that the top segment labeled $j$ must be mapped to the bottom segment labeled $j$, as no other saddle connection would have the appropriate length. We then show iteratively that $\pi(j) = d - j + 1$ as desired (see Figure \ref{fig.std_hyp}). By Lemma \ref{lem.pi_d_is_hyp}, this is hyperelliptic.
    \end{proof}

\sect{Calculation of Spin Parity}\label{sec.spin}

	The results in this section follow from Appendix C in \cite{c.Z08}. We refer the reader to that paper for
	details.\\
	
    To each $\pi\in\irr_d$ with all singularities of even degree, we can define the parity of the spin structure of the corresponding suspension surface $S$. To do so, we must find a \term{symplectic basis} of $H_1(S)$. This is a choice of closed cycles $\alpha_1,\beta_1,\dots,\alpha_g,\beta_g\in H_1(S)$, $g=g(\pi)$, with the following conditions: $\alpha_i\inx\alpha_j=\beta_i\inx\beta_j=0$ and $\alpha_i\inx\beta_j=\delta_{ij}$ where $\alpha\inx\beta$ is the algebraic intersection number. For a loop $\gamma$, the Gauss map is the lift of $\gamma$ to the unit tangent bundle, a map from $H_1(S)\rightarrow \sS^1$ (where $\sS^1\subset\RR^2$ is the unit circle), and let $\mathrm{ind}(\gamma)$ be the degree of the Gauss map. The spin parity of the surface can be calculated by:
	\begin{equation}\label{eq.prespin}
	{
		\Phi(S) := \sum_{i=1}^{g}(\mathrm{ind}(\alpha_i)+1)(\mathrm{ind}(\beta_i)+1)~~(mod~2).
	}
	\end{equation}

    Let $\phi(\gamma):=\mathrm{ind}(\gamma)+1$. This value is independent of choice of suspension surface $S$. Therefore we may instead speak of the parity of $\pi$ itself. Using these conventions, the previous equation becomes
	\begin{equation}\label{eq.spin}
	{
		\Phi(\pi) :=\sum_{i=1}^{g}\phi(\alpha_i)\phi(\beta_i)~~(mod~2).
	}
	\end{equation}
	
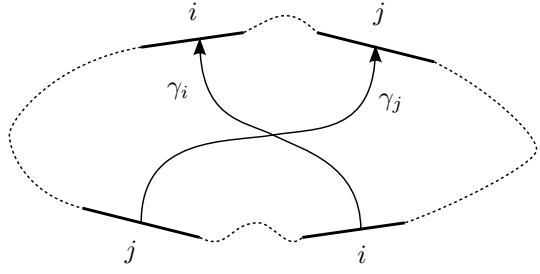
\begin{figure}
        \begin{center}
           \setlength{\unitlength}{200pt}
            \begin{picture}(1,.46)
                \put(0,0){\includegraphics[width=\unitlength]{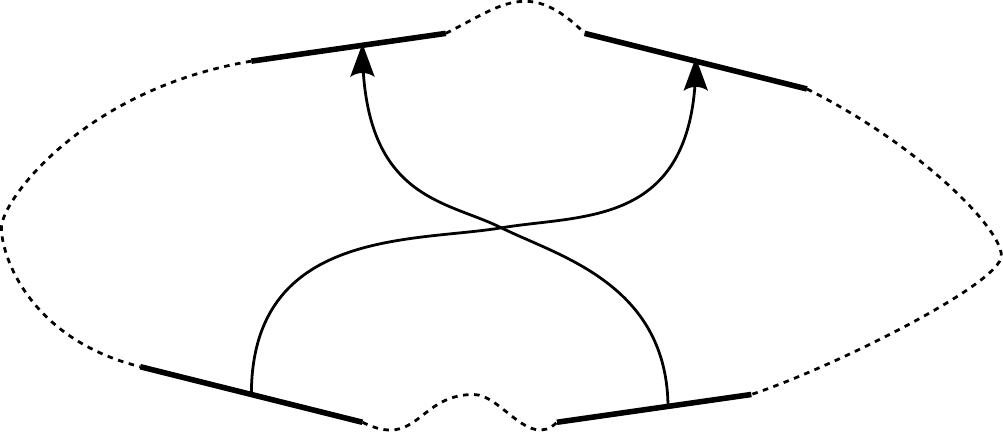}}

                \put(.34,.42){$i$}
                \put(.69,.42){$j$}

                \put(.22,-.03){$j$}
                \put(.66,-.04){$i$}

                \put(.30,.28){$\gamma_i$}
                \put(.70,.26){$\gamma_j$}

            \end{picture}
       \end{center}
        \caption{Loops $\gamma_i$ and $\gamma_j$ intersect in suspension $S$.}\label{fig.gamma}
    \end{figure}

    For a surface $S=S(\pi)$, we will define for each $i\in\dset{d}$ a loop $\gamma_i$. Start with any point on the embedded subinterval $I_i$ in $I_S$. The loop will move in the positive vertical direction until it returns to $I_S$. Then close the loop by a horizontal line. Now deform the loop continuously so it becomes smooth and everywhere transverse to the horizontal direction. Call this loop $\gamma_i$. Let $c_i=[\gamma_i]$ be the cycle representative of $\gamma_i$ in $H_1(S)$. See Figure \ref{fig.gamma}.
	
    It is clear that $\mathrm{ind}(\gamma_i)=0$ as $\gamma_i$ is always transverse to the horizontal direction. Therefore $\phi(c_i)=1$. From the definition of $\gamma_i$ and $\Omega=\Omega_\pi$ (from Section \ref{sec.IET}),
	\begin{equation}\label{eq.inx}
	{
		c_i\inx c_j = \Omega_{i,j}
	}
	\end{equation}
    and that the span of the $c_i$'s is $H_1(S)$. Because the above calculations \eqref{eq.prespin} and \eqref{eq.spin} are $(mod~2)$, the following calculations are over $\ZZ_2$. Note that now $\Omega$ is a symmetric matrix of zeros and ones. We may still keep the definition $\phi(c_i):H_1(S)\rightarrow\ZZ_2$. It is a well defined quadratic form on the intersection and has the following relationship as a direct result from \cite{c.J80}: for $c,c'\in H_1(S)$,
	\begin{equation}\label{eq.quadform}
		\phi(c+c')=\phi(c)+\phi(c')+c\inx c'.
	\end{equation}
	We recall the following relationship for $a,b,c\in H_1(S)$ on the intersection number:
	\begin{equation}\label{eq.addinx}
		(a+b)\inx c = a\inx c + b\inx c.
	\end{equation}
    We now describe the iterative process to choose our symplectic basis from the $c_i$'s. First let $\alpha_1:=c_1$. Let $\beta_1:=c_j$ for some $j$ such that $\Omega_{1j}=1$. We adjust each $c_i$, $i=2,3,\dots,j-1,j+1,\dots,d$ by the following rule: the remaining vectors must be adjusted so that they have trivial intersection number with $c_1$ and $c_j$. So we consider $c_i':=c_i+\eps_1c_1+\eps_jc_j$. Then $c_i'\inx c_1=0\Rightarrow\eps_j=c_i\inx c_1=\Omega_{1,i}$, and $c_i'\inx c_j=0\Rightarrow\eps_1=c_i\inx c_j=\Omega_{i,j}$. Now we use \eqref{eq.quadform} and \eqref{eq.inx} to calculate
	$$ \arry{rcl}{
			\phi(c_i') &=& \phi(c_i+\Omega_{i,j}c_1+\Omega_{1,i}c_j)\\
				&=& \phi(c_i)+\phi(\Omega_{i,j}c_1+\Omega_{1,i}c_j)\\
				&=& \phi(c_i)+\Omega_{i,j}\phi(c_1) + \Omega_{1,i}\phi(c_j) + \Omega_{i,j}\Omega_{1,i}.}$$
	And using \eqref{eq.addinx}, for $i,k\in\{2,3,\dots,j-1,j+1,\dots,d\}$,
	$$ \arry{rcl}{
			c_i'\inx c_k' &=& (c_i+\Omega_{i,j}c_1+\Omega_{1,i}c_j)\inx(c_k+\Omega_{k,j}c_1+\Omega_{1,k}c_j)\\
				&=& c_i\inx(c_k+\Omega_{k,j}c_1+\Omega_{1,k}c_j)+
					\Omega_{i,j}c_1\inx(c_k+\Omega_{k,j}c_1+\Omega_{1,k}c_j)\\
				& &	+ \Omega_{1,i}c_j\inx(c_k+\Omega_{k,j}c_1+\Omega_{1,k}c_j)\\
				&=& \Omega_{i,k}+\Omega_{k,j}\Omega_{1,i}+\Omega_{1,k}\Omega_{i,j}
					+\Omega_{i,j}\Omega_{1,k}+\Omega_{i,j}\Omega_{1,k}
					+\Omega_{1,i}\Omega_{j,k}+\Omega_{1,i}\Omega_{k,j}\\
				&=& \Omega_{i,k}+\Omega_{k,j}\Omega_{1,i}+\Omega_{1,k}\Omega_{i,j}.}$$
	So we restate these results together for reference in later calculations,
    \begin{subequations}\label{eq.iterate}
	   \begin{align}
		  c_i'&:=c_i+\Omega_{i,j}c_1+\Omega_{i,1}c_j,\\
    		\phi(c_i')&:=\phi(c_i)+\Omega_{i,j}\phi(c_1)+\Omega_{i,1}\phi(c_j)+\Omega_{i,1}\Omega_{i,j},\\
    		c_i'\inx c_k'&:=\Omega_{i,k}+\Omega_{i,1}\Omega_{k,j}+\Omega_{i,j}\Omega_{k,1}.
	   \end{align}
	\end{subequations}
	
    We now have a new set of remaining cycles $c_i'$ with intersection matrix defined by $\Omega_{i,k}'=c_i'\inx c_k'$. We then pick a pair of intersecting cycles and name them $\alpha_2$ and $\beta_2$. We then alter the remaining cycles again by Equations \eqref{eq.iterate}. This process terminates when all pairs $\alpha_1,\beta_1,\dots,\alpha_g,\beta_g$ are chosen. Now we can calculate the parity by \eqref{eq.spin}.

    \begin{exam}
        We will calculate the spin parity of $\pi=(4,3,6,1,5,2)$. This is not hyperelliptic and has one singularity of degree $4$. So we first consider the initial conditions,
        $$ \Omega = \tbl{|r|cccccc|}{
            \hline
            & 1&2&3&4&5&6\\
            \hline
            1& 0&0&1&1&0&1\\
            2& 0&0&1&1&1&1\\
            3& 1&1&0&1&0&0\\
            4& 1&1&1&0&0&0\\
            5& 0&1&0&0&0&1\\
            6& 1&1&0&0&1&0\\
            \hline } \mbox{ and } \phi(c_i)=1.$$
        We may choose initial basis pair $(\alpha_1,\beta_1)=(c_1,c_3)$ and $\phi(\alpha_1)=\phi(\beta_1)=1$. We use the Equations \ref{eq.iterate} to derive
        $$ \Omega' = \tbl{|r|cccc|}{
            \hline
            & 2&4&5&6\\
            \hline
            2& 0&0&1&0\\
            4& 0&0&0&1\\
            5& 1&0&0&1\\
            6& 0&1&1&0\\
            \hline }
            \mbox{ for  }
            \arry{|rcl|rcl|}{
                \hline
                c_2' & = & c_1+c_2      &\phi(c_2') & = & 0\\
                c_4' & = & c_1+c_3+c_4  &\phi(c_4') & = & 0\\
                c_5' & = & c_5          &\phi(c_5') & = & 1\\
                c_6' & = & c_3+c_6      &\phi(c_6') & = & 0\\
                \hline } $$
        From these remaining vectors, we choose $(\alpha_2,\beta_2) =(c_2',c_5')$ where $\phi(\alpha_2)=0$ and $\phi(\beta_2)=1$. We then modify the remaining vectors to derive
        $$  \Omega'' = \tbl{|r|cc|}{
                \hline
                    & 4&6 \\
                \hline
                    4& 0&1 \\
                    6& 1&0 \\
                \hline }
            \mbox{ for }
            \arry{|rcl|rcl|}{
                \hline
                c_4'' & = & c_4'        & \phi (c_4'') &=& 0 \\
                c_6'' & = & c_2'+c_6'   & \phi (c_6'') &=& 0 \\
                \hline }$$
        Our only remaining choice is $(\alpha_3,\beta_3)=(c_4'',c_6'')$ with $\phi(\alpha_3)=\phi(\beta_3)=0$. So by Equation \eqref{eq.spin},
        $$ \Phi(\pi) = \dusum{i=1}{3}\phi(\alpha_i)\phi(\beta_i) = 1.$$
    \end{exam}

\sect{Linear Involutions and Generalized Permutations}\label{sec.gen.perm}

    From the previous sections, we see that interval exchange transformations are realized as the first return to transversals of vertical foliations on flat surfaces. In \cite{c.DN90}, Danthony and Nogueira discuss linear involutions as first returns of vertical foliations on surfaces with non-trivial $\ZZ_2$ holonomy. Such surfaces are Riemannian surfaces with transition functions that are compositions of translations as well as order 2 rotations. These surfaces may be realized as surfaces with quadratic differentials that are not squares of Abelian differentials. Consider such surface $S$ with horizontal segment $I$. This time, an ``upward" vertical trajectory from $I$ may return to $I$ pointing in the ``downward" direction. As a result, we have have the following:

    \begin{defn} \label{def.gen.linear_involutions} Let $I=(0,L)$, $L>0$, $\hI = I \times \{0,1\}$ and $\Sigma_0,\Sigma_1\subset\hI$ be finite sets.
        A \term{linear involution} $T$ is a bijection $T:\hI\setminus\Sigma_0\rightarrow \hI\setminus\Sigma_1$ such that
        \begin{enumerate}
            \item $T=f\circ \tT$, where $f(x,\eps)=(x,1-\eps)$.
            \item $\tT$ is an involution, or $\inv\tT=\tT$.
            \item $\tT$ has no fixed point.
            \item $\tT$ is smooth.
            \item For $(x,\eps)\in \hI\setminus \Sigma_0$,
                $$\frac{\mathrm{d}}{\mathrm{d}x}\tT(x,\eps) = \RHScase{1,& \mbox{if }\tT(x,\eps)\in I \times\{1-\eps\}, \\ -1, & \mbox{if }\tT(x,\eps)\in I \times\{\eps\}}.$$
        \end{enumerate}
    \end{defn}

    From this definition, we see that a linear involution $T$ takes a finite set of open intervals and reorders them by isometry. However, unlike \IET's, some intervals are mapped by order reversing isometries.

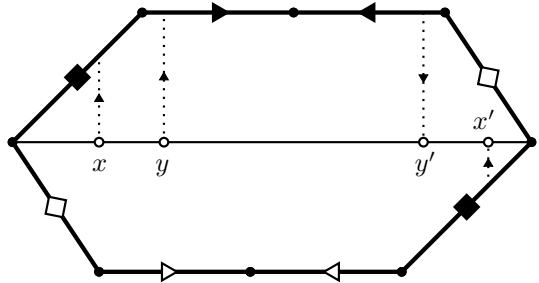
\begin{figure}[h]
        \begin{center}
           \setlength{\unitlength}{200pt}
            \begin{picture}(1,.54)
                \put(0,0){\includegraphics[width=\unitlength]{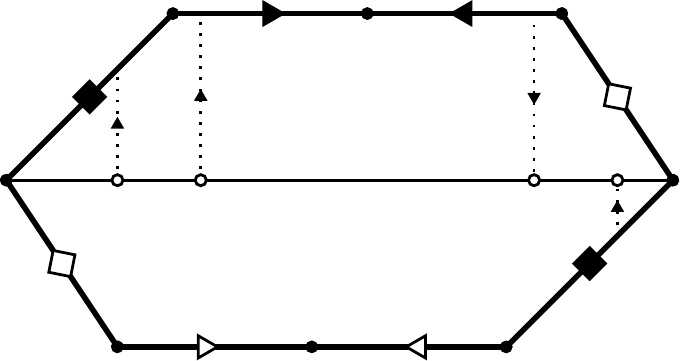}}

                \put(.16,.21){$x$}
                \put(.28,.21){$y$}
                \put(.77,.21){$y'$}
                \put(.88,.29){$x'$}

            \end{picture}
       \end{center}
       \caption{For this linear involution, $T(x,0)=(x',0)$ and $T(y,0)=(y',1)$.}\label{fig.Linear_Involution}
    \end{figure}

    \begin{rem}\label{rem.gen.interpret}
        The interpretation of $T$ can be viewed as acting on $\hI$ when viewed as a ``position" times ``direction" space. In other words, if $T(x,\eps)\in I\times \{\eps\}$, then the vertical trajectory leaving $I$ in $S$ will return to $I$ pointing in the same direction, while $T(x,\eps)\in I\times \{1-\eps\}$ would indicate that the direction would be opposite, as shown in Figure \ref{fig.Linear_Involution}.
    \end{rem}

    We will now encode a linear involution $T$ in a similar way to \IET's. Let $\AAA$ be an alphabet on $d$ letters as before.
    \begin{defn}\label{def.gen.perm}
        A \term{generalized permutation}, on $\AAA$, $\pi = (\pi,d_0,d_1)$ is a triple such that $d_0,d_1 \in \NN$, $d_0+d_1=2d$ and $\pi:\dset{2d}\rightarrow \AAA$ is a two-to-one map. A generalized permutation shall be denoted as
        $$ (\pi,d_0,d_1) = \cmtrx{\pi(1), \dots ,\pi(d_0) \\ \pi(d_0+1), \dots, \pi(2d)}.$$
        The set of all generalized permutations on $\AAA$ shall be denoted as $\genperm_\AAA$.
    \end{defn}

    \begin{rem}
        By the notation given in the previous definition, the values of $d_0$ and $d_1$ are understood. Therefore, we will assume the values $d_0$ and $d_1$ when possible, referring to the generalized permutation only as $\pi$.
    \end{rem}

    We consider $\pi\in\genperm_\AAA$ and $\pi'\in\genperm_{\AAA'}$ equivalent, denoted $\pi=\pi'$, if there exists a bijection $b:\AAA'\to\AAA$ such that $\pi\equiv b\circ\pi'$ and, naturally, $d_i=d_i'$ for $i\in\{0,1\}$.

    \begin{exam}\label{ex.gen.equiv}
        Let $\AAA=\{a,b,c,d\}$ and $\AAA'=\{1,2,3,4\}$. Then
            $$ \pi=\cmtrx{\LL{}{c} ~ \LL{a}{d} ~ \LL{b}{c} ~ \LL{b}{d} ~ \LL{}{a}}=\cmtrx{\LL{}{3}~\LL{1}{4}~\LL{2}{3}~\LL{2}{4}~\LL{}{1}}.$$
        In this case $d_0=3$, $d_1=5$.
    \end{exam}

    For generalized permutation $\pi$, we shall define $D_\pi:\AAA\to\RR$ as
    \begin{equation}
        D_\pi(\alpha) := (i-(d_0+\frac{1}{2}))(j-(d_0+\frac{1}{2}))
    \end{equation}
    for $\alpha\in\AAA$ and $\inv\pi(\alpha) = \{i,j\}$. The sign of $D_\pi$ has the following interpretation: if $D_\pi(\alpha)>0$ then both occurrences of $\alpha$ in $\pi$ appear on the same row, and if $D_\pi(\alpha)<0$ the symbol $\alpha$ appears once on each row. In the representation of $\pi$ over $\AAA$ in Example \ref{ex.gen.equiv}, $D_\pi(a)<0$ and $D_\pi(b),D_\pi(c),D_\pi(d)>0$.

    \begin{rem}
        If for each $\alpha\in\AAA$, $D_\pi(\alpha)<0$, then $\pi$ is a true permutation with $d_0=d_1=d$. Considering the notation in Section \ref{sec.IET}, we see that $\pi=(\pi_0,\pi_1)$ by $\pi_0 \equiv \inv{(\pi|_{\{1,\dots,d\}})}$ and $\pi_1 \equiv \inv{(\pi|_{\{d+1,\dots,2d\}})}-d$.
    \end{rem}

    \begin{exam}
        Consider the permutation
        $$ \pi=\cmtrx{\LL{a}{d}~\LL{b}{a}~\LL{c}{c}~\LL{d}{b}}.$$
        According to Definition \ref{def.yoccoz} (as a true permutation), $\pi=(\pi_0,\pi_1)$ where
        $$ \begin{array}{llll}
            \pi_0(a) = 1, & \pi_0(b)=2, & \pi_0(c)=3, & \pi_0(d)=4, \\
            \pi_1(a) = 2, & \pi_1(b)=4, & \pi_1(c)=3, & \pi_1(d)=1.
            \end{array}$$
        We also see that by Definition \ref{def.gen.perm}, $\pi=(\pi,4,4)$ where
        $$ \begin{array}{llll}
            \inv\pi(a) = \{1,6\}, & \inv\pi(b) = \{2,8\},& \inv\pi(c) = \{3,7\},& \inv\pi(d) = \{4,5\}.
            \end{array}$$
    \end{exam}

    Let $\pi\in \genperm_\AAA$ be given and $\lambda\in \RR^\AAA_+$ satisfy
    \begin{equation}\label{eq.gen.lengths}
        L=L(\lambda):=\sum_{k=1}^{d_0} \lambda_{\pi(k)} = \sum_{j=d_0+1}^{2d} \lambda_{\pi(j)}.
    \end{equation}
    Define $\RR^{\AAA,\pi}_+:=\{\lambda \in \RR^\AAA_+ :~ \eqref{eq.gen.lengths}\makebox{ is satisfied}\}$. Assume that $\RR^{\AAA,\pi}_+$ is non-empty (we will address this assumption later). For $\lambda\in\RR_+^{\AAA,\pi}$, we shall define a linear involution $T$ (see Definition \ref{def.gen.linear_involutions}) associated to $(\pi,\lambda)$: Let $L$ be as in \eqref{eq.gen.lengths}. Let $I = (0,L)$ and $\hI=I\times\{0,1\}$. Define the points
    $$ p_i = \RHScase{ (\dsum{1\leq j < i} \lambda_{\pi(j)},0), & 1\leq i \leq d_0 \\ (\dsum{d_0< j < i} \lambda_{\pi(j)},1), & d_0<i\leq 2d}$$
    for $i\in\dset{2d}$. Then $\Sigma_0:=\underset{i\neq 1,d_0+1}\bigcup \{p_i\}$, and $\Sigma_1 = f(\Sigma_0)$ where $f(x,\eps)=(x,1-\eps),$ as in Definition \ref{def.gen.linear_involutions}. We consider two cases. First let $i=i(x)\in\dset{d_0}$ be defined for $x\in I$ by
        $$ \dusum{k=1}{i-1} \lambda_{\pi(k)}\leq x<\dusum{k=1}{i}\lambda_{\pi(k)}$$
    and $\delta(x) = x - \dusum{k=1}{i-1}\lambda_{\pi(k)}$.
    In this case
        $$ \tT(x,0) = \RHScase{([\dusum{k=1}{\ti}\lambda_{\pi(k)}]-\delta(x),0), & \mbox{if }\ti \leq d_0,\\ ([\dusum{k=d_0+1}{\ti-1}\lambda_{\pi(k)}] + \delta(x),1) & \mbox{if }\ti > d_0,}$$
    for $(x,0)\in \hI\setminus\Sigma_0\cap I\times\{0\}$ where $\ti$ is defined by $\{i,\ti\}=\inv{\pi}(\pi(i))$. Now let $i=i(x)\in\{d_0+1,\dots,2d\}$ be defined by
        $$ \dusum{k=d_0+1}{i-1}\lambda_{\pi(k)} \leq x < \dusum{k=d_0+1}{i}\lambda_{\pi(k)}$$
    and $ \delta(x) = x - \dusum{k=d_0+1}{i-1}\lambda_{\pi(k)}.$ Then
        $$ \tT(x,1) = \RHScase{([\dusum{k=1}{\ti-1}\lambda_{\pi(k)}]+\delta(x),0), & \mbox{if }\ti\leq d_0, \\ ([\dusum{k=d_0+1}{\ti}\lambda_{\pi(k)}] - \delta(x), 1), & \mbox{if }\ti > d_0,}$$
    for $(x,1)\in \hI\setminus\Sigma_0\cap I\times\{1\}$. Then finally let $T=f\circ \tT$.

    \begin{defn}
        We say that $T\sim(\pi,\lambda)$ if $T$ agrees with the linear involution defined by $(\pi,\lambda)$.
    \end{defn}

    \begin{exam}
        Let
        $$ \pi = \cmtrx{\LL{a}{c}~\LL{b}{c}~\LL{b}{a}}$$
        and $\lambda=(\lambda_a,\lambda_b,\lambda_c)=(\frac{1}{2},\frac{1}{4},\frac{1}{4})$. Then $T\sim(\pi,\lambda)$ acts by
        $$ T(x,\eps) = \RHScase{(x+\frac{1}{2},0), & \mbox{if }x\in(0,\frac{1}{2}), \eps=0,\\
                        (\frac{3}{2}-x,1) &       \mbox{if }x\in(\frac{1}{2},\frac{3}{4})\cup(\frac{3}{4},1), \eps=0,\\
                        (\frac{1}{2}-x,0) &     \mbox{if }x\in(0,\frac{1}{4})\cup(\frac{1}{4},\frac{1}{2}), \eps=1,\\
                        (x-\frac{1}{2},1) &     \mbox{if }x\in(\frac{1}{2},1), \eps=1.}$$
    \end{exam}

    So we see that, similar to \IET's, a linear involution may be encoded by combinatorial data $\pi$ and length data $\lambda$. We will now address the existence of $\lambda\in\RR_+^{\AAA,\pi}$ for given $\pi\in\genperm_\AAA$.

    \begin{defn}
        A generalized permutation $\pi$ is \term{proper} when
        $$ \exists \alpha\in\AAA\mbox{ s.t. }\inv{\pi}(\alpha)\subset\{1,\dots,d_0\} \iff \exists \beta\in\AAA\mbox{ s.t. }\inv{\pi}(\beta)\subset\{d_0+1,\dots,2d\}.$$
    \end{defn}

    \begin{prop}
        Let $\pi\in\genperm_\AAA$.
        $\RR^{\AAA,\pi}_+\neq\emptyset$ if and only if $\pi$ is proper.
    \end{prop}

    \begin{proof}
        We partition $\AAA$ by the following sets:
        $$ \begin{array}{rcl}
            \AAA_0^\pi &:=&\{\alpha\in\AAA: \inv\pi(\alpha)\subset\dset{d_0}\},\\
            \AAA_1^\pi &:=&\{\alpha\in\AAA:\inv\pi(\alpha)\subset\ddset{d_0+1}{2d}\},\\
            \AAA_{0,1}^\pi &:=& \{\alpha\in\AAA: D_\pi(\alpha)<1\}.
            \end{array}$$
        Equation \eqref{eq.gen.lengths} then is equivalent to the equation
        \begin{equation}\label{eq.gen.lengths2}
            \sum_{\alpha\in\AAA^\pi_0} \lambda_\alpha = \sum_{\beta\in\AAA^\pi_1} \lambda_\beta.
        \end{equation}
        Assume $\pi$ is not proper with $\eps\in\{0,1\}$ such that $\AAA^\pi_\eps \neq\emptyset$ and $\AAA^\pi_{1-\eps}=\emptyset$. Then \eqref{eq.gen.lengths2} would imply
        $$ 0 = \sum_{\alpha\in\AAA^\pi_{1-\eps}}\lambda_\alpha = \sum_{\beta\in\AAA^\pi_\eps}\lambda_\beta > 0$$
        for any $\lambda\in\RR_+^{\AAA,\pi}$. Therefore, $\RR_+^{\AAA,\pi}$ is empty.\\

        Now suppose that $\pi$ is proper. Then the sets $\AAA_0^\pi$ and $\AAA_1^\pi$ are either both empty or both nonempty. If they are empty then $\pi$ is a true permutation, \eqref{eq.gen.lengths2} is vacuously true and $\RR_+^{\AAA,\pi}=\RR_+^\AAA\neq\emptyset$. So assume both sets are not empty with $n_\eps:=\#\AAA_\eps^\pi>0$. We can then construct $\lambda\in\RR_+^\AAA$ by
        $$ \lambda_\alpha = \RHScase{\frac{1}{n_0}, & \mbox{if }\alpha \in \AAA_0^\pi, \\ 1, & \mbox{if }\alpha\in\AAA_{0,1}^\pi, \\ \frac{1}{n_1}, & \mbox{if }\alpha\in\AAA_1^\pi.}$$
        So $\lambda\in\RR_+^{\AAA,\pi}$ as $\lambda$ satisfies \eqref{eq.gen.lengths2}. Therefore $\RR_+^{\AAA,\pi}\neq\emptyset$.
    \end{proof}

    \begin{rem}
        Throughout this paper, we will assume that all generalized permutations are proper, as they are the only generalized permutations realizable by linear involutions.
    \end{rem}

    \begin{defn}
        Let $F_i\subset \AAA$, $i\in\{1,\dots,4\}$, be finite subsets, $n_i=\# F_i$, and $\pi\in\genperm_\AAA$. We say that $\pi$ admits a \term{corner decomposition} with corners $F_1,\dots,F_4$ if
            $$ \begin{array}{rcl}
                \pi(\{1,\dots,n_1\}) &=& F_1, \\
                \pi(\{d_0 - n_2 + 1,\dots,d_0\}) & = & F_2, \\
                \pi(\{d_0+1,\dots,d_0+n_3\}) & = & F_3,\mbox{ and }\\
                \pi(\{2d-n_4+1,\dots,2d\}) & = & F_4.\end{array}$$
        If this is the case, we denote the decomposition as
            $$ \pi = \cmtrx{ F_1 & \dots & F_2 \\ F_3 & \dots & F_4}.$$
    \end{defn}

    \begin{defn}\label{def.gen.irr}
        A generalized permutation $\pi\in\genperm_\AAA$ is \term{reducible} if there exists disjoint $A,B,C,D\subsetneq \AAA$ not all empty such that
        $$ \pi = \cmtrx{A\sqcup B & \dots & A\sqcup C \\ B\sqcup D & \dots & C \sqcup D}$$
        and one of the following holds:
        \begin{itemize}
            \item No corner is empty.
            \item Only one corner is empty and is on the left.
            \item Exactly two corners are empty and are either both on the left or right.
        \end{itemize}
        If $\pi$ is not reducible, $\pi$ is \term{irreducible}. The set of all irreducible generalized permutations on $\AAA$ shall be denoted as $\genirr_\AAA$.
    \end{defn}

    \begin{exam}
        If $\pi\in\genperm_\AAA$ is a true permutation, then $A=D=\emptyset$. It follows that Definition \ref{def.gen.irr} coincides with the definition of irreducible at the beginning of Section \ref{sec.IET} for true permutations. Therefore Definition \ref{def.gen.irr} is indeed a generalization of irreducibility when applied to generalized permutations.
    \end{exam}

    \begin{defn}\label{def.gen.inv}
        For $\pi\in\genperm_\AAA$, the \term{inverse} of $\pi$ is
        $$ \inv\pi = \cmtrx{\pi(d_0+1),\pi(d_0+2),\dots,\pi(2d) \\ \pi(1),\pi(2),\dots,\pi(d_0)}$$
        or denoted $(\inv{\pi},d_1,d_0)$ with
        $$ \inv\pi(i) = \RHScase{\pi(i+d_0), & \mbox{if }1\leq i\leq d_1, \\ \pi(i-d_1), & \mbox{if }d_1<i\leq 2d.}$$
    \end{defn}
    Again, it is clear that this is a natural generalization of inverse for true permutations.

\sect{Generalized Suspensions and Quadratic Differentials}\label{sec.gen.surfaces}

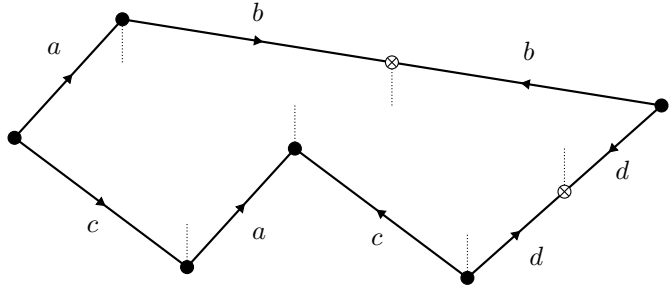
\begin{figure}[h]
        \begin{center}
           \setlength{\unitlength}{250pt}
            \begin{picture}(1,.46)
                \put(0,0){\includegraphics[width=\unitlength]{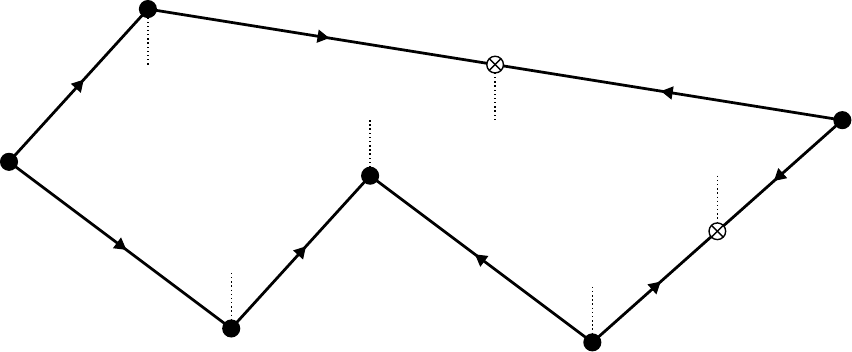}}

                \put(.06,.35){$a$} \put(.37,.40){$b$} \put(.78,.34){$b$}

                \put(.12,.08){$c$} \put(.37,.07){$a$} \put(.55,.06){$c$} \put(.79,.03){$d$} \put(.92,.16){$d$}
           \end{picture}
       \end{center}
       \caption{A suspension $S$ for $\pi=\{a~b~b~/~c~a~c~d~d\}$. $S$ has three singularities: one of order $2$ and two of order $-1$.}\label{fig.gen.sus_exam}
    \end{figure}

    Let $T\sim(\pi,\lambda)$ for $\pi\in\genperm_\AAA$ and $\lambda\in\RR^{\AAA,\pi}_+$ be as in the previous section. We perform a process to suspend a surface over $T$ analogous to the suspension of a regular \IET\ in Section \ref{sec.surface}. We define the cone
    \begin{equation}\label{eq.gen.Tpi}
        \begin{array}{rcl} \TTT_\pi:=\{\tau\in\RR^\AAA & : & \sum_{j=1}^{k}\tau_{\pi(j)} > 0 \mbox{ for all }1\leq k<d_0,\\
            & &\sum_{j=d_0+1}^{k} \tau_{\pi(j)}<0\mbox{ for all }d_0<k<2d,\\
            & &\sum_{j=1}^{d_0} \tau_{\pi(j)} = \sum_{j=d_0+1}^{2d}\tau_{\pi(j)}\}.\end{array}
    \end{equation}
    We have the following from Boissy and Lanneau (\cite{c.BL07}, Theorem 3.2),
    \begin{prop}
        [Boissy-Lanneau 07] Let $\pi\in\genperm_\AAA$. $\TTT_\pi\neq\emptyset$ if and only if $\pi\in\genirr_\AAA$.
    \end{prop}
    So assume the triple $(\pi,\lambda,\tau)$ for $\tau\in\TTT_\pi$. We now construct a surface with these constraints. Define the vectors $\zeta_\alpha = (\lambda_\alpha,\tau_\alpha)$ for each $\alpha\in\AAA$. Start at point $(0,0)$ and attach these vectors in order determined by the first row of $\pi$, i.e. $\zeta_{\pi(1)}$ followed by $\zeta_{\pi(2)}$ up until $\zeta_{\pi(d_0)}$. Beginning again at $(0,0)$, attach vectors $\zeta_{\pi(d_0+1)},\zeta_{\pi(d_0+2)},\dots,\zeta_{\pi(2d)}$ in order. If the two occurrences of $\zeta_\alpha$ appear both above or below the $x$-axis, identify their interiors by a function of the form $z\mapsto -z+c_\alpha$. If the two $\zeta_\alpha$'s appear on opposite sides of the $x$-axis, identify their interiors by translation, $z\mapsto z+c_\alpha$. Denote this suspension by $S:=S(\pi,\lambda,\tau)$. Figure \ref{fig.gen.sus_exam} gives an example suspension for $\pi=\cmtrx{\LL{\msp}{c}~\LL{a}{a}~\LL{b}{c}~\LL{b}{d}~\LL{\msp}{d}}$.\\

    Consider $\pi\in\genirr_\AAA$, $\lambda\in\RR_+^{\AAA,\pi}$ and $\tau\in\TTT_\pi$. Let $L = \dusum{i=1}{d_0} \lambda_{\pi(i)}$ as in Equation \eqref{eq.gen.lengths}. Given $S=S(\pi,\lambda,\tau)$, let $I_S=(0,L)\times\{0\}$ be the natural embedding of $I=(0,L)$. We see that the first return of the upward and downward trajectories on $I_S$ realizes the linear involution $T\sim(\pi,\lambda)$, under the interpretation in Remark \ref{rem.gen.interpret}.\\

    If $\pi\in\genirr_\AAA$ is proper and not a true permutation, the suspended surface $S$ admits a quadratic differential that is not the square of an Abelian differential. The vertex classes of $S$ represent singularities of the differential. We consider the vertex class of $v$ a singularity of order $k$ if the total angle around $[v]$ in $S$ is $\pi(k+2)$. We use the term \term{order} to distinguish from the term \term{degree} in the Abelian differential case. We also remark that the number and orders of the singularities is again independent of choice of suspension over $S$. Using Gauss-Bonnet, if there are $k$ singularities of orders $\ell_1,\dots,\ell_k$, then the genus $g$ of the surface $S$ (another constant of $\pi$) is related to the $\ell_i$'s by the equation
    $$ \sum_{i=1}^k \ell_i = 4g-4.$$

    We shall denote by $\QQQ_g$ the moduli space of all quadratic differentials on a surface of genus $g\geq 0$ that are not squares of Abelian differentials. These spaces are further stratified by the number and orders of singularities, denoted by $\QQQ(\ell_1,\dots,\ell_k)$ for $\ell_i\geq -1$ such that $\ell_1+\dots+\ell_k=4g-4$. In Section \ref{sec.gen.classify}, we shall discuss which strata are non-empty and how many connected components each strata contains.\\

    We now give a counting statement similar to the Abelian differential case. Consider $\pi\in\genirr_\AAA$ and suspension $S$. We label each segment's endpoints by $(i,\jmath)$, where $i\in\dset{2d}$ and $\jmath\in\{L,R\}$, denoting the left or right endpoint of the segment related to $\pi(i)$. We have the following equivalencies:
    \begin{enumerate}
        \item $(i,R)\sim (i+1,L)$ for $i\in\{1,\dots,d_0-1,d_0+1,\dots,2d-1\}$.
        \item $(1,L)\sim(d_0+1,L)$ and $(d_0,R)\sim(2d,R)$.
        \item Let $\alpha=\pi(i)=\pi(\ti)$ and $\jmath\in\{L,R\}$. If $D_\pi(\alpha)<0$ (if $i$ and $\ti$ are in the opposite rows), then
            $$ (i,\jmath)\sim(\ti,\jmath).$$
            If $D_\pi(\alpha)>0$ (if $i$ and $\ti$ are both in the same row), then
            $$ (i,\jmath)\sim(\ti,\tilde\jmath)$$
            where $\tilde\jmath$ is the opposite side of $\jmath$.
    \end{enumerate}
    In order to determine the order of a given singularity, we note that if $k$ is its order then there must be $k+2$ distinct vertical directions coming from the vertex class of the singularity in $S$. In other words, if there are $n$ vertices identified with our singularity (excluding the rightmost and leftmost vertices), then our singularity is of order $n-2$ (we note that we are now counting the top and bottom rows). For example, we see from Figure \ref{fig.gen.sus_exam} that the suspension for $\pi=\cmtrx{\LL{\msp}{c}~\LL{a}{a}~\LL{b}{c}~\LL{b}{d}~\LL{\msp}{d}}$ has three singularities. One singularity is of order $2$, while the other two singularities are \term{poles} (of order $-1$).

    \begin{rem}
        Just as in the case for true permutations, the number and orders of the singularities for $\pi$ do not vary by choice of $\lambda$ and $\tau$. It follows that the genus, $g(\pi)$ is well defined.
    \end{rem}

\sect{Generalized Rauzy Classes}\label{sec.gen.classes}

    Similar to IETs, we define Rauzy induction on linear involutions. Let $T\sim(\pi,\lambda)$ for $\pi\in\genirr_\AAA$ and $\lambda\in\RR^{\AAA,\pi}_+$ be a map on $I=(0,L(\lambda))$ and $\hI=I\times\{0,1\}$ (see Section \ref{sec.gen.perm}). Let map $T'$ be the first return map on $\hI':=I'\times\{0,1\}$ where $I':=(0,L-\max\{\lambda_{\pi(d_0)},\lambda_{\pi(2d)}\})$. We recall for $i\in\dset{2d}$ the value $\ti$ given by $\{i,\ti\}=\inv\pi(\pi(i))$.

    \begin{defn}\label{def.gen.RV}
        Let $\pi=(\pi,d_0,d_1)\in\genirr_\AAA$. We define Ruazy induction of each type as follows:
        \begin{itemize}
            \item Let $\rb\pi$ be the move of type $\rtb$ on $\pi$. If $\td_0>d_0$, then
            $$ \rb\pi(i) := \RHScase{\pi(i), & \mbox{if }i\leq\td_0,\\ \pi(2d), & \mbox{if }i = \td_0+1,\\ \pi(i-1), & \mbox{if }i > \td_0+1,}$$
            and $\rb\pi = (\rb\pi,d_0,d_1)$. If $\td_0<d_0$, then
            $$ \rb\pi(i) := \RHScase{\pi(i), & \mbox{if }i<\td_0,\\ \pi(2d), & \mbox{if }i=\td_0,\\ \pi(i-1), & \mbox{if }i>\td_0,}$$
            and in this case $\rb\pi=(\rb\pi,d_0+1,d_1-1)$ is defined only if it is proper.

            \item Let $\rt\pi$ be the move of type $\rtt$ on $\pi$. If $\tdd<d_0$, then
            $$ \rt\pi(i) := \RHScase{\pi(i), & \mbox{if }i\leq \tdd,\\ \pi(d_0), & \mbox{if }i=\tdd+1,\\ \pi(i-1), & \mbox{if }\tdd+1<i\leq d_0,\\ \pi(i), & \mbox{if }i > d_0,}$$
            and $\rt\pi=(\rt\pi,d_0,d_1)$. If $\tdd>d_0$, then
            $$ \rt\pi(i) := \RHScase{\pi(i), & \mbox{if }i<d_0,\\ \pi(i+1), & \mbox{if }d_0\leq i <\tdd-1,\\ \pi(d_0), & \mbox{if } i=\tdd-1,\\ \pi(i), & \mbox{if }\tdd \leq i \leq 2d,}$$
            and in this case $\rt\pi=(\rt\pi,d_0-1,d_1+1)$ is defined only if it is proper.
        \end{itemize}
    \end{defn}

    \begin{rem}\label{rem.gen.inv_conv}
        This definition agrees with the definition in Section \ref{sec.RC} if $\pi$ is a true permutation. It can also be shown similarly that for $\pi\in\genirr_\AAA$ and $\eps\in\{\rtb,\rtt\}$, $\inv{(\reps\pi)}=\ropp\inv{\pi}$.
    \end{rem}

    So we see that if $\lambda_{\pi(d_0)}>\lambda_{\pi(2d)}$, $T' = (\rb\pi,\lambda')$ where $\lambda'_{\pi(d_0)}=\lambda_{\pi(d_0)}-\lambda_{\pi(2d)}$ and $\lambda'_\alpha=\lambda_\alpha$ otherwise. Likewise if $\lambda_{\pi(d_0)}<\lambda_{\pi(2d)}$, $T'=(\rt,\lambda')$ where $\lambda'_{\pi(2d)}=\lambda_{\pi(2d)}-\lambda_{\pi(d_0)}$ and $\lambda'_\alpha=\lambda_\alpha$ otherwise. It can be shown that if $T'=(\pi',\lambda')$, then $\lambda'\in\RR^{\AAA,\pi'}_+$. Also, as in the case of regular IETs, $T'$ is not well defined if $\lambda_{\pi(d_0)}=\lambda_{\pi(2d)}$, which is a Lebesgue measure zero condition on $\RR^{\AAA}_+$ and $\RR^{\AAA,\pi}_+$.\\

    \begin{rem}
        Consider $\pi\in\genirr_\AAA$. Let $\eps\in\{0,1\}$ be such that $\reps\pi$ is not well defined. It follows that there does not exist $\lambda\in\RR_+^{\AAA,\pi}$ such that:
        \begin{itemize}
            \item $T'\sim(\pi',\lambda')$ exists, and
            \item $\pi'=\reps\pi$.
        \end{itemize}
        By this we see that if a move is not defined for $\pi$, no linear involution $T\sim(\pi,\lambda)$ realizes this move.
    \end{rem}

    The following definition and proposition (with proof) may be seen in Section $4$ of \cite{c.BL07}. This definition is the linear involution generalization of Definition \ref{def.Keane}.

    \begin{defn}\label{def.gen.KEANE}
        A linear involution $T\sim(\pi,\lambda)$ on $\hI$ has a \term{connection} if there exists $(x,\eps)\in I\times\{0,1\}$ and $r\geq 0$ such that
        \begin{itemize}
            \item $(x,\eps)$ is a singularity of $\inv{T}$, and
            \item $T^{r}(x,\eps)$ is a singularity of $T$.
        \end{itemize}
        If $T$ does not have a connection, it satisfies the \term{Keane Condition}.
    \end{defn}

    This proposition is the linear involution equivalent of Proposition \ref{prop.Keane_is_inducible}.

    \begin{prop}\label{prop.gen.KEANE_is_Inducible}
        Let $T\sim(\pi,\lambda)$ be a linear involution. Then the following are equivalent:
        \begin{enumerate}
            \item $T$ satisfies the Keane Condition.
            \item $T^{(n)}$, the $n^{th}$ iteration of R-V induction on $T$, is defined for all $n\geq 0$. Also, if $T^{(n)}=(\pi^{(n)},\lambda^{(n)})$, then $\lambda^{(n)}\to 0$ as $n\to \infty$.
        \end{enumerate}
        If these conditions are satisfied, $T$ is minimal.
    \end{prop}

    \begin{prop}\label{prop.gen.irr_ind}
        Let $\pi\in\genirr_\AAA$ and $\eps\in\{\rtb,\rtt\}$. If $\reps\pi$ is defined, then $\reps\pi\in\genirr_\AAA$.
    \end{prop}

    \begin{proof}
        By Remark \ref{rem.gen.inv_conv}, it suffices to prove the claim for $\pi':=\rb\pi$ and assume it exists, for $\pi\in\genirr_\AAA$. Assume $\pi'$ is reducible. Define the sets $A,B,C,D\subset\AAA$ by the corner decomposition
        $$ \pi'=\cmtrx{A\sqcup B & \dots & A\sqcup C \\ B\sqcup D & \dots & C\sqcup D}.$$
        Let $\alpha:=\pi'(d_0')$ be the last letter of the first row in $\pi'$. We consider each case from Definition \ref{def.gen.irr}:
        \begin{itemize}
            \item Suppose no corner is empty. Then either $\alpha\in A$ or $\alpha\in C$. If $\alpha\in C$, then $\pi$ has a corner decomposition with the same sets $A,B,C,D$, a contradiction to the irreducibility of $\pi$. Now if $\alpha\in A$, let $\beta$ appear before the first occurrence of $\alpha$ or $\beta:=\pi'(d_0^{\pi'}-1)$. Then
                    $$ \pi = \cmtrx{A'\sqcup B' & \dots & A'\sqcup C' \\ B'\sqcup D' & \dots & C'\sqcup D'}$$
                where
                $$ \begin{matrix}
                    {A' = \RHScase{A\setminus\{\beta\}, & \mbox{if }\beta\in A, \\ A, & \mbox{if }\beta\in B,}} &
                    {B' = \RHScase{B, & \mbox{if }\beta\in A, \\ B\setminus\{\beta\}, & \mbox{if }\beta\in B,}} \\
                    {C' = \RHScase{C\sqcup\{\beta\}, & \mbox{if }\beta\in A, \\ C, & \mbox{if }\beta\in B,}} &
                    {D' = \RHScase{D, & \mbox{if }\beta\in A, \\ D\sqcup\{\beta\}, & \mbox{if }\beta\in B,}}
                    \end{matrix}$$
                which again shows that $\pi$ is reducible.
            \item Suppose the top left corner is empty, or $A=B=\emptyset$. Then $\pi$ has the same corner decomposition as $\pi'$.
            \item Suppose the bottom left corner is empty, or $B=D=\emptyset$. If $\alpha\in C$, then $\pi$ has the same corner decomposition as $\pi'$ and is therefore reducible. Now if $\alpha\in A$, then
                    $$ \pi = \cmtrx{A' & \dots & A'\sqcup C' \\ \emptyset & \dots & C'}$$
                where $A':= A\setminus\{\beta\}$ and $C':=C\sqcup\{\beta\}$ for $\beta:=\pi'(d_0^{\pi'}-1)$ as before. Then we see again that $\pi$ is reducible.
            \item Suppose two corners on either the right or left are empty, or either $A=B=D=\emptyset$ or $B=C=D=\emptyset$. In either case, $\pi$ has the same corner decomposition as $\pi'$ and is therefore reducible.
        \end{itemize}
        So after exhausting all cases, we see that if $\pi'$ exists and is reducible, then $\pi$ must be reducible as well. This is a contradiction of the assumption on $\pi$.
    \end{proof}

    \begin{exam}
        Consider the following general permutations:
        $$ \pi = \cmtrx{\LL{a}{\msp}~\LL{a}{d}~\LL{b}{c}~\LL{b}{d}~\LL{c}{\msp}} \mbox{ and } \pi'=\rt\pi = \cmtrx{\LL{a}{c}~\LL{a}{d}~\LL{b}{c}~\LL{b}{d}}.$$
        We see that while $\pi'\in\genirr_\AAA$, $\pi$ is reducible. This example shows that the converse of Proposition \ref{prop.gen.irr_ind} is in general false.
    \end{exam}

\begin{figure}[h]
        \begin{center}
           \setlength{\unitlength}{350pt}
            \begin{picture}(1,.30)
                \put(0,0){\includegraphics[width=\unitlength]{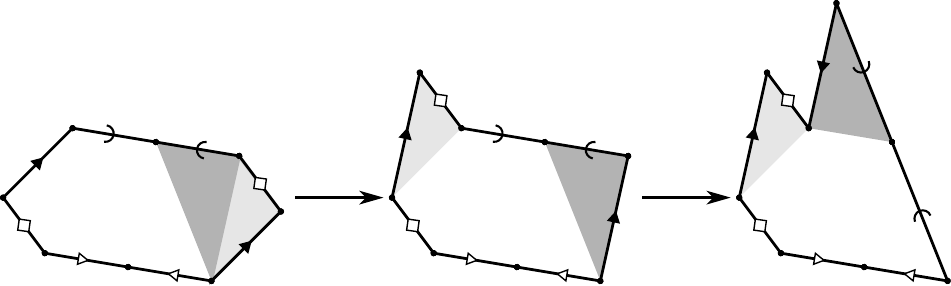}}

                \put(.14,.075){$S$}
                \put(.345,.115){$\rtt$}
                \put(.52,.075){$S'$}
                \put(.71,.115){$\rtb$}
                \put(.87,.075){$S''$}

            \end{picture}
       \end{center}
       \caption{Two successive iterations of Rauzy-Veech induction on a surface.}\label{fig.gen.RVSurface}
    \end{figure}
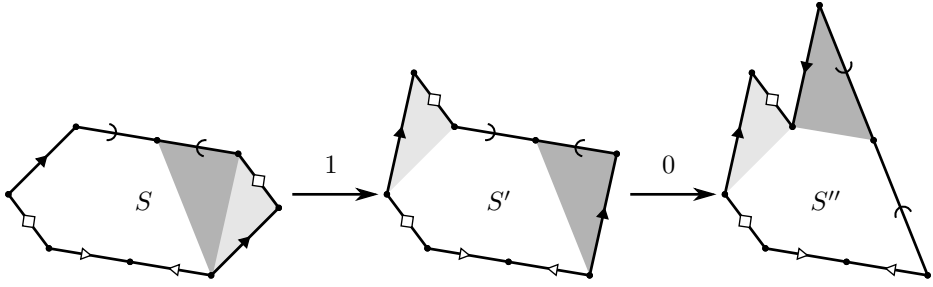

    \begin{rem}
        R-V Induction may be defined for suspensions as well. If $S=S(\pi,\lambda,\tau)$ is a suspension for $\pi\in\genirr_\AAA$, then $S'=S(\pi',\lambda',\tau')$ is a suspension for $\pi'$ where $\pi'$ and $\lambda'$ are from Definition \ref{def.gen.RV} and $\tau'\in\TTT_{\pi'}$ is given by
            $$ \tau_\alpha' = \RHScase{\tau_\beta - \tau_{\beta'}, & \mbox{if }\alpha = \beta,\\ \tau_\alpha, & \mathrm{otherwise,}}$$
        where $\pi'=\reps\pi$, $\beta=\alpha_\eps$ and $\beta'=\alpha_{1-\eps}$ for $\alpha_0 = \pi(d_0)$, $\alpha_1 = \pi(2d)$. See Figure \ref{fig.gen.RVSurface}. We note that this procedure is still a ``cut and paste," like in the true permutation case. However, the pieces may be moved by translation composed with an order two rotation.
    \end{rem}

    \begin{defn}\label{def.gen.RC}
        Let $\pi\in\genirr_\AAA$, then the \term{Rauzy Class} of $\pi$, $\RClass(\pi)$, is the smallest set of generalized permutations containing $\pi$ and closed under both types of inductive moves. The \term{Rauzy Graph} of $\RClass(\pi)$ is the directed graph with vertices in $\RClass(\pi)$ and edges representing the (valid) moves of induction.
    \end{defn}

    \begin{rem}
        An equivalent definition for a generalized Rauzy class on $\pi$ is the set of all forward iterates of induction on $\pi$, as in \cite{c.BL07}. In that paper, Boissy and Lanneau showed the following: if there exists a series of inductive steps from $\pi$ to $\pi'$ then there exists a series of steps from $\pi'$ to $\pi$. Their proof does not follow in the same spirit of Claim \ref{cor.1}, as there are $\pi\in\genirr_\AAA$ such that a series of moves of the same type terminate and do not form a cycle.
    \end{rem}

\begin{figure}[h]
        \begin{center}
           \setlength{\unitlength}{250pt}
            \begin{picture}(1,.52)
                \put(0,0){\includegraphics[width=\unitlength]{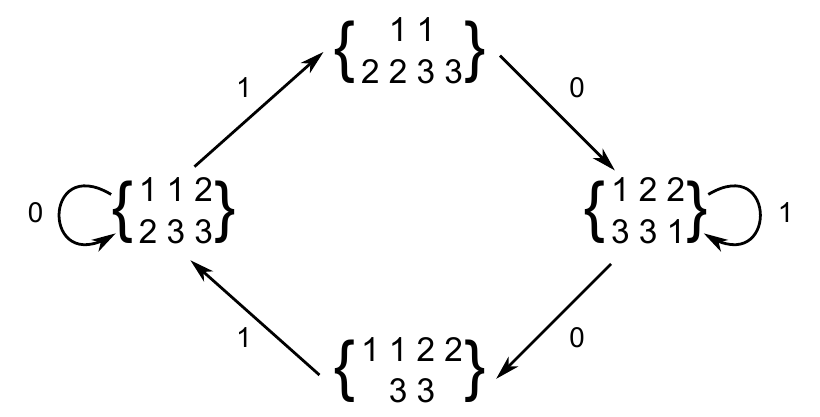}}
            \end{picture}
        \end{center}
        \caption{A generalized Rauzy Graph.}\label{fig.gen.graph_of_122_331}
    \end{figure}

    \begin{exam}
        Consider $ \pi = \cmtrx{\LL{1}{3}~\LL{2}{3}~\LL{2}{1}}.$ We see that $\rt\pi=\pi$ and $ \pi':=\rb\pi = \cmtrx{\LL{1}{\msp}~\LL{1}{3}~\LL{2}{3}~\LL{2}{\msp}}.$ $\rb\pi'$ is not defined, as it would not be proper. Let $ \pi'' := \rt\pi' = \cmtrx{\LL{1}{2}~\LL{1}{3}~\LL{2}{3}}.$ Note that $\rb\pi'' = \pi''$, and let $ \pi''' := \cmtrx{\LL{\msp}{2}~\LL{1}{2}~\LL{1}{3}~\LL{\msp}{3}}.$ Finally, we see that $\rt\pi'''$ is not well defined and $ \rb\pi''' = \cmtrx{\LL{3}{2}~\LL{1}{2}~\LL{1}{3}} = \pi.$ So $\RClass(\pi) = \{\pi,\pi',\pi'',\pi'''\}$, and its Rauzy Graph is represented in Figure \ref{fig.gen.graph_of_122_331}.
    \end{exam}

    \begin{rem}
        R-V induction acting on a surface preserves the number and orders of singularities of the surface. Therefore the number and orders of singularities remains fixed throughout the entire Rauzy Class (and consequently the genus). The leftmost singularity is again fixed by R-V Induction, and is still referred to as the \term{marked singularity}. However the leftmost letters $\pi(1)$ and $\pi(d_0+1)$ do not remain fixed within a class, a contrast to the case of true permutations.
    \end{rem}

    \begin{defn}
        Given Rauzy Class $\RRR$, the \term{signature} of $\RRR$, $\sig(\RRR)$, is the $m$-tuple $(\ell_1,\dots,\ell_m)$ where $m$ is the number of singularities and $\ell_i$'s are the orders with multiplicity. The value of $\ell_1$ denotes the order of the marked singularity. For any $\pi\in\RClass$, $\sig(\pi):=\sig(\RClass)$.
    \end{defn}

\sect{Classification of Generalized Rauzy Classes}\label{sec.gen.classify}

    We note that as opposed to the strata of Abelian differentials, E. Lanneau shows in \cite{c.L04Spin} that the parity of the spin structure is determined by the genus and therefore cannot distinguish connected components in each stratum $\QQQ(\ell_1,\dots,\ell_m)$.

    \begin{nonum}(M. Kontsevich \cite{c.K97}; H. Masur and J. Smillie \cite{c.MS}; E. Lanneau \cite{c.L08})
    All connected components of any stratum of meromorphic quadratic differentials with at most simple poles on a complex curve of genus $g\leq 2$ are described as follows:\\

    \begin{itemize}
        \item
        For $g=0$, all strata are nonempty and connected.
        \item
        For $g=1$, $\QQQ(\emptyset)$ and $\QQQ(1,-1)$ are empty. All other strata are nonempty and connected.
        \item
        For $g=2$, $\QQQ(3,1)$ and $\QQQ(4)$ are empty. The strata $\QQQ(6,-1^2)$ and $\QQQ(3,3,-1^2)$ contain exactly two connected components, one is hyperelliptic and one is not. All other strata are nonempty and connected.
    \end{itemize}
    \end{nonum}

    The result for $g=0$ is due to M.\ Kontsevich. The results on the empty strata are due to H.\ Masur and J.\ Smillie. The remaining results for $g=1,2$ are by E.\ Lanneau.

    \begin{nonum}(E. Lanneau \cite{c.L08}) All connected components of any stratum of meromorphic quadractic differentials with at most simple poles on a complex curve of genus $g\geq 3$ are described as follows:\\

    \begin{itemize}
    \item
        The four exceptional strata
        $$ \QQQ(9,-1), \QQQ(6,3,-1), \QQQ(3,3,3,-1), \QQQ(12)$$
        have two connected components, denoted $\QQQ^{reg}$ and $\QQQ^{irr}$.

    \item
        The following strata
        $$ \begin{array}{lll}
                \QQQ(2j-1, 2j-1, 2k-1, 2k-1), & j,k\geq 0, & j+k=g\\
                \QQQ(2j-1,2j-1,4k+2), & j,k\geq 0, & j+k=g-1\\
                \QQQ(4j+2, 4k+2), & j,k\geq 0, &  j+k=g-2
            \end{array}$$
        have exactly two connected components, one is hyperelliptic and one is not.

    \item
        All other such strata are nonempty and connected.

    \end{itemize}
    \end{nonum}

    The above theorems only categorize the connected components. As in Section \ref{sec.classify}, C.\ Boissy in \cite{c.B09} provides the final step in classification of all generalized Rauzy classes:

    \begin{nonum}(C.\ Boissy \cite{c.B09})
        Two generalized permutations belong to the same Rauzy class if and only if they belong to the same connected component and have marked singularities of the same order.
    \end{nonum}

    We conclude that Corollary \ref{cor.classify} applies in the case of generalized permutations as well.

\chap{True Rauzy Classes}\label{chapTrue}

	\begin{thm}\label{thm.main}
	(Main Result) Every (true) Rauzy Class contains a permutation $\pi$ such that $\pi=\inv{\pi}$.
	\end{thm}

    We prove this by using Corollary~\ref{cor.classify} to identify each special case of Rauzy Class. The hyperelliptic case is covered as the only standard element of each hyperelliptic class, $(d,d-1,\dots,2,1)$ as shown in Proposition \ref{prop.hyp}, is its own inverse. The remainder of the connected components shall be covered by Theorems~\ref{thm.g3},~\ref{thm.odds}, \ref{thm.even.odd}, \ref{thm.twos.even}, and \ref{thm.even.even}. The special case of component $\Hnonhyp(2\ell-1,2\ell-1)$ mentioned in Section \ref{sec.classify} is covered in Theorem \ref{thm.odds}, as the permutation constructed is not hyperelliptic so can apply to these components. Finally, singularities of degree zero are considered in Theorem \ref{thm.zero}. \\
	
    While the following fact can be deduced from the works \cite{c.B09} and \cite{c.KZ}, we can now state an alternate proof as a direct result.
	\begin{corollary}\label{cor.main}
		Every Rauzy Class is closed under taking inverses.
	\end{corollary}
	
	\begin{proof}
    For any Rauzy Class $\RClass$, we may pick $\pi'\in\RClass$ that is self inverse by Theorem \ref{thm.main}. Now choose any $\pi\in\RClass$. By Claim \ref{cor.1}, we may choose a series $\rsub{\eps_1}\dots \rsub{\eps_k}$, $\eps_i\in\{0,1\}$, such that $\pi=\rsub{\eps_k}\rsub{\eps_{k-1}}\dots \rsub{\eps_2}\rsub{\eps_1}\pi'$. By Claim \ref{cor.2},
		$$ \begin{array}{rcl}
			\inv{\pi}&=&\inv{(\rsub{\eps_k}\rsub{\eps_{k-1}}\dots \rsub{\eps_2}\rsub{\eps_1}\pi')}\\
				&=&\rsub{1-\eps_k}\rsub{1-\eps_{k-1}}\dots \rsub{1-\eps_2}\rsub{1-\eps_1}\inv{\pi'}\\
				&=&\rsub{1-\eps_k}\rsub{1-\eps_{k-1}}\dots \rsub{1-\eps_2}\rsub{1-\eps_1}\pi'.
		\end{array}$$
		 So $\pi^{-1}\in\RClass$ as $\pi'\in\RClass$.
	\end{proof}

    \begin{corollary}\label{cor.main2}
        In every connected component of every stratum $\HHH(\ell_1,\dots,\ell_m)$, there exists a differential that allows an order two orientation reversing linear isometry.
    \end{corollary}

    \begin{proof}
        In every connected component $\CCC\subset\HHH(\ell_1,\dots,\ell_m)$, consider a Rauzy Class $\RClass$ contained in $\CCC$. Choose self-inverse $\pi\in\RClass$ by Theorem \ref{thm.main}. Let $\mone=(1,\dots,1)\in\RR_+^\AAA$ and $\tau=(1,0,\dots,0,-1)\in\TTT_\pi$. Let $S=S(\pi,\mone,\tau)$ and note that $h(x,y)=(x,-y)$ satisfies the claim.
    \end{proof}

\sect{Spin Parity for Standard Permutations}\label{sec.stdspin}

    When $\pi=(\pi_0,\pi_1)\in\irr_d$ is standard, the calculations mentioned in Section \ref{sec.spin} can be further refined. Just as in that section, the following calculations are over $\ZZ_2$. The matrix $\Omega=\Omega_\pi$ has the following form
	\begin{equation}
		\Omega = \pmtrx{	0& 1 & \dots&\dots & 1 & 1\\
					1& \mA_1& \mO& \dots &\mO& 1\\
					1& \mO& \mA_2 & \mO&\vdots & 1\\
					\vdots& \vdots&\ddots&\ddots&\ddots&\vdots\\
					1&\mO&\dots&\mO&\mA_p&1\\
					1&1&\dots&\dots&1&0}
	\end{equation}
    In other words, along the rows $1$ and $d$ and columns $1$ and $d$, the entries are all $1$ except entries $\Omega_{1,1}$ and $\Omega_{d,d}$ which are $0$. The interior $d-2$ by $d-2$ matrix is composed of $p$ square matrices, labeled $\mA_i$, along the diagonal with zeros otherwise. Note that $p=1$ is allowed. Each matrix $\mA_i$ corresponds to sub-alphabet $\AAA_i$ such that, for $\eps\in\{0,1\}$, $\pi_\eps(\AAA_i)=\{n_i,\dots,n_i+m_i-1\}$ where $n_i>1$ and $m_i=\#\AAA_i<d-2$.
    \begin{exam} Let $\pi=(7,3,2,6,5,4,1)$. We have that, for $p=2$
	$$ \Omega_\pi=\pmtrx{	0&1&1&1&1&1&1\\
				1&0&1&0&0&0&1\\
				1&1&0&0&0&0&1\\
				1&0&0&0&1&1&1\\
				1&0&0&1&0&1&1\\
				1&0&0&1&1&0&1\\
				1&1&1&1&1&1&0},~\mA_1=\pmtrx{0&1\\1&0},~\mA_2=\pmtrx{0&1&1\\1&0&1\\1&1&0}$$
    So if we write $\pi=\cmtrx{1&2&3&4&5&6&7\\7&3&2&6&5&4&1}$, we can assign to our $\mA_i$'s their corresponding blocks in $\pi$ as follows:
	$$ \mA_1\sim\cmtrx{2&3\\3&2},~\mA_2=\sim\cmtrx{4&5&6\\6&5&4}.$$
    \end{exam}
	
    Because our definitions are invariant under renaming, we have a unique correspondence between a matrix of the form $\mA_i$ and its block in the permutation $\pi=(\pi_0,\pi_1)$.
    \begin{defn} \label{def.what_is_block}
        We allow $\mA_i$ to refer to the matrix and the block in $\pi$, and we shall denote $\mA_i$ as a \term{block} in either case.
    \end{defn}
    We now show the significance of this definition by showing how it aids in determining the spin parity for a standard permutation.
    Because $\Omega_{1,d}=1$, we may choose $\alpha_1=c_1$ and $\beta_1=c_d$. Recalling that $\phi(c_1)=\phi(c_d)=1$, the Equations \eqref{eq.iterate}, for $i,k=2,\dots,d-1$, are now
	\begin{subequations}\label{eq.iterate2}
		\begin{align}
			c_i'&:=c_i+c_1+c_d,\\
			\phi(c_i')&:=1+1+1+1=0,\\
			c_i'\inx c_k'&:=\Omega_{i,k}+1+1=\Omega_{i,k}.
		\end{align}
	\end{subequations}
	The new matrix $\Omega'$ over the remaining $c_i$'s becomes
	\begin{equation}
		\Omega'= \pmtrx{
					\mA_1& \mO& \dots &\dots&\mO\\
					\mO& \mA_2 & \mO&\dots&\vdots\\
					\vdots&\mO&\ddots&\ddots&\vdots\\
					\vdots& \vdots&\ddots&\mA_{p-1}&\mO\\
					\mO&\dots&\dots&\mO&\mA_p		},
	\end{equation}

	and \eqref{eq.spin} becomes
	\begin{equation}\label{eq.spin2}
		\Phi(\pi)=1+\sum_{i=2}^g\phi(\alpha_i)\phi(\beta_i).
	\end{equation}

    Next we notice that if $c_i$ belongs to the block associated to $\mA_j$, which we shall denote as $c_i\in\mA_j$, and $c_k$ is associated to $\mA_m$, $j\neq m$, then $c_i\inx c_k=0$. So for any $\alpha_i$ we select in a given $\mA_j$, $\beta_i$ must also belong to $\mA_j$ as well. So once a pair $\alpha_i,\beta_i\in\mA_j$ has been chosen, for any $c_k\in\mA_m$, $j\neq m$, then by \eqref{eq.iterate},
	\begin{subequations}\label{eq.iterate3}
		\begin{align}
			c_k'&:=c_k,\\
			\phi(c_k')&:=\phi(c_k)=0.
		\end{align}
	\end{subequations}
	Beginning with the initial data in \eqref{eq.iterate2}, we can calculate the value of
	\begin{equation}\label{eq.spinblock}
		\phi(\mA_i):=\sum_{j=2,~\alpha_j,\beta_j\in\mA_i}^{g}\phi(\alpha_j)\phi(\beta_j)
	\end{equation}
	for each $\mA_i$ independently over each other $\mA_j$. So we are lead to our final equation
	\begin{equation}\label{eq.spinf}
		\Phi(\pi)=1+\sum_{i=1}^p \phi(\mA_i).
	\end{equation}
    This final equation is a crucial part to Theorems \ref{thm.even.odd}, \ref{thm.twos.even} and \ref{thm.even.even}.

\sect{Blocks}\label{sec.blocks}

	In this section, we define the necessary blocks to construct the permutations in the following sections.
	These blocks allow us to control the degrees of the singularities as well as the parity of spin.
	
	\begin{defn}\label{def.blocks}
		Let
		$$\SPACE:=\cmtrx{\LL{\alpha}{\alpha}}.$$
		For $n\geq 1$, let
            $$\EVEN_{2n}:=\cmtrx{ \LL{\alpha_1}{\beta_1} ~\LL{\beta_1}{\alpha_1} ~\LL{\alpha_2}{\beta_2} ~\LL{\beta_2}{\alpha_2} ~\LL{\dots}{\dots} ~\LL{\alpha_n}{\beta_n} ~\LL{\beta_n}{\alpha_n} }.$$
		For $n>1$, let
            $$ \ODD_{2n}:=\cmtrx{	\LL{\alpha_1}{\beta_1} ~\LL{\beta_1}{\alpha_1} ~\LL{\alpha_2}{\beta_2} ~\LL{\beta_2}{\alpha_2} ~\LL{\dots}{\dots} ~\LL{\alpha_{n-1}}{\beta_n} ~\LL{\beta_{n-1}}{\alpha_n} ~\LL{\alpha_n}{\beta_{n-1}} ~\LL{\beta_n}{\alpha_{n-1}}}.$$
		Let
            $$ \ODD_{2,2}:=\cmtrx{	\LL{\alpha_1}{\alpha_5} ~\LL{\alpha_2}{\alpha_4} ~\LL{\alpha_3}{\alpha_3} ~\LL{\alpha_4}{\alpha_2} ~\LL{\alpha_5}{\alpha_1}}.$$
		For $m,n\geq 0$, let
            $$ \PAIR_{2m+1,2n+1}:= \cmtrx{\LL{\alpha_1}{\beta_1} ~\LL{\beta_1}{\alpha_1} ~\LL{\dots}{\dots} ~\LL{\alpha_m}{\beta_m} ~\LL{\beta_m}{\alpha_m} ~ \LL{\eps}{\eta} ~\LL{\zeta}{\zeta} ~\LL{\gamma_1}{\delta_1} ~\LL{\delta_1}{\gamma_1} ~\LL{\dots}{\dots} ~\LL{\gamma_n}{\delta_n} ~\LL{\delta_n}{\gamma_n} }.$$
	\end{defn}
	
    For the remainder of the paper, when we speak of concatenating the blocks above, we assume that each block is defined over its own unique subalphabet. For example, $\PAIR_{2m+1, 2n+1}=\EVEN_{2m}\PAIR_{1,1}\EVEN_{2n}$ and $\EVEN_{2(m+n)}=\EVEN_{2m}\EVEN_{2n}$. Also note that all of these blocks contribute a self-inverse portion of a permutation. When we say a block appears \term{inside} a standard permutation $\pi$, we mean that it is a block in $\pi$ and does not include the letters on the outside, i.e. the letters $\pi_0^{-1}(1)=\pi_1^{-1}(d)$ and $\pi_1^{-1}(1)=\pi_0^{-1}(d)$.
    \begin{defn}\label{def.blockConstructed}
        A permutation $\pi\in\irr_d$ is \term{block-constructed} if it is standard and every block that appears inside comes from Definition \ref{def.blocks}. In other words, if $\pi$ is block constructed then
        $$ \pi=\stdperm{\mB_1\cdots\mB_k}$$
        where each $\mB_i$ is from Definition \ref{def.blocks}.
    \end{defn}
    We now show the desired properties of the defined blocks.

	\begin{lem}\label{lem.blocks}
        For the blocks in Definition \ref{def.blocks}, assuming they appear inside a standard permutation, the following are true:
		\begin{itemize}
			\item Assuming its leftmost (or rightmost) top and bottom singularities are identified,
				$\EVEN_{2n}$ contributes a singularity of degree $2n$ and $\phi(\EVEN_{2n})=0$.
			\item Assuming its leftmost (or rightmost) top and bottom singularities are identified,
				$\ODD_{2n}$, $n>1$, contributes a singularity of degree $2n$ and $\phi(\ODD_{2n})=1$.
			\item $\ODD_{2,2}$ contributes two singularities of degree $2$ and $\phi(\ODD_{2,2})=1$.
			\item $\PAIR_{2m+1,2n+1}$ contributes two singularities, one of degree $2m+1$ and one of degree $2n+1$.
		\end{itemize}
		Also, any block-constructed permutation is its own inverse.
	\end{lem}
	
	\begin{proof}
		We make note of a few relationships between our defined blocks:
		\begin{equation*}
			\begin{split}
				\EVEN_{2(n+1)}&=\EVEN_2\EVEN_{2n},\\
				\ODD_{2(n+1)}&=\EVEN_2\ODD_{2n},\\
				\PAIR_{2m+1,2n+1}&=\EVEN_{2m}\PAIR_{1,1}\EVEN_{2n}.
			\end{split}
		\end{equation*}
		
        We first prove the statement for $\EVEN_2$. This block has the structure in the suspended surface for $\pi$ as seen in Figure \ref{fig.blocks}. Because the two singularities are identified by assumption, this is one singularity of degree $2$. To calculate $\phi(\EVEN_2)$, we observe that the matrix associated to $\EVEN_2$ is just $\pmtrx{0&1\\1&0}$. So we choose $\alpha_1=c_1,\beta_1=c_2$ as our \term{canonical basis}. Using Equations \eqref{eq.iterate2} and \eqref{eq.spinblock},
		$$\phi(\EVEN_2)=\phi(\alpha_1)\phi(\beta_1)=0.$$

        \begin{figure}[h]
            \begin{center}
               \setlength{\unitlength}{350pt}
                \begin{picture}(1,.49)
                    \put(0,0){\includegraphics[width=\unitlength]{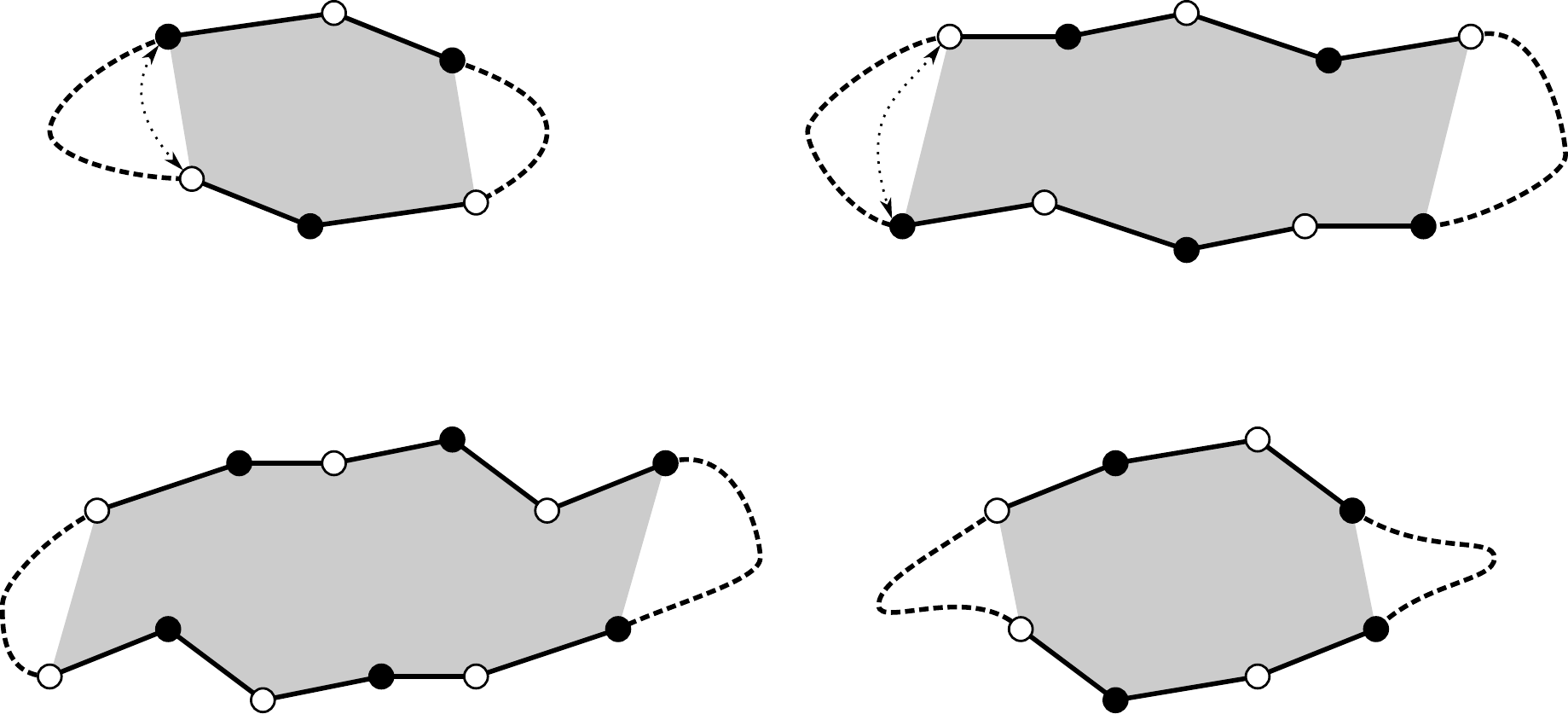}}

                    \put(.14,.46){$\alpha_1$} \put(.24,.46){$\beta_1$}
                    \put(.14,.29){$\beta_1$} \put(.24,.29){$\alpha_1$}
                    \put(.18,.37){$\EVEN_2$}

                    \put(.63,.46){$\alpha_1$} \put(.71,.47){$\beta_1$} \put(.79,.46){$\alpha_2$} \put(.87,.46){$\beta_2$}
                    \put(.61,.28){$\beta_2$} \put(.69,.28){$\alpha_2$} \put(.79,.26){$\beta_1$} \put(.86,.28){$\alpha_1$}
                    \put(.73,.36){$\ODD_{4}$}

                    \put(.08,.16){$\alpha_1$} \put(.16,.18){$\alpha_2$} \put(.235,.19){$\alpha_3$} \put(.32,.17){$\alpha_4$} \put(.38,.17){$\alpha_5$}
                    \put(.05,0){$\alpha_5$} \put(.11,0){$\alpha_4$} \put(.195,-.02){$\alpha_3$} \put(.26,-.01){$\alpha_2$} \put(.34,.01){$\alpha_1$}
                    \put(.19,.08){$\ODD_{2,2}$}

                    \put(.66,.15){$\eps$} \put(.75,.19){$\zeta$} \put(.83,.17){$\eta$}
                    \put(.66,.00){$\eta$} \put(.75,-.02){$\zeta$} \put(.84,.01){$\eps$}

                    \put(.73,.08){$\PAIR_{1,1}$}
                \end{picture}
            \end{center}
            \caption{The blocks $\EVEN_2$, $\ODD_4$, $\ODD_{2,2}$ and $\PAIR_{1,1}$ in suspensions. The two singularities in the top surfaces are identified by assumption in Lemma \ref{lem.blocks}.}\label{fig.blocks}
        \end{figure}
		
        For $\EVEN_{2n}$, see that this is nothing more than $n$-$\EVEN_2$ blocks concatenated, each block contributing a singularity of degree $2$. Because of concatenation, these singularities are all identified to form one of degree $2n$. To calculate $\phi(\EVEN_{2n})$, notice that its matrix is just $n$ $\EVEN_2$ blocks along the diagonal. So by our reasoning in Section \ref{sec.stdspin},
    		$$ \phi(\EVEN_{2n})=\sum_{i=1}^n \phi(\EVEN_2)=0.$$
		
        We now prove the claim for $\ODD_4$. This block has the structure in the surface as indicated by Figure \ref{fig.blocks}. As in the case for $\EVEN_2$, it is clear that this contributes a singularity of degree $4$. The matrix for $\ODD_4$ is
    		$$ \ODD_4=\pmtrx{0&1&1&1\\1&0&1&1\\1&1&0&1\\1&1&1&0}$$
        We now calculate $\phi(\ODD_4)$. Using again Equations \eqref{eq.iterate2} and \eqref{eq.spinblock} with choices $\alpha_1=c_1$, $\beta_1=c_4$, $\alpha_2=c_2'$, and $\beta_2=c_3'$,	
    		$$ \phi(\ODD_4)=\phi(\alpha_1)\phi(\beta_1)+\phi(\alpha_2)\phi(\beta_2)=0+1=1.$$

        For $\ODD_{2n}$, we note that it is $(n-2)$-copies of $\EVEN_2$ followed by $\ODD_4$. As above the degree of the singularity is then $(n-2)\cdot 2+4=2n$. The matrix of $\ODD_{2n}$ is $(n-2)$-$\EVEN_2$ blocks and one $\ODD_4$ block along the diagonal. So
    		$$ \phi(\ODD_{2n})=\sum_{i=1}^{n-2}\phi(\EVEN_2)+\phi(\ODD_4)=0+1=1$$
		
        We now prove the theorem for $\ODD_{2,2}$. This follows from the block's portion in the surface (see Figure \ref{fig.blocks}). The matrix for $\ODD_{2,2}$ is
    		$$ \ODD_{2,2}=\pmtrx{0&1&1&1&1\\1&0&1&1&1\\1&1&0&1&1\\1&1&1&0&1\\1&1&1&1&0}$$
		all of the statements in the theorem for $\ODD_{2,2}$ follow immediately as they have for the previous blocks.\\
		
        For $\PAIR_{2m+1,2n+1}$, it suffices to prove the statement for $\PAIR_{1,1}$, as $\PAIR_{2m+1,2n+1} = \EVEN_{2m}\PAIR_{1,1}\EVEN_{2n}$, making the singularities degrees $2m+1$ and $2n+1$ as desired. The portion of the surface determined by $\PAIR_{1,1}$ is shown in Figure \ref{fig.blocks}. So again, counting verifies that there are two singularities with degree $1$.\\
		
		The final statement of the theorem is clear as each block places all of its letters in self-inverse positions and the
		outside $\lA$ and $\lB$ letters (making the permutation standard) are in self-inverse position as well.
	\end{proof}
	
	Before proving the main theorem, we remark that the block $\SPACE$ is designed to keep singularities
	of neighboring blocks separate. To illustrate this point, notice that $\EVEN_{2n}\EVEN_{2m}$ contributes
	\term{one} singularity of degree $2(m+n)$, while $\EVEN_{2m}\SPACE\EVEN_{2n}$ contributes the desired
	\term{two} singularities. The block $\SPACE$ also causes any neighboring block's leftmost (or rightmost)
	top and bottom singularities to be identified, as required in Lemma \ref{lem.blocks}.

\sect{Self-Inverses for $g\leq 3$}\label{sec.gleq3}

    \begin{thm}\label{thm.g3}
            Given $\tilde{\pi}\in\irr$ such that $g(\tilde{\pi})\leq 3$, There exists $\pi\in\RClass(\tilde{\pi})$ such that $\pi=\pi^{-1}$.
	\end{thm}

    We shall prove this result simply by stating such an element for each class according to the following tables. When possible, we use the construction methods in the theorems for higher genera to make our example. To consider additional removable singularities, refer to Theorem \ref{thm.zero}.\\

	\begin{center}
        \begin{tabular}{c}
		Genus 1 \\
		\tbl {|c |c| c|}{
			\hline
			Signature & Type & Self-inverse \\
			\hline
			$(0)$ 	& hyperelliptic	& $(2,1)$\\
			\hline
			$(0,0)$	& hyperelliptic	& $(3,2,1)$\\
			\hline }
        \end{tabular}
	\end{center}

	\begin{center}
        \begin{tabular}{c}
		Genus 2 \\
		\tbl {|c |c| c|}{
			\hline
			Signature & Type & Self-inverse \\
			\hline	
			$(2)$ 	& hyperelliptic	& $(4,3,2,1)$\\
			\hline
			$(1,1)$ & hyperelliptic	& $(5,4,3,2,1)$\\
			\hline }
        \end{tabular}
	\end{center}

	\begin{center}
        \begin{tabular}{c}
		Genus 3 \\
		\tbl{|c |c| c|}{
			\hline
			Signature & Type & Self-inverse \\
			\hline
			$(4)$ & hyperelliptic	& $(6,5,4,3,2,1)$\\
			\hline
			$(4)$ & odd	& $(6,3,2,5,4,1)$\\
			\hline
			$(3,1)$ & -	& $(7,4,3,2,6,5,1)$\\
			\hline
			$(1,3)$ & -	& $(7,3,2,6,5,4,1)$\\
			\hline
			$(2,2)$ & hyperelliptic	& $(7,6,5,4,3,2,1)$\\
			\hline
			$(2,2)$ & odd	& $(7,3,2,4,6,5,1)$\\
			\hline
			$(1,1,2)$ & -	& $(8,3,2,4,7,6,5,1)$\\
			\hline
			$(2,1,1)$ & -	& $(8,4,3,2,5,7,6,1)$\\
			\hline
			$(1,1,1,1)$ & -	& $(9,4,3,2,5,8,7,6,1)$\\
			\hline }
        \end{tabular}
	\end{center}

\sect{Self-Inverses for $g\geq 4$}\label{sec.ggeq4}

	\begin{thm}\label{thm.odds} Let Rauzy Class $\RClass$ have signature $(\ell_1,\dots,\ell_m)$ such that $\ell_i$ is odd for
		some $i$. Then there exists $\pi\in\RClass$ such that $\pi=\inv{\pi}$.
	\end{thm}
	
	\begin{proof}
        We shall give an explicit construction of such a $\pi$. As there are singularities of odd degree, we only need to verify that our constructed permutation has the appropriate signature. We do this considering two cases: $\ell_1$ is even or $\ell_1$ is odd.\\
		
        First assume that $\ell_1$ is odd. Then we can rearrange our $\ell_i$'s such that $\ell_i$ is odd for $1\leq i\leq k$ and $\ell_i$ is even for $k<i\leq m$. Notice that $k$ must be even as the sum over all $\ell_i$'s is even. So we define $\pi$ by
            $$ \pi=\stdperm{ \EVEN_{\ell_m} \SPACE \EVEN_{\ell_{m-1}} \SPACE \cdots \SPACE \EVEN_{\ell_{k+1}} \SPACE \PAIR_{\ell_k,\ell_{k-1}} \SPACE \cdots \SPACE \PAIR_{\ell_2,\ell_1}}.$$
        By Lemma \ref{lem.blocks}, this permutation has the appropriate singularities. Since the singularity of degree $\ell_1$ is the one immediately to the left of $\lA$ and $\lB$, it is the marked singularity.\\
		
        The second case is to assume that $\ell_1$ is even. We then make a division such that $\ell_i$ is even for all $1\leq i\leq k$ and is odd for $k< i \leq m$, noticing that this time $m-k$ must be even. Then our desired $\pi$ is
            $$ \pi=\stdperm{\PAIR_{\ell_m,\ell_{m-1}} \SPACE \PAIR_{\ell_{m-2},\ell_{m-3}} \SPACE \cdots \SPACE \PAIR_{\ell_{k+2},\ell_{k+1}} \SPACE \EVEN_{\ell_k} \SPACE \cdots \SPACE \EVEN_{\ell_1}}.$$
		Just as in the previous case, this has the desired signature.
	\end{proof}

	\begin{thm}\label{thm.even.odd}
        Let Rauzy Class $\RClass$ have signature $(2\ell_1,\dots,2\ell_m)$ and odd spin. Then there exists $\pi\in\RClass$ such that $\pi=\inv{\pi}$.
	\end{thm}
	
	\begin{proof}
        As opposed to the proof of Theorem \ref{thm.odds}, we must construct a permutation that not only has the appropriate signature $(2\ell_1,\dots,2\ell_m)$ but also satisfies $\Phi(\pi)=1$. Let $\pi$ be defined as
    		$$ \pi = \stdperm{\EVEN_{2\ell_m} \SPACE \EVEN_{2\ell_{m-1}} \SPACE \cdots \SPACE \EVEN_{2\ell_2} \SPACE \EVEN_{2\ell_1}}.$$
		Again, as in the proof of Theorem \ref{thm.odds}, this has the desired signature. By Equation \eqref{eq.spinf},
		$$ \Phi(\pi) = 1 + \sum_{i=1}^{m}\phi(\EVEN_{2\ell_i})=1$$
		as $\phi(\EVEN_{2\ell_i})=0$ by Lemma \ref{lem.blocks}.
	\end{proof}

	\begin{thm}\label{thm.twos.even}
        Let Rauzy Class $\RClass$ have signature $(2,\dots,2)$ and even spin. Then there exists $\pi\in\RClass$ such that $\pi=\inv{\pi}$.
	\end{thm}
	
	\begin{proof}
        We will construct our desired $\pi$, show that it has the appropriate signature and verify that $\Phi(\pi)=0$. Let $m>1$ be the number singularities of degree $2$. Then we may define $\pi$ as
            $$ \pi = \stdperm{\mB_{m-1} \SPACE \mB_{m-2} \SPACE \cdots \SPACE \mB_1},~\mB_i =\RHScase{\ODD_{2,2},& i=1 \\ \EVEN_2,& \mathrm{otherwise.}}$$
        By Lemma \ref{lem.blocks}, this has the appropriate signature. We also know that $\phi(\EVEN_2)=0$ and $\phi(\ODD_{2,2})=1$. So by \eqref{eq.spinf}
    		$$ \Phi(\pi)=1+\sum_{i=1}^{m-1}\phi(\mB_i)=1+1=0.$$
	\end{proof}

	\begin{thm}\label{thm.even.even}
         Let Rauzy Class $\RClass$ have even spin and signature $(2\ell_1,\dots,2\ell_m)$ such that $\ell_i>1$ for some $i$. Then there exists $\pi\in\RClass$ such that $\pi=\inv{\pi}$.
	\end{thm}
	
	\begin{proof}
        We must again construct a $\pi$ with signature $(2\ell_1,2\ell_2,\dots,2\ell_m)$ and such that $\Phi(\pi)=0$. Let $j$ be chosen such that $\ell_j>1$. Then we define $\pi$ as
            $$ \pi = \stdperm{\mB_m \SPACE \mB_{m-1} \SPACE \cdots \SPACE \mB_1},~\mB_i = \RHScase{\ODD_{2\ell_i}& \mbox{if }i=j, \\ \EVEN_{2\ell_i}& \mathrm{otherwise.}}$$
        By Lemma \ref{lem.blocks}, this has the appropriate signature. We know that $\phi(\EVEN_{2\ell_i})=0$ and $\phi(\ODD_{2\ell_j})=1$. So by Equation \eqref{eq.spinf},
    		$$ \Phi(\pi)=1+\sum_{i=1}^m\phi(\mB_i)=1+1=0.$$	
	\end{proof}

\sect{Self-inverses with Removable Singularities}\label{sec.removable}

    Singularities are called removable if they have degree $0$. A suspension of a permutation with removable singularities will move by induction to consider such singularities, but they are not actually zeroes of the corresponding Abelian differential. The results for Corollary \ref{cor.classify} extend to Rauzy Classes with removable singularities. Given a permutation $\pi$ with removable singularities, we first must consider which connected component of which stratum $\RClass(\pi)$ belongs to, called $\CCC$. We then choose a Rauzy Class $\RClass'\subseteq \CCC$ in such that:
    \begin{itemize}
        \item If the genus of $\RClass$ is $1$, then we define $\RClass'$ to be the irreducible permutation on $2$ letters.
        \item If the genus of $\RClass$ is greater than $1$, no singularity of $\RClass'$ is removable.
        \item If the marked singularity of $\RClass$ is not removable and of degree $n$, then the marked singularity of $\RClass'$ is of degree $n$ as well. If the marked singularity of $\RClass$ is removable, $\RClass'$ has no restriction on which singularity is marked.
    \end{itemize}
    In another context, consider a suspension with differential, $(S,q)$, on any representative $\pi\in\RClass$. We complete the differential $q$ at any removable singularity and call this new differential $q'$ on the same surface. As long as the marked singularity of $\pi$ wasn't removable, then the same singularity is marked by $\pi'$, the permutation resulting from $(S,q')$. If the marked singularity of $\pi$ was removable, then we may choose a new marked singularity of what remains in $(S,q')$ and define $\pi'$ by this choice.

    \begin{defn}
        We call such $\RClass'$ an \term{underlying Rauzy Class}, or underlying class, of $\RClass$.
    \end{defn}

    \begin{exam}
        For $\pi=(7,4,5,2,6,3,1)$, $\RClass=\RClass(\pi)$ has signature $(4,0)$ and odd spin (non hyperelliptic). So we need to find $\RClass'\subseteq\Hodd(4)$ with no removable singularities. Our only choice is, by theorem \ref{thm.even.odd}, $\RClass'=\RClass(\pi')$ for $\pi'=(6,3,2,5,4,1)$.
    \end{exam}

    \begin{exam}
        For $\pi=(7,6,1,4,3,2,8,5)$, $\RClass=\RClass(\pi)$ has signature $(0,3,1)$ so it belongs to the connected stratum $\HHH(3,1)$. There are two choices of underlying Rauzy classes, $\RClass'=\RClass(7,4,3,2,6,5,1)$ and $\RClass''=\RClass(7,3,2,6,5,4,1)$, based on the choice of the new marked singularity.
    \end{exam}

    \begin{rem}
        If $\RClass$ has a marked singularity that is not removable, the choice of $\RClass'$ is unique. One can also check that if $\RClass$ is a Rauzy class on $d$ letters with $k$ removable singularities, an underlying class $\RClass'$ is a Rauzy class on $d-k$ letters.
    \end{rem}

    \begin{thm}\label{thm.zero}
        Given a Rauzy Class $\RClass$ with at least one removable singularity, there exists $\pi\in\RClass$ such that $\pi=\inv{\pi}$.
    \end{thm}

    \begin{proof}
        By the above discussion and Theorems \ref{thm.g3}-\ref{thm.even.even}, we have a $\tilde{\pi}=\inv{\tilde{\pi}}$ that belongs to the underlying class $\RClass'$ on $d-k$ letters. Denote this by
            $$ \tilde{\pi}=\stdperm{\BLAH}$$
        We now only need to confirm two cases, either $\ell_1=0$ or $\ell_1\neq 0$ for $\sig(\RClass)=(\ell_1,\dots,\ell_m)$. In the first case, we have the permutation
            $$ \pi=\stdperm{ \BLAH \underbrace{\SPACE\cdots\SPACE}_k}.$$
        In the second case, the permutation is of the form
            $$ \pi=\stdperm{\underbrace{\SPACE\cdots\SPACE}_k \BLAH}.$$
    \end{proof}

\begin{figure}[h]
        \begin{center}
           \setlength{\unitlength}{350pt}
            \begin{picture}(1,.23)
                \put(0,0){\includegraphics[width=\unitlength]{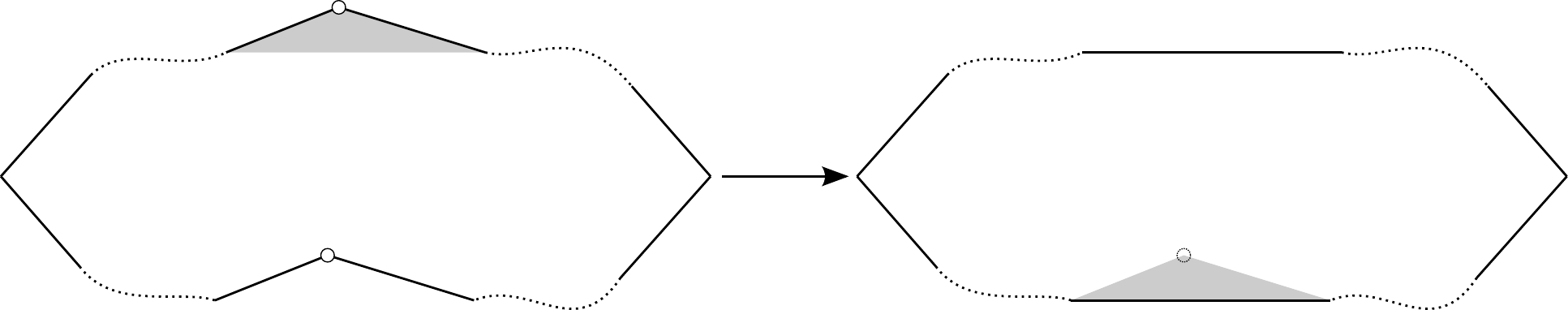}}

                \put(.01,.13){$\alpha$} \put(.16,.19){$\gamma$} \put(.26,.19){$\delta$} \put(.425,.13){$\beta$} \put(.005,.02){$\beta$} \put(.16,-.01){$\gamma$} \put(.25,-.02){$\delta$} \put(.42,.02){$\alpha$}

                \put(.20,.08){$S$} \put(.75,.08){$S'$}

                \put(.56,.13){$\alpha$} \put(.76,.18){$\gamma$} \put(.555,.02){$\beta$} \put(.975,.13){$\beta$} \put(.76,-.02){$\gamma$} \put(.97,.02){$\alpha$}
           \end{picture}
       \end{center}
       \caption{The differential for $S'$ completes the removable singularity in the differential for $S$. All other singularities remain unchanged.}\label{fig.remove_ex}
    \end{figure}
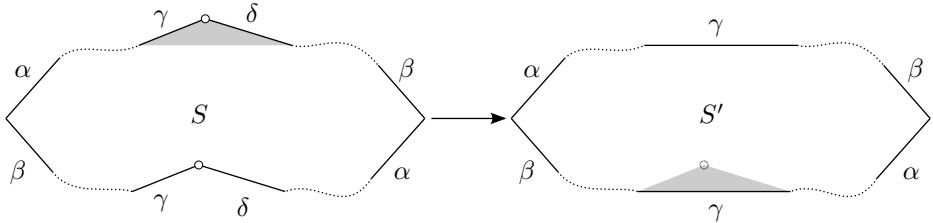

    \begin{rem}\label{rem.reduce_removable}
        Consider any permutation $\tilde \pi$. Move to standard $\pi\in\RClass(\tilde \pi)$ (see Claim \ref{cor.std_in_RC}). If $\pi$ has removable singularities then one of the following must be satisfied:
        \begin{enumerate}
            \item There exists $\gamma,\delta\in\AAA$ such that $\pi_\eps(\gamma)+1=\pi_\eps(\delta)$ for each $\eps\in\{0,1\}$.
                Consider the map
                $$ \pi = \cmtrx{\LL{\alpha}{\beta}~\LL{\dots\gamma\delta}{\dots}~\LL{\dots}{\gamma\delta\dots}~\LL{\beta}{\alpha}}
                    \to
                \cmtrx{\LL{\alpha}{\beta}~\LL{\dots\gamma}{\dots}~\LL{\dots}{\gamma\dots}~\LL{\beta}{\alpha}} = \pi'$$
                that ``forgets" $\delta$. Figure \ref{fig.remove_ex} shows how this map eliminates the removable singularity between $\gamma$ and $\delta$.
            \item There exists $\gamma\in\AAA$ such that $\pi_0(\gamma)=\pi_1(\gamma)=2$. The map
                $$ \pi = \cmtrx{\LL{\alpha}{\beta}~\LL{\gamma}{\gamma}~\LL{\delta\dots}{\eta\dots}~\LL{\beta}{\alpha}}
                    \to
                \cmtrx{\LL{\alpha}{\beta}~\LL{\delta\dots}{\eta\dots}~\LL{\beta}{\alpha}}=\pi'$$
                ``forgets" $\gamma$. This map eliminates the removable singularity to the left of the segments labeled $\gamma$.
            \item There exists $\gamma\in\AAA$ such that $\pi_0(\gamma)=\pi_1(\gamma) = d-1$. Consider the map
                $$ \pi = \cmtrx{\LL{\alpha}{\beta}~\LL{\dots\delta}{\dots\eta}~\LL{\gamma}{\gamma}~\LL{\beta}{\alpha}}
                    \to
                \cmtrx{\LL{\alpha}{\beta}~\LL{\dots\delta}{\dots\eta}~\LL{\beta}{\alpha}}=\pi'$$
                that ``forgets" $\gamma$. This new permutation no longer has a removable marked singularity. The new marked singularity was originally on the left of the segments labeled $\gamma$.
        \end{enumerate}
        By performing these maps, we may explicitly derive a standard representative of $\RClass'$.
    \end{rem}

    \begin{exam}
        Begin with $\pi = (7,4,5,1,6,2,3)$, with signature $\sig(\pi)=(2,0,0,0)$. We then consider standard $\pi'=\rb^3\pi = (7,6,2,3,4,5,1)$. By performing the first reduction listed in Remark \ref{rem.reduce_removable}, we get
        $$ \cmtrx{\LL{1}{7}\LL{2}{6}\LL{3}{2}\LL{4}{3}\LL{5}{4}\LL{6}{5}\LL{7}{1}}
            \to \cmtrx{\LL{1}{7}\LL{2}{6}\LL{4}{2}\LL{5}{4}\LL{6}{5}\LL{7}{1}}
            \to \cmtrx{\LL{1}{7}\LL{2}{6}\LL{5}{2}\LL{6}{5}\LL{7}{1}}
            \to \cmtrx{\LL{1}{7}\LL{2}{6}\LL{6}{2}\LL{7}{1}} = (4,3,2,1) = \tilde\pi.$$
        In this case, $\tilde\pi$ is in the hyperelliptic class and is self-inverse with one singularity of degree $2$. So by the proof of Theorem \ref{thm.zero}, we derive,
        $$ \pi'' = \cmtrx{\LL{1}{4}~\LL{a}{a}~\LL{b}{b}~\LL{c}{c}~\LL{2}{3}~\LL{3}{2}~\LL{4}{1}}=(7,2,3,4,6,5,1).$$
        One may verify that $\pi'' = \rt^3\rb^2\rt^2\rb\pi$. Therefore $\pi''\in\RClass(\pi)$ and is self-inverse.
    \end{exam} 
\chap{Generalized Rauzy Classes}\label{chapGen}

    \begin{thm}\label{thm.q.main}
        (Main Result) A generalized Rauzy Class $\RClass\subseteq\genirr_\AAA$ with signature $(\ell_1,\dots,\ell_m)$ contains a permutation $\pi$ such that $\pi=\inv{\pi}$ if and only if the following conditions hold:
        \begin{subequations}\label{eq.q.nec}
            \begin{equation}\label{eq.q.nec1}
                \ell_1=0,
            \end{equation}
            \begin{equation}\label{eq.q.nec2}
                \#\{i:\ell_i=j,~1\leq i \leq m\}\mbox{ is even}, \mbox{ for all odd }j\geq-1.
            \end{equation}
        \end{subequations}
    \end{thm}

\sect{Necessary Conditions}

    \begin{defn}
        A generalized permutation $\pi$ is \term{balanced} if $d_0=d_1$, i.e. if the top and bottom rows are the same length.
    \end{defn}

    Let $\AAA$ be an alphabet on $d$ letters. Let $h:\dset{2d}\to\dset{2d}$ be the involution defined by
    $$ h(i) := \RHScase{i+d, & \mbox{if }1\leq i\leq d, \\ i-d, & \mbox{if }d+1\leq i \leq 2d.}$$
    A self-inverse $\pi\in\genirr_\AAA$ must be balanced and satisfy $\pi = \pi\circ h$ up to renaming.

    \begin{defn}
        Letters $\alpha,\beta\in\AAA$ are \term{paired} in self-inverse $\pi\in\genirr_\AAA$ if $h(\inv\pi(\alpha)) = \inv\pi(\beta)$. In other words, $\alpha$ and $\beta$ appear in the same columns in our representation of $\pi$ (it is possible that $\alpha=\beta$).
    \end{defn}

    \begin{exam}
        For
        $$ \pi=\cmtrx{\LL{a}{d}~\LL{b}{e}~\LL{c}{c}~\LL{b}{e}~\LL{d}{a}}$$
        the letters $a$ and $d$ are paired, as are $b$ and $e$. The letter $c$ is paired with itself.
    \end{exam}

    If $\pi\in\genirr_\AAA$ is balanced, $\mone:=(1,\dots,1)\in\RR_+^{\AAA,\pi}$.
    \begin{defn}
        For a balanced $\pi\in\genirr_\AAA$, consider any $\tau\in\TTT_\pi$. Denote
        $$ S_1=S(\pi,\mone,\tau)$$
        as a \term{unit suspension} for $\pi$.
    \end{defn}

    \begin{rem}
        Given such a $\pi$ and $S_1$, every vertical geodesic is either a closed saddle connection (a geodesic with singularities at each endpoint) or a closed loop avoiding any singularities. This fact makes the following definition possible.
    \end{rem}

    \begin{defn}
        Let $\pi\in\genirr_\AAA$ be balanced and $S_1$ be a unit suspension. A singularity represented by vertex class $[v]$ is \term{paired} to singularity $[v']$ if there exists a vertical saddle connection between $[v]$ and $[v']$ ($[v]=[v']$ is possible).
    \end{defn}

    \begin{lem}\label{lem.q.pairs}
        Let $\pi\in\genirr_\AAA$ be self-inverse with unit suspension $S_1$. Then every singularity $[v]$ is paired to a unique singularity $[v']$.
    \end{lem}

    \begin{proof}
        Consider any occurrence of $v$ in $S_1$ and its paired vertex $v'$. The other vertices in the class of $v$ are found by identifying $v$ by the identified edges connected to $v$. Let $\alpha$ represent an edge connected to $v$, and let $\beta$ be the paired letter with edge attached to $v'$ (see Section \ref{sec.gen.surfaces}). Then the occurrence of $v$ attached to the other edge labeled $\alpha$ is paired to a vertex. However, this vertex must be an occurrence of $v'$ as it is attached to the second edge labeled $\beta$. We conclude if $v$ and $v'$ are endpoints of one vertical saddle connection in $S_1$, every vertical saddle attached to any copy of $v$ must have a copy of $v'$ as its other endpoint.
    \end{proof}

    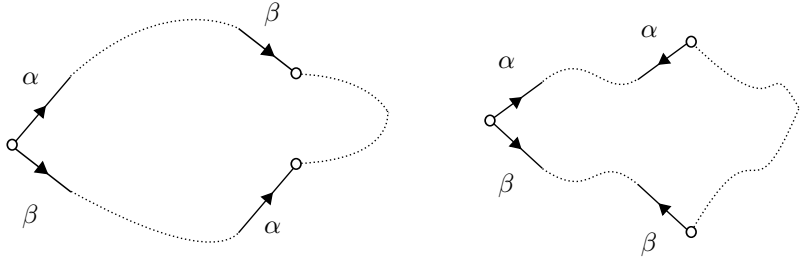
\begin{figure}[h]
        \begin{center}
            \setlength{\unitlength}{300pt}
            \begin{picture}(1,.35)
                \put(0,0){\includegraphics[width=\unitlength]{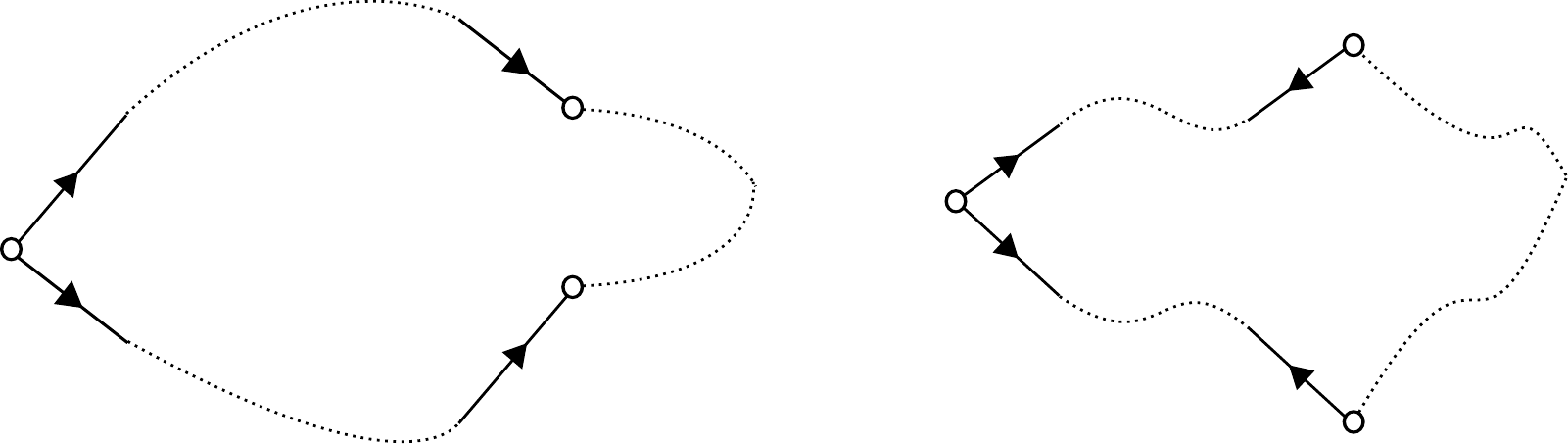}}
                \put(.02,.025){$\beta$}
                \put(.02,.2){$\alpha$}
                \put(.325,.28){$\beta$}
                \put(.325,.0125){$\alpha$}

                \put(.62,.22){$\alpha$}
                \put(.62,.064){$\beta$}
                \put(.8,.26){$\alpha$}
                \put(.8,-.01){$\beta$}
            \end{picture}
        \end{center}
        \caption{The marked singularity of self-inverse $\pi$ must be paired with itself.}\label{fig.marked_paired}
    \end{figure}

    \begin{lem}\label{lem.q.nec}
        Let $\pi\in\genirr$ with signature $\sigma=(\ell_1,\dots,\ell_m)$. If $\pi=\inv{\pi}$, then Equations \eqref{eq.q.nec1} and \eqref{eq.q.nec2} hold.
    \end{lem}

    \begin{proof}
        Let $\alpha = \pi(1)$ and $\beta = \pi(d+1)$ be the first letters in each row and $v_1$ the leftmost singularity. Because $\alpha$ and $\beta$ are paired, there exists another appearance of $v_1$ attached to $\alpha$. However, $\beta$ also is attached to an occurrence of $v_1$, implying that $v_1$ is paired to itself, see Figure \ref{fig.marked_paired}. By Lemma \ref{lem.q.pairs} and counting, we conclude that $\ell_1$ is even.\\

        Consider any odd $k\geq -1$ and $j$ such that $\ell_j=k$ corresponding to vertex $v_j$. By our previous observation, $v_j$ can not be paired to itself. So there exists unique singularity $v_{j'}$ paired to $v_j$. By this we can conclude that $\ell_j=\ell_{j'}$. Because every odd ordered singularity has a unique and distinct paired singularity of the same order, \eqref{eq.q.nec2} follows.
    \end{proof}

    \begin{rem}\label{rem.q.signatureform}
        For the remainder of this \chapword, we assume that Equations \eqref{eq.q.nec} hold for any generalized Rauzy Class $\RClass$. As a result, we assume that $\RClass$ has a signature of the form $(2\ell_1,\dots,2\ell_m,j_1^2,\dots,j_n^2,-1^p)$ for $\ell_i\geq 0$, $j_i\geq 1$ odd, and $p$ the number of poles.
    \end{rem}

\sect{Blocks and Insertions}\label{sec.gen.blocks}

    We recall Definition \ref{def.what_is_block} and continue to use the term \term{block} to refer a pattern of letters that appear in a balanced permutation. In this section, we define blocks needed to construct self-inverses in Sections \ref{sec.gen.g01}, \ref{sec.gen.g2hyp} and \ref{sec.gen.g2}.

    \begin{defn}
        Fix balanced $\pi\in\genirr_\AAA$ and let $\alpha_i=\pi(i)$, $i\in\dset{2d}$. Consider paired singularities $[v]$ and $[v']$ and $i\in\dset{d-1}$, $i$, such that:
        \begin{itemize}
            \item The vertex $v$ in between edges $\alpha_i$ and $\alpha_{i+1}$ belongs to $[v]$.
            \item The vertex $v'$ in between edges $\alpha_{i+d}$ and $\alpha_{i+d+1}$ belongs to $[v']$.
        \end{itemize}
        For any block of $d'$ letters $\mB$ over alphabet $\AAA'$, let $(\pi',d+d',d+d')$ be a generalized permutation over alphabet $\AAA\sqcup\AAA'$ defined by
        $$ \pi' = \cmtrx{\LL{\alpha_1}{\alpha_{d+1}} ~\LL{\dots}{\dots} ~\LL{\alpha_i}{\alpha_{i+d}} ~\LL{b_1}{b_{d'+1}} ~\LL{\dots}{\dots} ~\LL{b_{d'}}{b_{2d'}} ~\LL{\alpha_{i+1}}{\alpha_{i+d+1}} ~\LL{\dots}{\dots} ~\LL{\alpha_d}{\alpha_{2d}}}$$
        where
        $$\mB=\cmtrx{\LL{b_1}{b_{d'+1}} ~\LL{\dots}{\dots} ~\LL{b_{d'}}{b_{2d'}}}.$$
        We say that $\pi'$ is made by an \term{insertion} of $\mB$ at paired singularities $[v]$ and $[v']$.
    \end{defn}

    \begin{defn}\label{def.POLES}
        Let $\POLES_0$ be defined to be empty. We then define $\POLES_i$, $i\geq 1$, recursively by:
        \begin{itemize}
            \item $\POLES_1:=\cmtrx{\LL{\alpha}{\beta}&\LL{\alpha}{\beta}}$ and
            \item $\POLES_{i+1}:=\POLES_1\POLES_i$.
        \end{itemize}
    \end{defn}

    \begin{lem}\label{lem.POLES}
        Let $\pi\in\genirr_\AAA$ be self-inverse with fixed paired singularities of odd order $j\geq -1$. When block $\POLES_{2i}$ is inserted at that singularity pair, the resulting permutation $\pi'$ has replaced the fixed singularities with paired singularities of odd order $j+2i$ and has $4i$ more poles than $\pi$. When block $\POLES_i$ is inserted at an even ordered singularity of order $2\ell\geq 0$, the resulting permutation $\pi'$ has replaced the fixed singularity with one of order $2(\ell+i)$ and has $2i$ more poles than $\pi$. In either case, $\pi'$ is self-inverse and $\pi'\in\genirr_\AAA$.
    \end{lem}

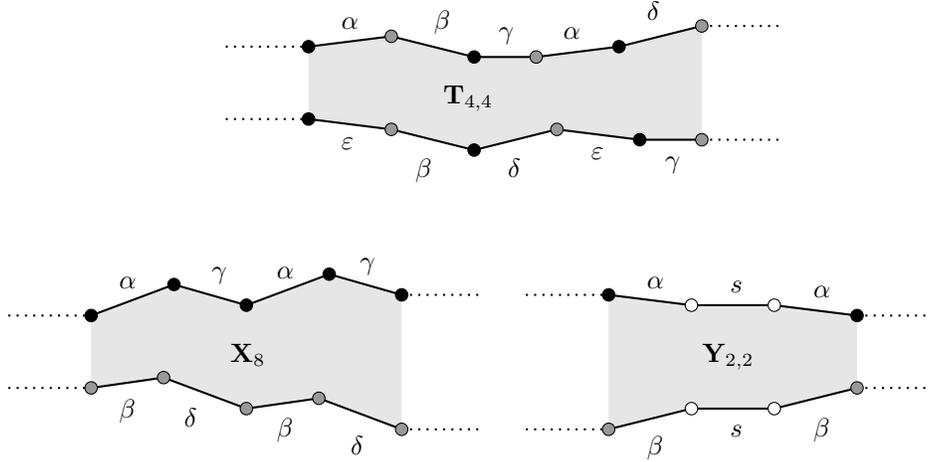
\begin{figure}[h]
        \begin{center}
           \setlength{\unitlength}{350pt}
            \begin{picture}(1,.48)
                \put(0,0){\includegraphics[width=\unitlength]{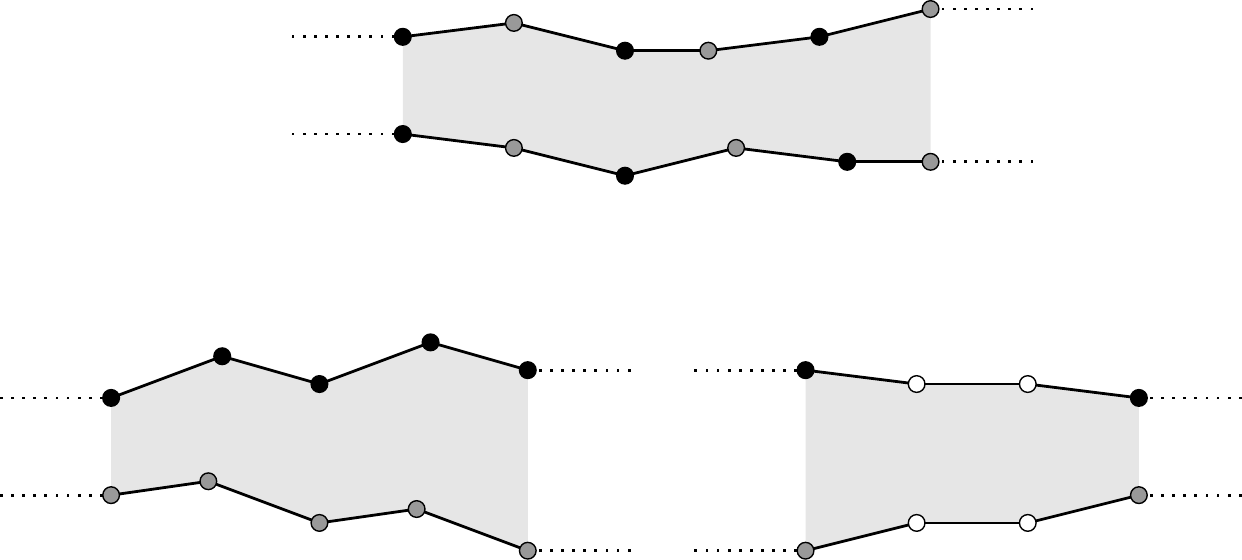}}

                \put(.47,.36){$\TWOFOURS_{4,4}$}
                \put(.24,.08){$\JUSTFOURS_8$}
                \put(.75,.08){$\TWOTWOS_{2,2}$}

                \put(.36,.44){$\alpha$} \put(.46,.44){$\beta$} \put(.53,.43){$\gamma$} \put(.60,.43){$\alpha$} \put(.69,.45){$\delta$}
                \put(.36,.31){$\eps$} \put(.44,.28){$\beta$} \put(.54,.28){$\delta$} \put(.63,.30){$\eps$} \put(.71,.29){$\gamma$}

                \put(.12,.16){$\alpha$} \put(.22,.17){$\gamma$} \put(.29,.17){$\alpha$} \put(.38,.18){$\gamma$}
                \put(.12,.02){$\beta$} \put(.19,.01){$\delta$} \put(.29,.00){$\beta$} \put(.37,-.02){$\delta$}

                \put(.69,.16){$\alpha$} \put(.78,.155){$s$} \put(.87,.15){$\alpha$}
                \put(.69,-.02){$\beta$} \put(.78,.00){$s$} \put(.87,.00){$\beta$}
            \end{picture}
       \end{center}
       \caption{Blocks $\TWOFOURS_{4,4}$, $\JUSTFOURS_8$ and $\TWOTWOS_{2,2}$ in suspensions. The bottom two blocks have their leftmost (or rightmost) singularities identified by assumption.}\label{fig.gen_blocks}
    \end{figure}

    \begin{proof}
        We note that the proof follows in a similar manner to the proof of Lemma \ref{lem.blocks}. See Figure \ref{fig.gen_blocks} to verify the counting argument for $\POLES_1$. The result for $\POLES_i$ follows by induction.
    \end{proof}

    \begin{figure}[h]
        \begin{center}
           \setlength{\unitlength}{250pt}
            \begin{picture}(1,.50)
                \put(0,0){\includegraphics[width=\unitlength]{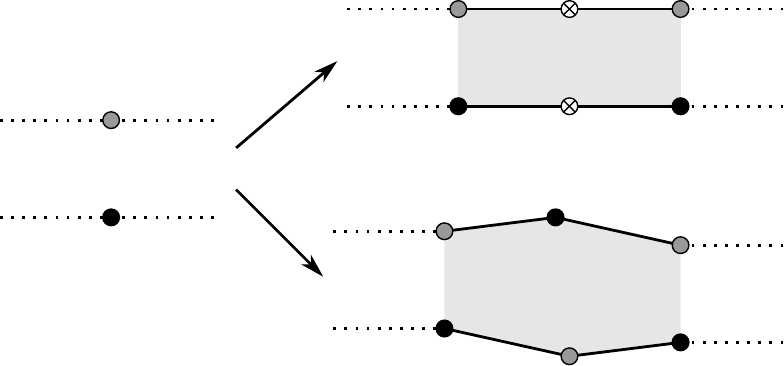}}

                \put(.70,.38){$\POLES_1$} \put(.69,.09){$\EVEN_2$}
                \put(.65,.47){$\alpha$} \put(.78,.47){$\alpha$} \put(.65,.28){$\beta$} \put(.78,.28){$\beta$}
                \put(.60,.20){$\alpha_1$} \put(.78,.19){$\beta_1$} \put(.61,-.02){$\beta_1$} \put(.78,-.02){$\alpha_1$}
            \end{picture}
       \end{center}
       \caption{Inserting blocks $\POLES_1$ and $\EVEN_2$ at paired singularities.}\label{fig.gen_insert}
    \end{figure}

    \begin{rem}
        We will use $\EVEN$ from Definition \ref{def.blocks} adding the convention that $\EVEN_0$ is the empty block.
    \end{rem}

    \begin{lem}\label{lem.q.EVENS}
        Let $\pi\in\genirr_\AAA$ be self-inverse and $i\geq 0$. When $\EVEN_{2i}$ is inserted at an even ordered (self-paired) singularity of order $2\ell\geq 0$, the singularity is now of order $2\ell+4i$. When $\EVEN_{2i}$ is inserted at fixed paired singularities of odd order $j\geq -1$, the paired singularities are now of order $j+2i$. This new permutation is irreducible and self-inverse.
    \end{lem}

    \begin{proof}
        This again follows from counting, induction and Figure \ref{fig.gen_insert}.
    \end{proof}

    \begin{defn}\label{def.FOURS}
        $$\begin{array}{rcl}
            \TWOFOURS_{4,4}&:=&\cmtrx{\LL{\alpha}{\eps}~\LL{\beta}{\beta}~\LL{\gamma}{\delta}~\LL{\alpha}{\eps}~\LL{\delta}{\gamma}}\\
            \JUSTFOURS_8&:=&\cmtrx{\LL{\alpha}{\beta}~\LL{\gamma}{\delta}~\LL{\alpha}{\beta}~\LL{\gamma}{\delta}}\\
            \JUSTFOURS_{4i+8}&:=&\EVEN_{2i}\JUSTFOURS_{8}\mbox{ for }i>0\\
            \TWOTWOS_{4i+2,4j+2}&:=&\cmtrx{\LL{\alpha}{\beta}&\SPACE\EVEN_{2i}&\LL{\alpha}{\beta}&\EVEN_{2j}}
        \end{array}$$
    \end{defn}

    \begin{lem}\label{lem.FOURS}
        $\TWOFOURS_{4,4}$ contributes two singularities of order four. If its leftmost (or rightmost) singularities are identified, $\JUSTFOURS_{4k}$ contributes a singularity of order $4k$. If its leftmost (or rightmost) singularities are identified, $\TWOTWOS_{4i+2,4j+2}$ contributes two singularities: one of order $4i+2$ and one of order $4j+2$.
    \end{lem}

    \begin{proof}
        We refer to Figure \ref{fig.gen_blocks} for $\TWOFOURS_{4,4}$, $\JUSTFOURS_8$ and $\TWOTWOS_{2,2}$. The other $\JUSTFOURS$ and $\TWOTWOS$ blocks follow inductively by Lemma \ref{lem.q.EVENS}.
    \end{proof}

    \begin{defn}\label{def.INSERT}
        For $i,j,k\geq 0$,
        $$\INSERT_{2i+1,(2j-1,2k-1)}:=\cmtrx{\LL{a}{c}&\EVEN_{2j}&\LL{a}{c}&\LL{b}{d}&\EVEN_{2k}&\LL{b}{d}&\EVEN_{2i}}$$
    \end{defn}

    \begin{lem}\label{lem.INSERT}
        Let $\pi\in\genirr_\AAA$ be self-inverse. When $\INSERT_{r,(p,q)}$, $p,q\geq -1$ and $r\geq 1$ all odd, is inserted at paired singularities of odd order $j\geq -1$, the paired singularities are now of order $r+j+1$ and two new pairs, of orders $p$ and $q$, are added. Specifically, if $j=-1$ (the fixed paired singularities are poles), the initially paired singularities are now of order $r$. The resulting permutation is self-inverse and irreducible.
    \end{lem}

    \begin{proof}
        This is an immediate result of Lemmas \ref{lem.POLES} and \ref{lem.q.EVENS} as an insertion of a $\INSERT$ block may be achieved by a sequence of $\POLES$ and $\EVEN$ insertions.
    \end{proof}

\sect{Strata with $g\in\{0,1\}$}\label{sec.gen.g01}

    \begin{lem}\label{lem.q.g0g1}
        Let $\RClass$ be a generalized Rauzy Class of genus $g\in\{0,1\}$ such that Conditions \eqref{eq.q.nec} hold. Then $\RClass$ contains a self-inverse element.
    \end{lem}

    \begin{proof}
        Let us begin with $g=0$. By Remark \ref{rem.q.signatureform}, the signature of $\RClass$ is of the form
            $$\sig(\RClass)=(2\ell_1,\dots,2\ell_m,j_1^2,j_2^2,\dots,j_n^2, -1^{2\ell_1+\cdots+2\ell_m+2j_1+\cdots+2j_n+4})$$
        where each $j_i$ is odd. In this case, only the signature must be verified. Start with the base permutation
            $$ \pi_0=\cmtrx{\LL{a}{c} ~\LL{a}{c} ~\POLES_{\ell_1}\SPACE\POLES_{\ell_2}\SPACE\cdots\SPACE\POLES_{\ell_m} ~\LL{b}{d} ~\LL{b}{d}}.$$
        By Lemma \ref{lem.POLES}, this permutation has signature $(2\ell_1,\dots,2\ell_m, -1^{2\ell_1+\cdots+2\ell_m+4})$. We then proceed inductively to derive $\pi_{i+1}$ by inserting $\POLES_{j_{i+1}+1}$ at any paired set of pole (order $-1$) singularities in $\pi_i$. Again by Lemma \ref{lem.POLES}, $\pi_i$ has signature
            $$\sig(\pi_i)=(2\ell_1,\dots,2\ell_m,j_1^2,j_2^2,\dots,j_i^2, -1^{2\ell_1+\cdots+2\ell_m+2j_1+\cdots+2j_i+4}).$$
        $\pi=\pi_n$ is self-inverse and has the desired signature, so $\pi\in\RClass$.\\

        As in the $g=0$ case, all that must be verified for $g=1$ is the signature. In this instance, the signature for $\RClass$ is
            $$ \sig(\RClass)=(2\ell_1,\dots,2\ell_m,j_1^2,\dots,j_n^2, -1^{2\ell_1+\cdots+2\ell_m+2j_1+\cdots+2j_n})$$
        for positive odd $j_i$'s. First assume $\ell_i>0$ for some $i$. Let
            $$ \pi_0=\cmtrx{\AB ~ \POLES_{\ell_m}\SPACE\cdots\SPACE\POLES_{\ell_1}~\BA}.$$
        Again, as in the $g=0$ case, we may perform insertions $\POLES_{j_i+1}$ on paired poles for each $i\in\{1,\dots,n\}$ to construct self-inverse $\pi$ with the appropriate signature. If $\ell_1=\dots=\ell_m=0$, then by the classification theorems in Section \ref{sec.gen.classify}, $n>0$. We begin with
            $$ \pi_1 = \cmtrx{\LL{a}{c} &\EVEN_{2}\POLES_{2j_1-2} &\LL{a}{c}&\underbrace{\SPACE\cdots\SPACE}_{m-1} &\LL{b}{d} &\LL{b}{d}}.$$
        By Lemmas \ref{lem.POLES} and \ref{lem.q.EVENS}, $\pi_1$ has signature $(0^m,j_1^2,-1^{2j_1})$. As $\pi_1$ contains poles, we may proceed as in the above cases to create self-inverse $\pi$ with the desired signature.
    \end{proof}

    \begin{exam}
        Consider the Rauzy Class $\RClass$ with signature $\sig(\RClass)=(2,2,1,1,1,1,-1^{12})$. The genus of $\RClass$ is $0$. We give each step as presented in the above proof:
        $$ \begin{array}{rcl}
            \pi_0&=&\cmtrx{\LL{1}{10}~\LL{1}{10}~\LL{6}{15}~\LL{6}{15}~\LL{7}{7}~\LL{8}{16}~\LL{8}{16}~\LL{9}{17}~\LL{9}{17}}\\
            \pi_1&=&\cmtrx{\LL{1}{10}~\LL{\mathbf{2}}{\mathbf{11}}~\LL{\mathbf{2}}{\mathbf{11}}~\LL{\mathbf{5}}{\mathbf{14}}~\LL{\mathbf{5}}{\mathbf{14}}~\LL{1}{10}~\LL{6}{15}~\LL{6}{15}~\LL{7}{7}~\LL{8}{16}~\LL{8}{16}~\LL{9}{17}~\LL{9}{17}}\\
            \pi=\pi_2&=&\cmtrx{\LL{1}{10}~\LL{2}{11}~\LL{\mathbf{3}}{\mathbf{12}}~\LL{\mathbf{3}}{\mathbf{12}}~\LL{\mathbf{4}}{\mathbf{13}}~\LL{\mathbf{4}}{\mathbf{13}}~\LL{2}{11}~\LL{5}{14}~\LL{5}{14}~\LL{1}{10}~\LL{6}{15}~\LL{6}{15}~\LL{7}{7}~\LL{8}{16}~\LL{8}{16}~\LL{9}{17}~\LL{9}{17}}.
        \end{array}$$
        The blocks inserted in each step are in bold. Note that $\sig(\pi_0)=(2,2,-1^8)$, $\sig(\pi_1)=(2,2,1,1,-1^{10})$ and $\sig(\pi_2)=(2,2,1,1,1,1,-1^{12})=\sig(\RClass)$. Therefore $\pi=\pi_2\in\RClass$.
    \end{exam}

\sect{Strata with a hyperelliptic component, $g\geq 2$}\label{sec.gen.g2hyp}

    \begin{rem}\label{rem.q.removable}
        The generalized permutations we construct in genus $g\geq 2$ come in one of two forms,
        $$ \pi=\cmtrx{\AB~\BLAH_0~\BA}\mbox{ or }\pi=\cmtrx{\LL{a}{c}~\BLAH_1~\LL{a}{c}~\BLAH_2~\LL{b}{d}~\BLAH_3~\LL{b}{d}}.$$
        We can add removable singularities to these permutations $\pi$ in a manner similar to Theorem \ref{thm.zero}. Suppose $k$ removable singularities are added to our permutation $\pi$ and call this new permutation $\pi'$. Then
        $$ \pi'=\cmtrx{\AB~\BLAH_0\SPACE^k~\BA}\mbox{ or }\pi'=\cmtrx{\LL{a}{c}~\BLAH_1~\LL{a}{c}~\SPACE^k\BLAH_2~\LL{b}{d}~\BLAH_3~\LL{b}{d}}$$
        if the marked singularity of $\pi'$ is removable, and
        $$ \pi'=\cmtrx{\AB~\SPACE^k\BLAH_0~\BA}\mbox{ or }\pi'=\cmtrx{\LL{a}{c}~\BLAH_1~\LL{a}{c}~\BLAH_2\SPACE^k~\LL{b}{d}~\BLAH_3~\LL{b}{d}}$$
        otherwise. The lemmas in this and the following sections all may use the above changes to allow for removable singularities.
    \end{rem}

    \begin{lem}\label{lem.q.hyp}
        The following results hold for strata that contain hyperelliptic components:
        \begin{itemize}
            \item A Rauzy class in the stratum $\QQQ(2j-1,2j-1,4k+2)$ contains a self-inverse element if and only if the marked singularity is of even order.
            \item A Rauzy class in the stratum $\QQQ(2j-1,2j-1,2k-1,2k-1)$ contains a self-inverse element if and only if the marked singularity is of even order.
             \item Any Rauzy class in the stratum $\QQQ(4j+2,4k+2)$ contains a self-inverse element.
       \end{itemize}
    \end{lem}

    \begin{proof}
        We refer to Section \ref{sec.hyperelliptic} to determine the hyperellipticity of a permutation. For a permutation $\pi\in\genirr_\AAA$ to be hyperelliptic, there must be a map that is piecewise of the form $z\mapsto -z+c$ on the interior of the polygon for any suspension for $\pi$. This involution must map non-removable singularities to non-removable singularities and fix $2g+2$ points, where $g$ is the genus of $\pi$.\\

        Suppose that $\sig(\RClass)=(0,2j-1,2j-1,2k-1,2k-1)$ and $\RClass$ is hyperelliptic. Let
        $$ \pi=\cmtrx{\LL{\alpha}{\gamma}~\LL{a_1}{a_{2j}}~\LL{a_2}{a_{2j-1}}~\LL{\cdots}{\cdots}~\LL{a_{2j-1}}{a_2} ~\LL{a_{2j}}{a_1}~\LL{\alpha}{\gamma}~\LL{\beta}{\delta}~\LL{c_1}{c_{2k}}~\LL{c_2}{c_{2k-1}}~\LL{\cdots}{\cdots} ~\LL{c_{2k-1}}{c_2} ~\LL{c_{2k}}{c_1}~\LL{\beta}{\delta}}.$$
        \begin{figure}[h]
            \begin{center}
               \setlength{\unitlength}{300pt}
                \begin{picture}(1,.48)
                    \put(0,0){\includegraphics[width=\unitlength]{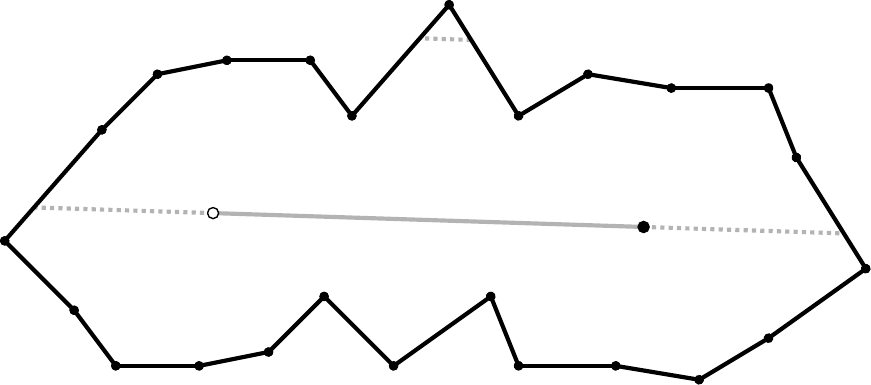}}

                    \put(0.03,0.26){$\alpha$} \put(0.44,0.41){$\alpha$} \put(0.55,0.41){$\beta$} \put(0.97,0.22){$\beta$}
                    \put(0.02,0.08){$\gamma$} \put(0.39,0.01){$\gamma$} \put(0.49,0.00){$\delta$} \put(0.95,0.06){$\delta$}
                    \put(.11,.34){$a_1$} \put(.20,.38){$\cdots$} \put(.28,.39){$\cdots$} \put(.38,.36){$a_{2j}$}
                    \put(.06,.03){$a_{2j}$} \put(.16,-.01){$\cdots$} \put(.25,.00){$\cdots$} \put(.35,.05){$a_1$}
                    \put(.60,.35){$c_1$} \put(.70,.37){$\cdots$} \put(.81,.36){$\cdots$} \put(.91,.31){$c_{2k}$}
                    \put(.53,.05){$c_{2k}$} \put(.63,-.02){$\cdots$} \put(.73,-.03){$\cdots$} \put(.85,.00){$c_1$}
                \end{picture}
           \end{center}
           \caption{There are two interior points that must be fixed by a hyperelliptic involution. If the involution fixes one, the other is fixed as well.}\label{fig.gen.hyp}
        \end{figure}
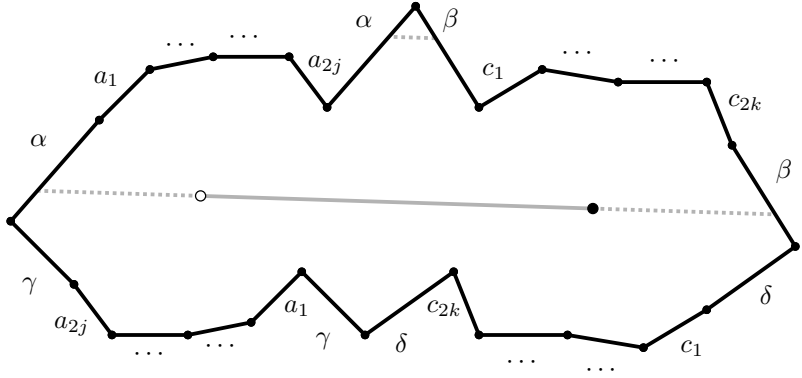
        From counting, it follows that $\sig(\pi)=\sig(\RClass)$. It remains to show that $\pi$ is hyperelliptic. Refer to Figure \ref{fig.gen.hyp}. Two points in the interior of the polygon, labeled as the white dot and black dot, must be fixed by the hyperelliptic involution, in order to fix the midpoints of each edge $a_i$ and $c_i$. We see that the grey geodesic is mapped to the dotted grey geodesic. Therefore an involution that fixes, say, the black dot fixes the white one as well. Because these two points are fixed, we see that there are $2+2j+2k = 2g+2$ fixed points. Therefore $\pi$ is hyperelliptic.\\

        Now assume that $\sig(\RClass)=(0,2j-1,2j-1,2k-1,2k-1)$ and $\RClass$ is non-hyperelliptic. We assume that either $j>1$ or $k>1$, as the stratum $\QQQ(1,1,1,1)$ is connected. Then
        $$ \pi=\cmtrx{\LL{\alpha}{\gamma}&\EVEN_{2j}&\LL{\alpha}{\gamma}&\LL{\beta}{\delta}&\EVEN_{2k}&\LL{\beta}{\delta}}$$
        belongs to $\RClass$ and is self-inverse.\\

        Consider hyperelliptic $\RClass$ with $\sig(\RClass)=(4k+2,2j-1,2j-1)$. Then a self-inverse element of $\RClass$ is
        $$ \pi=\cmtrx{\AB~\LL{c_1}{c_{2k}}~\LLcdots~\LL{c_k}{c_{k+1}}~\LL{\alpha}{\beta} ~\LL{a_1}{a_{2j}}~\LLcdots~\LL{a_{2j}}{a_1} ~\LL{\alpha}{\beta}~\LL{c_{k+1}}{c_k}~\LLcdots ~\LL{c_{2k}}{c_1}~\BA}.$$
        A similar argument for $\sig=(0,2j-1,2j-1,2k-1,2k-1)$ shows that $\pi\in\RClass$.\\

        Now consider non-hyperelliptic $\RClass$ such that $\sig(\RClass)=(4k+2,2j-1,2j-1)$. The stratum $\QQQ(2,-1,-1)$ is connected, so assume that $(j,k)\neq(1,0)$. Then a self-inverse element $\pi\in\RClass$ is
        $$ \pi=\cmtrx{\AB ~\EVEN_{2k}~\LL{\alpha}{\beta} ~\EVEN_{2j}~\LL{\alpha}{\beta}~\BA}.$$

        Finally, let $\RClass$ have signature $\sig(\RClass)=(4k+2,4j+2)$. Then there is self-inverse $\pi\in\RClass$ where
        $$ \pi=\cmtrx{\AB~\LL{c_1}{c_{2k}}~\LLcdots~\LL{c_k}{c_{k+1}}~\LL{\alpha}{\beta} ~\LL{a_1}{a_{2j+1}}~\LLcdots~\LL{a_{2j+1}}{a_1} ~\LL{\alpha}{\beta}~\LL{c_{k+1}}{c_k}~\LLcdots ~\LL{c_{2k}}{c_1}~\BA}$$
        if $\RClass$ is in the hyperelliptic component and
        $$ \pi=\cmtrx{\AB ~\LL{\alpha}{\beta} ~\SPACE~\EVEN_{2j} ~\LL{\alpha}{\beta}~\EVEN_{2k}~\BA} = \cmtrx{\AB~\TWOTWOS_{4j+2,4k+2}~\BA}$$
        if $\RClass$ is in the nonhyperelliptic component. These two constructions coincide if $j=k=0$. However, the stratum in this case, $\QQQ(2,2)$, is connected.\\
    \end{proof}

\sect{Other strata with $g\geq 2$}\label{sec.gen.g2}

    \begin{lem}\label{lem.q.exc}
        The four exceptional strata mentioned in Section \ref{sec.gen.classify}, $\QQQ(9,-1)$, $\QQQ(6,3,-1)$, $\QQQ(3,3,3,-1)$ and $\QQQ(12)$, each have two connected components. Any Rauzy Class $\RClass$ in the stratum $\QQQ(12)$ contains a self-inverse element. The other three strata do not contain self-inverses.
    \end{lem}

    \begin{proof}
        Any class contained in one of the three strata $\QQQ(9,-1)$, $\QQQ(6,3,-1)$ or $\QQQ(3,3,3,-1)$ cannot have self-inverses, as it would violate the necessary condition \eqref{eq.q.nec2}. So the only likely candidate is $\QQQ(12)$. Consider the two connected components labeled $\QQQ^{reg}(12)$ and $\QQQ^{irr}(12)$. By calculation using Zorich's Mathematica software \cite{c.Z_Web},
        $$ \pi=\cmtrx{1&2&1&2&3&4&3&4 \\ 5&6&5&6&7&8&7&8}$$
        is a self-inverse element of $\QQQ^{reg}(12)$. Again, by direct calculation
        $$ \pi=\cmtrx{1&2&3&4&5&2&6&7 \\ 7&8&4&3&6&8&5&1}$$
        is self-inverse and belongs to $\QQQ^{irr}(12)$.
    \end{proof}

    \begin{lem}\label{lem.q.alleven}
        Let $\RClass$ belong to $\QQQ(2\ell_1,\dots,2\ell_m)$ with $\ell_i>0$ and $g(\RClass)\geq 2$ that satisfy Conditions \eqref{eq.q.nec}. Then $\RClass$ contains a self-inverse element.
    \end{lem}

    \begin{proof}
        Suppose there are $m>0$ singularities of order $4$ for $\RClass$. Then the self-inverse permutation is
            $$ \pi=\cmtrx{\AB&\mA_{m-1}\SPACE\cdots\mA_{2}\SPACE\mA_1&\BA}\mbox{, where } \mA_i=\RHScase{\TWOFOURS_{4,4}, & \mbox{if }i=1,\\ \EVEN_2, & \mathrm{otherwise}.}$$
        Now suppose the signature for $\RClass$ is $\sig(\RClass)=(4\ell_1,\dots,4\ell_m)$, $\ell_i\geq 1$, and there is some $j$ such that $\ell_j>1$. Then the self-inverse element is
            $$ \pi=\cmtrx{\AB~\mA_m\SPACE\cdots\SPACE\mA_1~\BA}\mbox{, where } \mA_i=\RHScase{\JUSTFOURS_{4l_j}, & \mbox{if }i=j, \\ \EVEN_{2l_i}, & \mathrm{otherwise}.}$$
        Now suppose the signature of $\RClass$ is $(2\ell_1,\cdots,2\ell_m)$ where $\ell_j$ is odd for some $j$. If $\sig(\RClass)=(4k+2,4j+2)$, we may construct a self-inverse by Lemma \ref{lem.q.hyp}. Because the sum of singularities is a multiple of four, the number of odd $\ell_i$'s is even, denoted $2k>0$. Let $n=m-2k$ be the number of even $\ell_i$'s. Suppose $\ell_1$ is even. Then we may write the signature as $(4\ell_1',\dots,4\ell_n',2\ell_{n+1},\dots,2\ell_{m})$, where $2\ell_i'=\ell_i$ for $i\in\dset{n}$. The constructed example is
        $$ \pi=\cmtrx{\AB ~ \TWOTWOS_{2\ell_m,2\ell_{m-1}} \SPACE \cdots \SPACE \TWOTWOS_{2\ell_{n+2},2\ell_{n+1}} \SPACE\EVEN_{2\ell_n'} \SPACE \cdots \SPACE \EVEN_{2\ell_1'} ~ \BA}.$$
        If $\ell_1$ is odd, then we write the signature as $(2\ell_1,\dots,2\ell_{2k},4\ell_1',\dots,4\ell_n')$, where $2\ell_i'=\ell_{i+2k}$ for $i\in\ddset{2k+1}{m}$. We then define
        $$ \pi=\cmtrx{\AB ~\EVEN_{2\ell_n'} \SPACE \cdots \SPACE \EVEN_{2\ell_1'} \SPACE \TWOTWOS_{2\ell_{2k},2\ell_{2k-1}} \SPACE \cdots \SPACE \TWOTWOS_{2\ell_2,2\ell_1} ~\BA}.$$
        By Lemmas \ref{lem.FOURS} and \ref{lem.q.EVENS}, the examples constructed have the desired signature.
    \end{proof}

    \begin{lem}\label{lem.q.nopoles}
        Let $\RClass$ belong to $\QQQ(\ell_1,\dots,\ell_m)$, $\ell_i \geq 0$, with $g(\RClass)\geq 2$, such that $\RClass$ satisfies Conditions \eqref{eq.q.nec}. Then $\RClass$ contains a self-inverse element.
    \end{lem}

    \begin{proof}
        By Remark \ref{rem.q.signatureform}, let the signature be of the form
            $$ \sig(\RClass)=(2\ell_1,\dots,2\ell_m,j_1^2,\dots,j_n^2)$$
        for $\ell_i\geq 0$ and odd $j_i>0$.\\

        First assume that $n$ is even. Then $2\ell_1+\dots+2\ell_m$ is a multiple of $4$, so we may construct a self-inverse $\pi_0$ with signature $(2\ell_1,\dots,2\ell_m)$ by Lemma \ref{lem.q.alleven}. Define block $\BLAH_0$ by
            $$ \pi_0=\cmtrx{\AB& \BLAH_0 &\BA}.$$
        Note we may instead use $\BLAH_0=\EVEN_2$ and Remark \ref{rem.q.removable} in the case $(2\ell_1,\dots,2\ell_m)\in\{(0^{m-1},4),(4,0^{m-1})\}$. If $n=0$, then this is the desired element. So assume that $n=2k\geq2$. Let $\BLAH_1$ be the block obtained by reflecting $\BLAH_0$ horizontally, i.e.
            $$ \BLAH_0 = \cmtrx{\LL{\alpha_1}{\beta_1}&\LL{\alpha_2}{\beta_2}&\LLcdots&\LL{\alpha_p}{\beta_p}}\Rightarrow
                \BLAH_1 = \cmtrx{\LL{\alpha_p}{\beta_p}&\LL{\alpha_{p-1}}{\beta_{p-1}}&\LLcdots&\LL{\alpha_1}{\beta_1}}.$$
        Note that now the \textbf{leftmost} singularity of $\BLAH_1$ is order $2\ell_1$. Then define
            $$ \pi_1 = \cmtrx{ \LL{\alpha}{\gamma} & \LL{\alpha}{\gamma} & \BLAH_1 & \LL{\beta}{\delta} & \LL{\beta}{\delta}}.$$
        So $\pi_1$ has signature $(2\ell_1,\dots,2\ell_m,-1^4)$. We now move from $\pi_i$ to $\pi_{i+1}$ by inserting a block of form $\INSERT_{j_{2i-1},(j_{2i},-1)}$ to $\pi_i$ at any choice of paired poles. So $\pi_i$ will have signature $(2\ell_1,\dots,2\ell_m,j_1^2,j_2^2,\dots,j_{2i-3}^2,j_{2i-2}^2,-1^4)$ for all $i$ such that $1 < i \leq k$, by Lemma \ref{lem.INSERT}. To derive $\pi_{k+1}$, we replace the four paired poles in $\pi_k$ with $\EVEN_{j_{2k-1}+1}$ for one pair and $\EVEN_{j_{2k}+1}$ for the other. By Lemma \ref{lem.q.EVENS}, the signature of $\pi:=\pi_{k+1}$ is as desired.\\

        Now assume that $n=2k+1$ is odd. Because the sum of the orders of singularities must be a multiple of $4$, there exists $p$ such that $\ell_p$ is odd. Fix this $\ell_p$. We may again from Lemma \ref{lem.q.alleven} construct self-inverse $\pi_0$ with signature $(2\ell_1,\dots,2\ell_{p-1},2\ell_{p+1},\dots,2\ell_m)$. Define $\BLAH_0$ and $\BLAH_1$ as above. If $p=1$, let
            $$ \pi_1:=\cmtrx{\AB&\BLAH_1\SPACE\EVEN_{2\ell_p-2}&\LL{\alpha}{\beta}&\LL{\alpha}{\beta}&\BA},$$
        and if $p\neq 1$, let
            $$ \pi_1:=\cmtrx{\AB&\LL{\alpha}{\beta}&\LL{\alpha}{\beta}&\EVEN_{2\ell_p-2}\SPACE\BLAH_1&\BA}.$$
        By Lemma \ref{lem.q.EVENS}, $\pi_1$ has signature $(2\ell_1,\dots,2\ell_m,-1^2)$. As before, insert block $\INSERT_{j_{2i-1},(j_{2i},-1)}$ into the paired poles of $\pi_i$ to get $\pi_{i+1}$. By Lemma \ref{lem.INSERT}, the signature of $\pi_i$ is $(2\ell_1,\dots,2\ell_m,j_1^2,j_2^2,\dots,j_{2i-3}^2,j_{2i-2}^2,-1^2)$ for all $i$ such that $1<i\leq k+1$. Let $\pi:=\pi_{k+2}$ be derived by inserting a block $\EVEN_{j_n+1}$ at the paired poles of $\pi_{k+1}$. By Lemma \ref{lem.q.EVENS}, $\pi$ has the desired signature.
    \end{proof}

    \begin{corollary}\label{cor.q.twopoles}
        A Rauzy Class $\RClass$ with signature $(2\ell_1,\dots,2\ell_m,j_1^2,\dots,j_n^2,-1^2)$, $j_i>0$ odd, contains a self-inverse element.
    \end{corollary}

    \begin{proof}
        We may perform the construction in the proof of Lemma \ref{lem.q.nopoles} but consider the final two odd ordered singularities to be the poles, i.e. follow the pervious proof with $\sig(\RClass) = (2\ell_1,\dots,2\ell_m,j_1^2,\dots,j_n^2,j_{n+1}^2)$, where $j_{n+1}=-1$. The final step of that proof would insert a block $\EVEN_{j_{n+1}+1}=\EVEN_0$, which leaves the two poles.
    \end{proof}

    \begin{lem}\label{lem.q.final}
        Any generalized Rauzy Class $\RClass$ contains a self-inverse element if $\RClass$ satisfies Conditions \eqref{eq.q.nec}.
    \end{lem}

    \begin{proof}
        Because of the assumption on $\RClass$ and by Remark \ref{rem.q.signatureform}, we write the signature as
            $$ \sig_0:=(2\ell_1,\dots,2\ell_m,j_1^2,\dots,j_n^2,-1^p).$$
        where $j_i>0$ is odd for all $i\in\dset{n}$. Notice that $p$ is even. If $p=0$ or $p=2$, we can construct a self-inverse $\pi\in\RClass$ by Lemma \ref{lem.q.nopoles} or Corollary \ref{cor.q.twopoles}. We will describe a procedure to ``reduce'' signature $\sigma_0$ to signature $\sigma_n$ with either $0$ or $2$ poles. We then construct self-inverse $\pi_n$ with signature $\sigma_n$ and perform appropriate insertions on $\pi_n$ to derive self-inverse $\pi_0$ with the desired signature (see Figure \ref{fig.reductions}). Note that each reduction maintains the genus of the original signature.\\

        \begin{figure}[h]
            $$ \xymatrix{\sigma=\sigma_0 \ar[rr]^{\mathrm{reduce}}& &\sigma_1 \ar[rr]^{\mathrm{reduce}}& &\dots\ar[rr]^{\mathrm{reduce}}& & \sigma_n \ar[d]^{\mathrm{construction}}\\
                \pi=\pi_0 & &\pi_1 \ar[ll]^{\mathrm{insert}} & &\dots \ar[ll]^{\mathrm{insert}} & & \pi_n \ar[ll]^{\mathrm{insert}}}$$
            \caption{Performing successive reductions on $\sigma_0$ to more desirable signature $\sigma_n$. Given $\pi_n$, repeated insertions give desired $\pi_0$.}\label{fig.reductions}
        \end{figure}
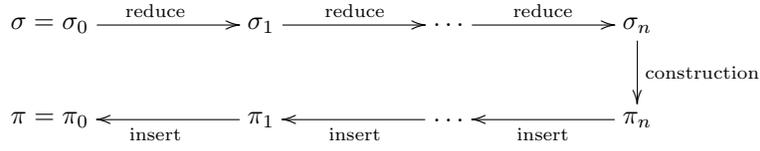

        Suppose first that $p>2$ and there exists some positive $\ell_i$ in signature $\sig_t=(2\ell_1,\dots,2\ell_m,j_1^2,\dots,j_n^2,-1^p)$. Then consider a signature of the form
            $$\sig_{t+1}:=(2\ell_1,\dots,2\ell_{i-1},2\ell_i-2,2\ell_{i+1},\dots,2\ell_m,j_1^2,\dots,j_n^2,-1^{p-2}).$$
        If we assume that we have constructed self-inverse $\pi_{t+1}$ with signature $\sig_{t+1}$, we construct $\pi_t$ by inserting a $\POLES_1$ block at the singularity of order $2\ell_i-2$. By Lemma \ref{lem.INSERT}, $\pi_t$ will have signature $\sig_t$. We may perform these reductions until either $p=2$ (ending the proof), or $\ell_i=0$ for all $i$.\\

        Now consider signature $\sig_t=(0^m,j_1^2,\dots,j_n^2,-1^p)$ with $p>2$ and $j_i>1$ for some $i$. We consider signature
            $$ \sig_{t+1}:=(0^m,j_1^2,\dots,j_{i-1}^2,(j_i-2)^2,j_{i+1},\dots,j_n^2,-1^{p-4}).$$
        If we construct self-inverse $\pi_{t+1}$ with signature $\sig_{t+1}$, we get $\pi_t$ by inserting a $\POLES_2$ block at the paired singularities of order $j_i-2$. Again by Lemma \ref{lem.INSERT}, $\pi_t$ has signature $\sig_t$. Proceed with these reductions until either $p\leq 2$, ending the proof, or $j_i=1$ for all $i$.\\

        So assume that we are working with signature $\sig_t=(0^m,1^n,-1^p)$. Consider signature
            $$ \sig_{t+1}:=(0^m,1^{n-2},-1^{p-2}).$$
        Once we create self-inverse $\pi_{t+1}$ with signature $\sig_{t+1}$, let $\pi_t$ be the result of inserting a $\POLES_2$ block at any paired set of poles in $\pi_{t+1}$. Then $\pi_t$ has signature $\sig_t$ by Lemma \ref{lem.INSERT}. This process will terminate when $p=2$, when such a self-inverse may be constructed by Corollary \ref{cor.q.twopoles}.
    \end{proof}

    We will end with a few examples using the procedure in the proof of the previous lemmas.

    \begin{exam}
        Consider Rauzy Class $\RClass\subset\genirr_\AAA$ with $\sig(\RClass)=\sig_0=(6,1^2,-1^4)$. We perform a reduction of the first kind to get $\sig_1=(4,1^2,-1^2)$. By Corollary \ref{cor.q.twopoles}, we construct
        $$ \pi_1 = \cmtrx{\LL{1}{8}~\LL{2}{3}~\LL{3}{2}~\LL{1}{8}~\LL{4}{5}~\LL{5}{4}~\LL{7}{10}~\LL{7}{10}}$$
        with $\sig(\pi_1)=\sig_1$. We then perform the insertion of a $\POLES_1$ block to get
        $$\pi_0=\cmtrx{\LL{1}{8} ~\LL{2}{3}~\LL{3}{2} ~\LL{1}{8} ~\LL{4}{5}~\LL{5}{4} ~\LL{\mathbf{6}}{\mathbf{9}}~\LL{\mathbf{6}}{\mathbf{9}} ~\LL{7}{10}~\LL{7}{10}}.$$
        If $\pi=\pi_0$ then $\sig(\pi) = \sig_0$. Therefore $\pi\in\RClass$.
    \end{exam}

    \begin{exam}
        Let $\RClass$ have signature $\sig(\RClass)=\sig_0=(0,3^2,1^4,-1^6)$. A reduction of the second kind will yield $\sig_1=(0,1^6,-1^2)$. Corollary \ref{cor.q.twopoles} will give
        $$ \pi_1 = \cmtrx{\LL{1}{9} ~\LL{2}{a} ~\LL{3}{4} ~\LL{4}{3} ~\LL{2}{a} ~\LL{5}{b} ~\LL{6}{7} ~\LL{7}{6} ~\LL{5}{b} ~\LL{1}{9} ~\LL{8}{c} ~\LL{8}{c}}$$
        with $\sig(\pi_1)=(0,1^6,-1^2).$ An insertion of a $\POLES_2$ block to get
        $$ \pi_0=\cmtrx{\LL{1}{9} ~\LL{\mathbf{d}}{\mathbf{f}} ~\LL{\mathbf{d}}{\mathbf{f}} ~\LL{\mathbf{e}}{\mathbf{g}} ~\LL{\mathbf{e}}{\mathbf{g}} ~\LL{2}{a} ~\LL{3}{4} ~\LL{4}{3} ~\LL{2}{a} ~\LL{5}{b} ~\LL{6}{7} ~\LL{7}{6} ~\LL{5}{b} ~\LL{1}{9} ~\LL{8}{c} ~\LL{8}{c}}.$$
        We let $\pi=\pi_0$ and note that $\sig(\pi) = (0,3^2,1^4,-1^6)$. So $\pi\in\RClass$.
    \end{exam}

    \begin{exam}
        Finally consider $\RClass$ with $\sig(\RClass)=\sig_0=(0,1^8,-1^4)$. The third type of reduction gets $\sig_1=(0,1^6,-1^2)$ which, by Corollary \ref{cor.q.twopoles}, allows us to make associated
        $$ \pi_1=\cmtrx{\LL{1}{9} ~\LL{2}{a} ~\LL{3}{4} ~\LL{4}{3} ~\LL{2}{a} ~\LL{5}{b} ~\LL{6}{7} ~\LL{7}{6} ~\LL{5}{b} ~\LL{1}{9} ~\LL{8}{c} ~\LL{8}{c}}.$$
        We insert a $\POLES_2$ block to get
        $$ \pi_0 = \cmtrx{\LL{1}{9} ~\LL{2}{a} ~\LL{3}{4} ~\LL{4}{3} ~\LL{2}{a} ~\LL{5}{b} ~\LL{6}{7} ~\LL{7}{6} ~\LL{5}{b} ~\LL{1}{9} ~\LL{8}{c} ~\LL{\mathbf{d}}{\mathbf{f}} ~\LL{\mathbf{d}}{\mathbf{f}} ~\LL{\mathbf{e}}{\mathbf{g}} ~\LL{\mathbf{e}}{\mathbf{g}} ~\LL{8}{c}}.$$
        $\pi=\pi_0$ has signature $\sig(\pi)=(0,1^8,-1^4)$, and $\pi\in\RClass$.
    \end{exam}

\chap{Explicit Lagrangian Subspaces in True Rauzy Classes}\label{chapLag}

    Let us now consider $\pi\in\irr_d$ and unit suspension $S_1=S(\pi,\mone,\tau)$, where $\mone=(1,\dots,1)$, of genus $g$. We will consider a natural question, when do the closed vertical loops in $S_1$ span a $g$-dimensional subspace in homology $H=H_1(S_1)$? Consider the symplectic space $(H,\omega)$ where $\omega$ is the (algebraic) intersection number. In Theorem \ref{t.Isotropic}, we verify algebraically that the vertical loops do not intersect. We then prove in Theorem \ref{t.Lag.Perm} that if $\pi=\inv{\pi}$, then there are $g$ vertical loops independent in homology. Self-inverses have transpositions, i.e. letters that are interchanged, and fixed letters. Theorem \ref{t.SLag.Perm} shows that the transposition pairs in block-constructed $\pi$ from \Chapword\ \ref{chapTrue} form the basis of this $g$-dimensional space.

\sect{Symplectic Space}
    In this section, we give the definition of a symplectic space and list some basic properties of such spaces. Using the well known result of Proposition \ref{prop.iso_dim}, we prove Lemma \ref{l.splitting} which we will use to prove Theorem \ref{t.Lag.Perm}.

    \begin{defn}
        A vector space and bilinear form $(H,\omega)$ is \term{symplectic} if for all $v\in H$
        \begin{itemize}
            \item $\omega(v,v)=0$, or $(H,\omega)$ is \term{isotropic}
            \item$\omega(u,v)=0$ for all $u\in H$ implies $v=0$, or $(H,\omega)$ is \term{non-degenerate}
        \end{itemize}
    \end{defn}

    \begin{defn}
        Given a symplectic space $(H,\omega)$ and subspace $V$, define $V^\omega$ as
        $$ V^\omega := \{u\in H: \omega(u,v)=0\mbox{ for every }v\in V\}.$$
    \end{defn}

    \begin{prop}\label{prop.iso_dim}
        Let $(H,\omega)$ be a symplectic space with subspace $V$. Then
        $$\dim V + \dim V^\omega = \dim H.$$
    \end{prop}

    \begin{defn}
        Given symplectic space $(H,\omega)$, $V\subset H$ is \term{isotropic} if $V\subset V^\omega$. $V$ is \term{Lagrangian} if $V=V^\omega$.
    \end{defn}

    \begin{lem} \label{l.splitting}
        Let $W,V$ be subspaces of symplectic space $(H,\omega)$ such that:
        \begin{itemize}
            \item $V$ is isotropic,
            \item $H = V + W$ and
            \item For each $w\in W$, $w\neq 0$, there exists $v\in V$ such that $\omega(w,v)\neq 0$.
        \end{itemize}
        Then $V$ is Lagrangian.
    \end{lem}

    \begin{proof}
        Suppose $V\subsetneq V^\omega$. Choose $q\in V^\omega\setminus V$. So $q\notin V$ and $\omega(q,u)=0$ for all $u\in V$. However,
        $q=v+w$ for $v\in V$, $w\in W$, $w\neq 0$. However there exists $v'\in V$ such that $\omega(q,v') = \omega(w,v') \neq 0$, a contradiction to the choice of $q$.
    \end{proof}

\sect{Lagrangian Permutations}

    Let $\pi=(\pi_0,\pi_1)\in\irr_d$ over alphabet $\AAA$, recall the definition of closed loops $\gamma_\alpha$ and their corresponding cycles $c_\alpha = [\gamma_\alpha]$ in homology (see Section \ref{sec.spin}). We recall that $(H^\pi,\omega)$, where $H^\pi = \Omega_\pi\RR^\AAA$ and $\omega(\Omega_\pi u, \Omega_\pi v)= u^t\Omega_\pi v$, form a symplectic space with dimension $2g(\pi)$. Let $\Omega=\Omega_\pi$, and recall that

    \begin{equation}\label{eq.omega_equation}
        \Omega_{\alpha,\beta} = \chi_{\pi_1(\alpha)\leq \pi_1(\beta)}-\chi_{\pi_0(\alpha)\leq \pi_0(\beta)}.
    \end{equation}

    \begin{defn}
        For alphabet $\AAA$ and $\pi=(\pi_0,\pi_1)\in\irr_d$, consider the natural action of $\pi$ on $\AAA$, $\pi_\AAA$, by
        $$ \pi_\AAA:= \inv{\pi_0}\circ \pi \circ \pi_0 = \inv{\pi_0}\circ (\pi_1\circ\inv{\pi_0}) \circ \pi_0 = \inv{\pi_0}\circ \pi_1.$$
        Denote the set of orbits of $\AAA$ of $\pi_\AAA$ by
        $$ \AA := \{\BBB\subseteq \AAA: \BBB= \OOO_{\pi_\AAA}(\alpha)\mbox{ for }\alpha\in\AAA\}.$$
        For each $k\in\dset{d}$, let
        $$ \AA_k:= \{\BBB\in\AA: \#\BBB=k\}\subseteq \AA.$$
    \end{defn}

    \begin{exam}\label{ex.perm_AAA}
        Consider
        $$ \pi=\cmtrx{\LL{a}{f}~\LL{b}{e}~\LL{c}{b}~\LL{d}{d}~\LL{e}{c}~\LL{f}{a}}\in\irr_6.$$
        In this case $\pi_\AAA(a) = f$, $\pi_\AAA(b) = c$, $\pi_\AAA(c)=e$, $\pi_\AAA(d) = d$, $\pi_\AAA(e) = b$ and $\pi_\AAA(f) = a$. Also, $\AA = \{\{a,f\},\{b,c,e\},\{d\}\}$, $\AA_1 = \{\{d\}\}$, $\AA_2 = \{\{a,f\}\}$ and $\AA_3 = \{\{b,c,e\}\}$.
    \end{exam}

    \begin{defn}
        Let $\{\ee_\alpha\}_{\alpha\in\AAA}$ form the standard orthonormal basis of $\RR^\AAA$. For $\BBB\subseteq\AAA$, let
        $$ \ee_\BBB := \sum_{\alpha\in\BBB}\ee_\alpha.$$
        Let the vector space of vertical cycles under $\pi$ be
        $$ V^\pi := \spaN\{\ee_\BBB: \BBB\in\AA\}\subseteq \RR^\AAA$$
        and for $k\in\dset{d}$, let $V^\pi_k := \spaN\{\ee_\BBB: \BBB\in\AA_k\}$.
        Let $W^\pi$ be naturally defined by $\RR^\AAA = V^\pi\oplus W^\pi$.
    \end{defn}

    \begin{defn}
        Let the image of the vertical cycles under homology be $\vv_\BBB:=\Omega\ee_\BBB$, their span be
        $$ H^\pi_V := \Omega V^\pi = \spaN\{\vv_\BBB: \BBB\in\AA\}$$
        and $H^\pi_{V_k} :=\Omega V^\pi_k = \spaN\{\vv_\BBB: \BBB\in\AA_k\}$ for $k\in\dset{d}$. Also, let $H^\pi_W = \Omega W^\pi$.
    \end{defn}

    \begin{exam}
        Consider again $\pi$ from Example \ref{ex.perm_AAA} with $\AA=\{\{a,f\},\{b,c,e\},\{d\}\}$. We have that
        $$ V^\pi = \spaN\{\ee_{a,f}, \ee_{b,c,e}, \ee_{d}\}\mbox{ and } W^\pi = \spaN\{\ee_{a}-\ee_f, \ee_b-\ee_c, \ee_c-\ee_e\}.$$
        Consider the definition of $\Omega_\pi$ from Equation \ref{eq.omega_pi2}. Then
        $$ \begin{array}{l}
                H_V^\pi = \spaN\{\underbrace{(1,0,0,0,0,-1)^t}_{\vv_{a,f}}, \underbrace{(3,1,1,0,-2,-3)^t}_{\vv_{b,c,e}}, \underbrace{(1,0,1,0,-1,-1)^t}_{\vv_d}\}\mbox{ and }\\
                H_W^\pi = \spaN\{\underbrace{(-1,-2,-2,-2,-2,-1)^t}_{\ww_1}, \underbrace{(0,0,0,1,0,0)^t}_{\ww_2}, \underbrace{(0,-1,-1,-2,-1,0)^t}_{\ww_3}\}.
            \end{array}$$
        One may then calculate that $\omega(\vv_\BBB,\vv_\CCC)=\omega(\ww_i,\ww_j)=0$ for each $\BBB,\CCC\in\AA$ and $i,j\in\{1,2,3\}$. The other values are given by the following table
        $$ \omega(\vv_\BBB,\ww_i)=\begin{array}{|r|ccc|}
            \hline
                \BBB ~ \setminus ~ i & 1 & 2 & 3\\
            \hline
                \{a,f\} & -2&0&0\\
                \{b,c,e\}& -6&0&-3\\
                \{d\}& -2&1&-2\\
            \hline
            \end{array}$$
        and that $\omega$ is antisymmetric.
    \end{exam}

    \begin{rem}
        For any $\BBB,\CCC\subseteq\AAA$, we note the following formula
        \begin{equation}\label{eq.cycle_form}
            \omega(\vv_\BBB,\vv_\CCC) = \sum_{\alpha\in\BBB, \beta\in\CCC} \Omega_{\alpha,\beta}.
        \end{equation}
    \end{rem}

    \begin{thm}\label{t.Isotropic}
        For any $\pi\in\irr_d$, $H_V^\pi$ is isotropic.
    \end{thm}

    \begin{proof}
        It suffices to show that for every $\BBB,\CCC\in\AA$, $\omega(\vv_\BBB,\vv_\CCC)=0$. Because each $\BBB\in\AA$ is an orbit of $\pi_\AAA$, for $\eps\in\{0,1\}$,
        \begin{equation}\label{eq.orbit_good}
            \forall \alpha\in\BBB, ~\exists \beta\in\BBB \mbox{ s.t. }\pi_\eps(\alpha) = \pi_{1-\eps}(\beta).
        \end{equation}
        So for each $k\in\dset{d}$,
        \begin{equation}\label{eq.orbit_result}
            \#\{\alpha\in\BBB:\pi_0(\alpha)\leq k\} = \#\{\alpha\in\BBB:\pi_1(\alpha)\leq k\}.
        \end{equation}
        The calculation follows:
        $$ \begin{array}{rcll}
                \omega(\vv_\BBB,\vv_\CCC) & = & \underset{\alpha\in\BBB, \beta\in\CCC}\sum \Omega_{\alpha,\beta} & \mbox{by \eqref{eq.cycle_form}}\\
                    & = & \underset{\alpha\in\BBB}\sum\underset{\beta\in\CCC}\sum \chi_{\pi_1(\alpha) \leq \pi_1(\beta)} - \chi_{\pi_0(\alpha) \leq \pi_0(\beta)} & \mbox{by \eqref{eq.omega_equation}}\\
                    & = & \underset{\alpha\in\BBB}\sum \#\{\beta\in\CCC:\pi_1(\alpha) \leq \pi_1(\beta)\}-\#\{\alpha\in\BBB:\pi_0(\alpha) \leq \pi_0(\beta)\} & \\
                    & = & \underset{\alpha\in\BBB}\sum \#\{\beta\in\CCC:\pi_1(\alpha) \leq \pi_1(\beta)\}-\#\{\alpha\in\BBB:\pi_0(\alpha) \leq \pi_1(\beta)\} & \mbox{by \eqref{eq.orbit_good}}\\
                    & = & 0 &\mbox{by \eqref{eq.orbit_result}}
            \end{array}$$
        Therefore $H^\pi_V$ is isotropic.
    \end{proof}

    \begin{defn}
        $\pi\in\irr_d$ is \term{Lagrangian} if $H_V^\pi$ is Lagrangian.
    \end{defn}

    \begin{exam}
        Let $\pi = (4,1,3,2)$. In this case $H^\pi_V$ is spanned by two vectors, $(1,1,0,-2)^t$ and $(0,1,0,-1)^t$. So $\pi$ is Lagrangian. On the other hand, if $\pi' = (3,1,4,2)$ then $H_V^{\pi'}$ is spanned by only the vector $(1,2,-2,-1)^t$. Therefore $\pi'$ is not Lagrangian.
    \end{exam}

    \begin{rem}
        Naturally
        $$V^\pi = \bigoplus_{k=1}^d V^\pi_k.$$
        When $\pi$ is self-inverse, $V^\pi_k=\{0\}$ for all $k>2$. So
        $$ V^\pi = V^\pi_2 \oplus V^\pi_1,$$
        where $V^\pi_1$ corresponds to fixed $\alpha$ under $\pi_\AAA$, and $V^\pi_2$ corresponds to \term{transpositions}, pairs of letters $\{\alpha,\beta\}$ that are switched under $\pi_\AAA$.
        In this case, let
        $$W^\pi= \spaN\{\vv_{\alpha,-\beta}: \{\alpha,\beta\}\in\AA_2\}$$
        where $\ee_{\alpha,-\beta}=\ee_\alpha-\ee_\beta$ and $\vv_{\alpha,-\beta} = \Omega\ee_{\alpha,-\beta}$. It follows that $\RR^\AAA=V^\pi\oplus W^\pi$.
    \end{rem}

    \begin{thm}\label{t.Lag.Perm}
        Suppose $\pi\in\irr_d$ is self-inverse. Then $\pi$ is Lagrangian.
    \end{thm}

    \begin{proof}
        Because $\RR^\AAA=V^\pi\oplus W^\pi$, $H^\pi = H^\pi_V+H^\pi_W$. Also by Theorem \ref{t.Isotropic}, $H^\pi_V$ is isotropic. Consider any pair $\{\alpha,\beta\}\in\AA_2$, then
        $$ \omega(\vv_{\alpha,-\beta},\vv_{\alpha,\beta}) = \Omega_{\alpha,\beta}- \Omega_{\beta,\alpha}= 2\Omega_{\alpha,\beta} = \pm 2\neq 0.$$
        As the $\vv_{\alpha,-\beta}$'s form a basis for $H_W^\pi$ and $\vv_{\alpha,\beta}\in H_V^\pi$, we conclude that $H_V^\pi$ is Lagrangian by Lemma \ref{l.splitting}.
    \end{proof}

    \begin{corollary}
        If $\pi$ is self-inverse, then $\pi$ has at least $g(\pi)$ transpositions.
    \end{corollary}

    \begin{proof}
        By construction $\#\AA_2 = \dim V^\pi_2 = \dim W^\pi \geq \dim H_W^\pi = g(\pi)$.
    \end{proof}
    
    The following provides an alternative proof of Lemma 4.4 in \cite{c.For02}.
    
    \begin{corollary}\label{cor.Lag.For}
        In every connected component $\CCC$ of every stratum of Abelian differentials, let $\LLL$ be the set of $q\in\CCC$ such that:
        \begin{itemize}
            \item the vertical trajectories defined by $q$ that avoid singularities are periodic,
            \item the span of these vertical trajectories span a Lagrangian subspace in homology.
        \end{itemize}
        Then the set $\LLL$ is dense in $\CCC$.
    \end{corollary}
    
    \begin{proof}
        By Theorem \ref{thm.main}, every Rauzy Class $\RClass$ in $\CCC$ contains a self-inverse permutation $\pi$. Theorem \ref{t.Lag.Perm} shows that any unit suspension $S_1=S(\pi,\mone,\tau)$ satisfies the conditions of the claim. It is known that the Teichm\"uller geodesic flow (see Section \ref{sec.surface}) is ergodic, and therefore, we may choose $\tau$ such that the inverse flow is dense. It follows from an argument similar to Proposition 2.11 in \cite{c.Ve84_II}, for example, that every differential in the inverse flow also satisfies the conditions of the claim.
    \end{proof}

\sect{Transposition Lagrangian Permutations}

    Theorem \ref{t.Lag.Perm} shows that the vertical cycles of any self-inverse permutation span a Lagrangian subspace in homology. In general, choosing a basis from these cycles still requires calculation. However the block-constructed permutations in Definition \ref{def.blockConstructed} enjoy an additional property: the transpositions cycles form a basis for the Lagrangian subspace. We make this definition explicit, and then prove this result in Theorem \ref{t.SLag.Perm}.

    \begin{defn}
        A self-inverse permutation $\pi$ is \term{transposition Lagrangian} if $\dim H_{V_2}^\pi = \dim V^\pi_2 = g(\pi)$.
    \end{defn}

    \begin{exam}
        The permutation $\pi = (7,5,3,6,2,4,1)$ is self-inverse with $g(\pi)=3$. There are $3$ transposition pairs, $\{1,7\}$,$\{2,5\}$ and $\{4,6\}$, and one fixed letter $\{3\}$. However, we see that $\vv_{2,5}=\vv_{4,6} = (2,1,0,1,-1,-1,-2)^t$, $\vv_{1,7}=(1,0,0,0,0,0,-1)^t$ and $\vv_3=(1,1,0,0,-1,0,-1)^t$. So $\dim H_{V_2}^\pi = 2 < 3 = \dim H_V^\pi$, and therefore $\pi$ is Lagrangian, but not transposition Lagrangian.
    \end{exam}

    \begin{thm}
        If $\pi=(d,d-1,\dots,2,1)\in\irr_d$ then $\pi$ is transposition Lagrangian.
    \end{thm}

    \begin{proof}
        Suppose
        $$ \pi = \cmtrx{a_1 & a_2 & \dots & a_{d-1} & a_d \\ a_d & a_{d-1} & \dots & a_2 & a_1}.$$
        We recall that
        $$\Omega_{a_i,a_j} = \RHScase{1, & \mbox{if }i<j, \\ 0, & \mbox{if }i=j, \\ -1, & \mbox{if }i>j.}$$
        There are exactly $g=g(\pi)$ transpositions, $\{a_i,a_{d+1-i}\}\in\AA_2$ for $i\in\dset{g}$. Because $\#\AA_2 = \dim V^\pi_2 = g$, we must now show that the vectors $\vv_{a_i,a_{d+1-i}}$ are linearly independent. We see that for $j,k\in\dset{g}$,
        $$ \omega(\vv_{a_k}, \vv_{a_j,a_{d+1-j}}) = \Omega_{a_k,a_j} + \Omega_{a_k,a_{d+1-j}}= \RHScase{0, & \mbox{if }k>j,\\ 1, & \mbox{if }k=j, \\ 2, & \mbox{if }k<j.}$$
        So consider any $c_1,\dots,c_g\in\RR$ such that
        $$ \ww = c_1\vv_{a_1,a_{d}}+ \dots + c_g\vv_{a_g,a_{d+1-g}}= 0.$$
        It is clear now that $\omega(\vv_{a_g},\ww) = c_g = 0$. Inductively, if $c_g=\dots = c_k=0$, then $\omega(\vv_{a_{k-1}},\ww) = 2c_g+\dots+2c_k+c_{k-1}=c_{k-1}=0$. So $c_1=\dots=c_g=0$, implying the vectors $\vv_{a_i,a_{d+1-i}}$ are linearly independent. So $\dim H^\pi_{V_2} = \dim V^\pi_2=g$.
    \end{proof}

    \begin{thm} \label{t.SLag.Perm}
        Let $\pi$ be a block-constructed permutation. Then $\pi$ is transposition Lagrangian.
    \end{thm}

    \begin{proof}
        We begin by showing that $\#\AA_2 = \dim V^\pi_2 = g(\pi)$. Recall the formula
        $$ g=g(\pi) = \frac{1}{2}\sum_i \ell_i + 1$$
        where the $\ell_i$'s are the degrees of the singularities of $\pi$. Every block constructed self-inverse is of the form
        $$ \pi = \stdperm{\mB_1\SPACE\cdots\SPACE\mB_k}.$$
        So $\{\lA,\lB\}\in\AA_2$. For each $i\in\dset{k}$, $$\mB_i\in\{\EVEN_{2m}=\EVEN_2^m,\ODD_{2,2},\ODD_{2n}=\ODD_2^n,\PAIR_{2m+1,2n+1}=\EVEN_2^m\PAIR_{1,1}\EVEN_2^n\}.$$
        The desired result then follows as $\EVEN_2$ and $\PAIR_{1,1}$ have exactly one transposition pair, and $\ODD_4$ and $\ODD_{2,2}$ have exactly two transposition pairs respectively.\\

        So there are $g$ pairs $\{\alpha_1,\beta_1\}\dots\{\alpha_g,\beta_g\}\in\AA_2$. Now we show that the set of $\vv_{\alpha_i,\beta_i}$'s are linearly independent. Suppose that there exists $c_1,\dots,c_g\in\RR$ such that
        $$ \ww = c_1\vv_{\alpha_1,\beta_1} + \dots + c_g\vv_{\alpha_g,\beta_g} = 0.$$
        Suppose $i$ is such that $\{\alpha_i,\beta_i\}=\{a,b\}\subseteq\EVEN_2$ or $\PAIR_{1,1}$ (as in the notation mentioned in \ref{sec.stdspin}), where
        $$\EVEN_2=\cmtrx{\LL{a}{b}~\LL{b}{a}} \mbox{ and } \PAIR_{1,1}=\cmtrx{\LL{a}{b}~\LL{c}{c}~\LL{b}{a}}.$$
        We see that
        $$ \omega(\vv_{a},\ww) = c_i\Omega_{a,b} = \pm c_i = 0,$$
        implying that $c_i=0$. Suppose $i$ is such that $\{\alpha_i,\beta_i\}=\{a,b\}\subseteq\ODD_{2,2}$ or $\ODD_4$ for
        $$\ODD_{2,2}=\cmtrx{\LL{a}{b} ~\LL{c}{d} ~\LL{e}{e} ~\LL{d}{c} ~\LL{b}{a}}\mbox{ and }\ODD_4=\cmtrx{\LL{a}{b} ~\LL{c}{d} ~\LL{d}{c} ~\LL{b}{a}}$$
        with $j$ such that pair $\{\alpha_j,\beta_j\}=\{c,d\}\in\AA_2$. Then
        $$ \begin{array}{rcl}
            \omega(\vv_{c},\ww) & = & c_j\Omega_{c,d}\\
                & = & \pm c_j=0 \Rightarrow c_j=0\\
            \omega(\vv_{a},\ww) & = & c_j(\Omega_{a,c}+\Omega_{a,d})+ c_i\Omega_{a,b}\\
                & = & \pm c_i=0 \Rightarrow c_i=0.
        \end{array}$$
        We now see that for $i$ such that $\{\alpha_i,\beta_i\}=\{\lA,\lB\}\in\AA_2$ (the outside letters of the permutation),
        $$ \ww = c_i\vv_{\lA,\lB} = 0 \Rightarrow c_i = 0.$$
        So $c_1=\dots=c_g=0$, implying that the $\vv_{\alpha_i,\beta_i}$'s are linearly independent.
    \end{proof} 


\bibliographystyle{abbrv}
\bibliography{bibfile}

\end{document}